\documentclass{article}
\usepackage{amsmath, amsthm, amssymb, amsfonts}
\usepackage{geometry}
\usepackage{graphicx}
\usepackage[titletoc]{appendix}
\usepackage{natbib}

\newtheorem{theorem}{Theorem}[section]
\newtheorem{lemma}[theorem]{Lemma}
\newtheorem{proposition}[theorem]{Proposition}
\newtheorem{corollary}[theorem]{Corollary}

\theoremstyle{definition}

\newtheorem{example}{Example}[section]
\newtheorem{assumption}{Assumption}
\newtheorem{remark}[theorem]{Remark}

\newcommand{\secret}[1]{}
\input{tcilatex}

\begin{document}

\title{On the Power of Invariant Tests for Hypotheses on a Covariance Matrix}
\author{David Preinerstorfer and Benedikt M. P\"{o}tscher\thanks{%
We thank Federico Martellosio and two referees for helpful comments on an
earlier version of the paper. Address for correspondence: Department of
Statistics, University of Vienna, Oskar-Morgenstern-Platz 1, A-1090 Vienna,
Austria. E-mail: \{david.preinerstorfer, benedikt.poetscher\}@univie.ac.at}}
\date{First version: February 2013\\
Second version: March 2014\\
This version: December 2014}
\maketitle

\begin{abstract}
The behavior of the power function of autocorrelation tests such as the
Durbin-Watson test in time series regressions or the Cliff-Ord test in
spatial regression models has been intensively studied in the literature.
When the correlation becomes strong, \cite{kramer1985} (for the
Durbin-Watson test) and \cite{kramer2005} (for the Cliff-Ord test) have
shown that the power can be very low, in fact can converge to zero, under
certain circumstances. Motivated by these results, \cite{Mart10} set out to
build a general theory that would explain these findings. Unfortunately, 
\cite{Mart10} does not achieve this goal, as a substantial portion of his
results and proofs suffer from serious flaws. The present paper now builds a
theory as envisioned in \cite{Mart10} in a fairly general framework,
covering general invariant tests of a hypothesis on the disturbance
covariance matrix in a linear regression model. The general results are then
specialized to testing for spatial correlation and to autocorrelation
testing in time series regression models. We also characterize the situation
where the null and the alternative hypothesis are indistinguishable by
invariant tests.


\end{abstract}

\section{Introduction}

Testing hypotheses on the covariance matrix of the disturbances in a
regression model is an important problem in econometrics and statistics, a
prime example being testing the hypothesis of uncorrelatedness of the
disturbances. Two particularly important cases are (i) testing for
autocorrelation in time series regressions and (ii) testing for spatial
autocorrelation in spatial models; for an overview see \cite{King1987a} and 
\cite{Anselin}. For testing autocorrelation in time series regressions the
most popular test is probably the Durbin-Watson test. While low power of
this test against highly correlated alternatives in some instances had been
noted earlier by \cite{Tillman1975} and \cite{King1985}, \cite{kramer1985}
seems to have been the first to show that the limiting power of the
Durbin-Watson test as autocorrelation goes to one can actually be zero. This
phenomenon has become known as the \textit{zero-power trap}. The work by 
\cite{kramer1985} has been followed up and extended in the context of
testing against autoregressive disturbances of order one in \cite{zeisel1989}%
, \cite{KramerZeisel1990}, and \cite{lobus2000}; see also \cite{small1993}
and \cite{bartels1992}. Loosely speaking, these results show that the power
of the Durbin-Watson test (and of a class of related tests) typically
converges to either one or zero (depending on whether a certain observable
quantity is below or above a threshold) as the strength of autocorrelation
increases, provided that there is no intercept in the regression (in the
sense that the vector of ones is not in the span of the regressor matrix);
in case an intercept is in the regression, the limit is typically neither
zero nor one. Some of these results were extended in \cite{kleiber2005} to
the case where the Durbin-Watson test is used, but the disturbances are
fractionally integrated. In the context of spatial regression models \cite%
{kramer2005} showed that the Cliff-Ord test can similarly be affected by the
zero-power trap. \cite{Mart10} set out to build a general theory for power
properties of tests of a hypothesis on the covariance matrix of the
disturbances in a linear regression, that would also uncover the mechanism
responsible for the phenomena observed in the before-cited literature. While
the intuition behind the general results in \cite{Mart10} is often correct,
the results themselves and/or their proofs are not. For example, the main
result (Theorem 1 in \cite{Mart10}), on which much of that paper rests, has
some serious flaws: Parts of the theorem are incorrect, and the proofs of
the correct parts are substantially in error. In particular, the proof in 
\cite{Mart10} is based on a "concentration" effect, which, however, is
simply not present in the setting of the proof of Theorem 1 in \cite{Mart10}%
, as the relevant distributions "stretch out" rather than "concentrate".
This has already been observed in \cite{mynbaev2012}, where a way to
circumvent the problems was suggested. Mynbaev's approach, which is based on
the "stretch-out effect", is somewhat cumbersome in that it requires the
development of tools dealing with the "stretch-out effect"; furthermore, the
treatment in \cite{mynbaev2012} is given only for a subclass of the tests
considered in \cite{Mart10} and under more restrictive distributional
assumptions than in \cite{Mart10}.

In the present paper we now build a theory as envisioned in \cite{Mart10} at
an even more general level. In particular, we allow for general invariant
tests including randomized ones, we employ weaker conditions on the
underlying covariance model as well as on the distributions of the
disturbances (e.g., we even allow for distributions that are not absolutely
continuous). One aspect of our theory is to show how invariance of the tests
considered can be used to convert Martellosio's intuition about the
"concentration" effect into a precise mathematical argument. Furthermore,
advantages of this approach over the approach in \cite{mynbaev2012} are that
(i)\ standard weak convergence arguments can be used (avoiding the need for
new tools to handle the "stretch-out" effect), (ii) more general classes of
tests can be treated, and (iii) much weaker distributional assumptions are
required. The general theory built in this paper is then applied to tests
for spatial autocorrelation, which, in particular, leads to correct versions
of the results in \cite{Mart10} that pertain to spatial models.\footnote{%
This involves more than just providing a correct version of Theorem 1, the
main result in \cite{Mart10}, and is not undertaken in \cite{mynbaev2012},
see his Remark 2.12.} A further contribution of the present paper is a
characterization of the situation where no invariant test can distinguish
the null hypothesis of no correlation from the alternative. This
characterization helps to explain, and provides a unifying framework for,
phenomena observed in \cite{Kadiyala}, \cite{Arnold79}, \cite{Kariya80JASA}, 
\cite{Mart10}, and \cite{Mart11}.

The paper is organized as follows: After laying out the framework in Section %
\ref{Framework}, the general theory is developed in Section \ref{main}. The
main results are Theorems \ref{MT}, \ref{BT}, and \ref{BT2}. Theorem \ref{MT}%
, specialized to nonrandomized tests, shows that under appropriate
assumptions the power of an invariant test converges to $0$ or $1$ as the
"boundary" of the alternative is approached. The limit is $0$ or $1$
depending on whether a certain observable vector $e$ (the "concentration
direction" of the underlying covariance model) belongs to the complement of
the closure or to the interior of the rejection region of the test. This
result constitutes a generalization of the correct parts of Theorem 1 in 
\cite{Mart10} (the proofs of which in \cite{Mart10} are incorrect). Theorems %
\ref{BT} and \ref{BT2} deal with the case where the concentration direction $%
e$ belongs to the boundary of the rejection region, a case excluded from
Theorem \ref{MT}, thus providing correct versions of the incorrect part of
Theorem 1 in \cite{Mart10}. The general results obtained in Theorems \ref{MT}%
, \ref{BT}, and \ref{BT2} are then specialized in Section \ref{T_B} to the
important class of tests based on test statistics that are ratios of
quadratic forms. The relationship between test size and the zero-power trap
is discussed in Section \ref{alpha_star}, before indistinguishability of the
null and alternative hypothesis by invariant tests is characterized in
Section \ref{IP}. Extensions of the general theory are discussed in Section %
\ref{gen}; in particular, we discuss ways of relaxing the distributional
assumptions. Section \ref{spatial} is devoted to applying the general theory
to testing for spatial correlation, while Section \ref{AR1C} contains an
application to testing for autocorrelation in time series regression models.
Whereas the problems with Theorem 1 in \cite{Mart10} are discussed in
Section \ref{main} as well as in Appendix \ref{Newapp}, problems with a
number of other results in \cite{Mart10} are dealt with in Appendix \ref{A2}%
. Proofs can be found in Appendices \ref{app_proofs} and \ref{app_proofs2}.
Some auxiliary results are collected in Appendix \ref{auxil}.

\section{The behavior of the power function: general theory\label{general}}

\subsection{Framework\label{Framework}}

As in \cite{Mart10}, we consider the problem of testing a hypothesis on the
covariance matrix of the disturbance vector in a linear regression model.
Given parameters $\beta \in \mathbb{R}^{k}$, $0<\sigma <\infty $, and $\rho
\in \lbrack 0,a)$, where $a$ is some prespecified positive real number, the
model is%
\begin{equation}
\mathbf{y}=X\beta +\mathbf{u},  \label{linmod}
\end{equation}%
where $X\in \mathbb{R}^{n\times k}$ is a non-stochastic matrix of rank $k$
with $0\leq k<n$ and $n\geq 2$. [In case $k=0$ we identify $\mathbb{R}%
^{n\times k}$, the space of real $n\times k$ matrices, with $\left\{
0\right\} \subseteq \mathbb{R}^{n}$ and $\mathbb{R}^{k}$ with $\left\{
0\right\} \subseteq \mathbb{R}$.] The disturbance vector $\mathbf{u}$ is
assumed to be an $n\times 1$ random vector with mean zero and covariance
matrix $\sigma ^{2}\Sigma (\rho )$, where $\Sigma (.)$ is a \emph{known}
function from $[0,a)$ to the set of symmetric and positive definite $n\times
n$ matrices. Without loss of generality (w.l.o.g.) $\Sigma (0)$ is assumed
to be the identity matrix $I_{n}$. [The case $a=\infty $ can be immediately
reduced to the case $a<\infty $ considered here by use of a transformation
like $\arctan (\rho )$.] We assume furthermore that, given $\beta $, $\sigma 
$, and $\rho $, the distribution of $\mathbf{u}$ is completely specified
(but see Remark \ref{rem_gen_1.5} in Section \ref{gen} for a relaxation of
this assumption). Note that this does not imply in general that the
distribution of $\sigma ^{-1}\Sigma ^{-1/2}(\rho )\mathbf{u}$ is independent
of $\rho $, $\sigma $, and $\beta $ (although this will often be the case in
important examples). In contrast to \cite{Mart10} we do not impose any
further assumptions on the distribution of $\mathbf{u}$ at this stage (see
Remark \ref{modelM10} below for a discussion of the additional assumptions
in \cite{Mart10}). All additional distributional assumptions needed later
will be stated explicitly in the theorems.

Under the preceding assumptions, model (\ref{linmod}) induces a \emph{%
parametric} family of distributions 
\begin{equation}
\mathfrak{P}=\left\{ P_{\beta ,\sigma ,\rho }:\beta \in \mathbb{R}%
^{k},0<\sigma <\infty ,\rho \in \lbrack 0,a)\right\}  \label{family}
\end{equation}%
on the sample space $(\mathbb{R}^{n},\mathcal{B}(\mathbb{R}^{n}))$ where $%
P_{\beta ,\sigma ,\rho }$ stands for the distribution of $\mathbf{y}$ under
the given parameters $\beta $, $\sigma $, and $\rho $, and where $\mathcal{B}%
(\mathbb{R}^{n})$ denotes the Borel $\sigma $-field on $\mathbb{R}^{n}$. The
expectation operator with respect to (w.r.t.) $P_{\beta ,\sigma ,\rho }\in 
\mathfrak{P}$ shall be denoted by $E_{\beta ,\sigma ,\rho }$. If $M$ is a
Borel-measurable mapping from $\mathbb{R}^{n}$ to $\mathbb{R}^{m}$, we shall
denote by $P_{\beta ,\sigma ,\rho }\circ M$ the pushforward measure of $%
P_{\beta ,\sigma ,\rho }$ under $M$, which is defined on $(\mathbb{R}^{m},%
\mathcal{B}(\mathbb{R}^{m}))$. As usual, a Borel-set $A$ will be said to be
a $\mathfrak{P}$-null set if it is a null set relative to every element of $%
\mathfrak{P}$.

\begin{remark}
\label{modelM10} \emph{(Comments on assumptions in \cite{Mart10})} (i) In 
\cite{Mart10}, p.154, additional assumptions on the distribution of $\mathbf{%
u}$ are imposed: for example, it is assumed that $\mathbf{u}$ possesses a
density which is positive everywhere on $\mathbb{R}^{n}$, is larger at $0$
than anywhere else, and satisfies a continuity property (the meaning of
which is not completely transparent). These assumptions are in general
stronger than what is needed; for example, as we shall see, some of our
results even hold for discretely distributed errors.

(ii) In \cite{Mart10} it is furthermore implicitly assumed that for fixed $%
\rho $, the distribution of $\sigma ^{-1}\mathbf{u}$ (or, equivalently, the
distribution of $\sigma ^{-1}\Sigma ^{-1/2}(\rho )\mathbf{u}$) does not
depend on $\beta $ and $\sigma $. This becomes apparent on p.~156, where it
is claimed that the testing problem under consideration is invariant w.r.t.
the group $G_{X}$ (defined below) in the sense of \cite{LR05}. In fact, \cite%
{Mart10} appears to even assume implicitly that the distribution of $\sigma
^{-1}\Sigma ^{-1/2}(\rho )\mathbf{u}$ is independent of \emph{all} the
parameters $\beta $, $\sigma $, and $\rho $; cf., e.g., the first line in
the proof of Theorem 1 on p. 182 in \cite{Mart10}.
\end{remark}

We consider the problem of testing $\rho =0$ against $\rho >0$. More
precisely, the null hypothesis and the alternative hypothesis are given by 
\begin{equation}
H_{0}:\rho =0,\beta \in \mathbb{R}^{k},0<\sigma <\infty \text{ \ against \ }%
H_{1}:\rho >0,\beta \in \mathbb{R}^{k},0<\sigma <\infty ,
\label{testproblem}
\end{equation}%
with the implicit understanding that always $\rho \in \lbrack 0,a)$. We note
that typically one would impose an additional (identifiability) condition
such as, e.g., $\sigma ^{2}\Sigma (\rho )\neq \tau ^{2}\Sigma (0)$ for every 
$\rho >0$ and every $0<\sigma ,\tau <\infty $ in order to ensure that $H_{0}$
and $H_{1}$ are disjoint, and hence that the test problem is meaningful.%
\footnote{%
Of course, even if $\sigma ^{2}\Sigma (\rho )=\tau ^{2}\Sigma (0)$ holds for
some $\sigma >0$, $\tau >0$ and some $\rho >0$, there may still be
additional identifiying information present in the distributions that goes
beyond the information contained in first and second moments.} The results
on the power behavior as $\rho \rightarrow a$ in the present paper are valid
without any such explicit identifiability condition, but note that one of
the basic assumptions (Assumption \ref{ASC}) underlying most of the results
automatically implies that $\sigma ^{2}\Sigma (\rho )\neq \tau ^{2}\Sigma
(0) $ for every $0<\sigma ,\tau <\infty $ holds at least for $\rho >0$ in a
neighborhood of $a$.

A (randomized) test is a Borel-measurable function $\varphi $ from the
sample space $\mathbb{R}^{n}$ to $[0,1]$, and a non-randomized test is the
indicator function of a set $\Phi \in \mathcal{B}(\mathbb{R}^{n})$, the
rejection region. A test statistic is a Borel-measurable function $T:\mathbb{%
R}^{n}\rightarrow \mathbb{R}$ which, together with a critical value $\kappa
\in \mathbb{R}$, gives rise to a rejection region $\left\{ y\in \mathbb{R}%
^{n}:T(y)>\kappa \right\} $.\footnote{%
The case of a test statistic $S$ taking values in the extended real line can
be easily accomodated in our framework by passing from $S$ to a real-valued
test statistic such as, e.g., $T=\arctan (S)$.} Note that the tests
(rejection regions, test statistics, critical values) may depend on the
sample size $n$ as well as on the design matrix $X$, but typically we shall
not show this in the notation. Recall that the size of a test $\varphi $ is
given by $\sup_{\beta \in \mathbb{R}^{k}}\sup_{0<\sigma <\infty }E_{\beta
,\sigma ,0}(\varphi )$, i.e., is the supremal rejection probability under
the null.

We shall also use the following terminology and notation: Random vectors and
matrices will always be denoted by boldface letters. All matrices considered
will be real matrices. The transpose of a matrix $A$ is denoted by $%
A^{\prime }$. The space spanned by the columns of $A$ is denoted by $%
\limfunc{span}(A)$. Given a linear subspace $L$ of $\mathbb{R}^{n}$, the
symbol $\Pi _{L}$ denotes orthogonal projection onto $L$, and $L^{\bot }$
denotes the orthogonal complement of $L$. Given an $n\times m$ matrix $Z$ of
rank $m$ with $0\leq m<n$, we denote by $C_{Z}$ a matrix in $\mathbb{R}%
^{(n-m)\times n}$ such that $C_{Z}C_{Z}^{\prime }=I_{n-m}$ and $%
C_{Z}^{\prime }C_{Z}=\Pi _{\limfunc{span}(Z)^{\bot }}$ where $I_{r}$ denotes
the identity matrix of dimension $r$. It is easily seen that every matrix
whose rows form an orthonormal basis of $\limfunc{span}(Z)^{\bot }$
satisfies these two conditions and vice versa, and hence any two choices for 
$C_{Z}$ are related by premultiplication by an orthogonal matrix. Let $l$ be
a positive integer. If $A$ is an $l\times l$ matrix and $\lambda \in \mathbb{%
R}$ is an eigenvalue of $A$ we denote the corresponding eigenspace by $%
\limfunc{Eig}\left( A,\lambda \right) $. The eigenvalues of a symmetric
matrix $B\in \mathbb{R}^{l\times l}$ ordered from smallest to largest and
counted with their multiplicities are denoted by $\lambda _{1}(B),\ldots
,\lambda _{l}(B)$. If $B$ is a symmetric and nonnegative definite $l\times l$
matrix, every $l\times l$ matrix $A$ that satisfies $AA^{\prime }=B$ is
called a \textit{square root} of $B$; with $B^{1/2}$ we denote its unique
symmetric and nonnegative definite square root. Note that every square root
of $B$ is of the form $B^{1/2}U$ for some orthogonal matrix $U$. A vector $%
x\in \mathbb{R}^{l}$ is said to be normalized if $\Vert {x}\Vert =1$, where $%
\Vert {.}\Vert $ denotes Euclidean norm on $\mathbb{R}^{l}$. The operators $%
\limfunc{bd}$, $\limfunc{int}$, and $\limfunc{cl}$ shall denote the
boundary, interior, and closure of a subset of $\mathbb{R}^{l}$ w.r.t. the
Euclidean topology. For $x\in \mathbb{R}^{l}$ the symbol $\delta _{x}$
denotes point mass at $x$. Lebesgue measure on $(\mathbb{R}^{l},\mathcal{B}(%
\mathbb{R}^{l}))$ shall be denoted by $\mu _{\mathbb{R}^{l}}$, while
Lebesgue measure on the Borel subsets of $(0,\infty )$ is denoted by $\mu
_{(0,\infty )}$. The uniform probability measure on the Borel subsets of $%
S^{n-1}$, the unit sphere in $\mathbb{R}^{n}$, is denoted by $\upsilon
_{S^{n-1}}$. We use $\Pr $ as a generic symbol for a probability measure,
with $E$ denoting the corresponding expectation operator.

\subsubsection{Groups of transformations, invariance, and maximal invariants 
\label{groups}}

Suppose that $G$ is a group of bijective Borel-measurable transformations $g:%
\mathbb{R}^{n}\rightarrow \mathbb{R}^{n}$, the group operation being
composition. A function $F$ defined on $\mathbb{R}^{n}$ is said to be
invariant w.r.t. $G$ if for every $y\in \mathbb{R}^{n}$ and every $g\in G$
we have $F(y)=F(g(y))$. A subset $A$ of $\mathbb{R}^{n}$ is said to be
invariant w.r.t. $G$ if for every $g\in G$ we have that $g(A)\subseteq A$.%
\footnote{%
The group structure implies that this is equivalent to $g(A)=A$ for every $%
g\in G$, and thus to invariance of the indicator function of $A$.} Of
course, invariance of $F$ implies invariance of $\left\{ y\in \mathbb{R}%
^{n}:F(y)>\kappa \right\} $.

Given a matrix $Z\in \mathbb{R}^{n\times m}$ such that $0\leq m<n$ with
column rank $m$, we will mainly work with the group%
\begin{equation*}
G_{Z}=\left\{ g_{\gamma ,\theta }:\gamma \in \mathbb{R}\backslash \left\{
0\right\} ,\theta \in \mathbb{R}^{m}\right\} ,
\end{equation*}%
where $g_{\gamma ,\theta }$ denotes the mapping $y\mapsto \gamma y+Z\theta $%
. The main reason for concentrating on invariance w.r.t. this group is that
the majority of tests for the hypothesis (\ref{testproblem}) considered in
the literature have this invariance property (for $Z=X$). Another reason is
that this is also the notion of invariance used in \cite{Mart10}.
Occasionally we shall consider invariance w.r.t. subgroups of $G_{Z}$, see
Remark \ref{Othgrps}.

The following is a maximal invariant w.r.t. $G_{Z}$%
\begin{equation*}
\mathcal{I}_{Z}(y)=%
\begin{cases}
\left\langle \Pi _{\limfunc{span}(Z)^{\bot }}y/\Vert {\Pi _{\limfunc{span}%
(Z)^{\bot }}y}\Vert \right\rangle & \text{ if }y\notin \limfunc{span}(Z), \\ 
0 & \text{ else},%
\end{cases}%
\end{equation*}%
where the function $\left\langle .\right\rangle :\mathbb{R}^{n}\rightarrow 
\mathbb{R}^{n}$ is defined as follows: $\left\langle y\right\rangle $ equals 
$y$ multiplied by the sign of the first nonzero coordinate of $y$ whenever $%
y\neq 0$, and $\left\langle y\right\rangle =0$ if $y=0$ (see \cite{PP13},
Section 5.1, where the group $G_{Z}$ is denoted as $G\left( \limfunc{span}%
(Z)\right) $). More generally, let $\zeta $ be any function from the unit
sphere in $\mathbb{R}^{n}$ into itself that satisfies $\zeta (y)=\zeta (-y)$
and has the property that $\zeta (y)$ is collinear with $y$; then defining $%
\mathcal{I}_{Z,\zeta }(y)$ in the same way as $\mathcal{I}_{Z}(y)$, but with 
$\left\langle .\right\rangle $ replaced by $\zeta $, provides another
maximal invariant w.r.t. $G_{Z}$. Obviously, given any normalized vector $e$%
, we can find a $\zeta $ as above that is additionally Borel-measurable and
is continuous in a neighborhood (in the unit sphere) of $e$. For such a $%
\zeta $ the maximal invariant $\mathcal{I}_{0,\zeta }$ w.r.t. $G_{0}$ is
then continuous in a neighborhood of $e$ (in $\mathbb{R}^{n}$), and hence in
a neighborhood of $\lambda e$ for any $\lambda \neq 0$ (in $\mathbb{R}^{n}$%
). Moreover, we can even choose $\zeta $ to be as before and also to satisfy 
$\zeta \left( e\right) =e$.\footnote{%
In fact, $\zeta $ then coincides with the identity in a neighborhood (in the
unit sphere) of $e$.} In the following we shall write $\zeta _{e}$ for any
such $\zeta $.\footnote{%
On p.~156 of \cite{Mart10} it is claimed that the quantity $\nu $ defined
there is a maximal invariant for the group $G_{X}$ (denoted by $F_{X}$ in 
\cite{Mart10}). First note that the author does not spell out how $\nu $ is
defined for $y\in \limfunc{span}(Z)$ and how $\limfunc{sgn}(0)$ is to be
interpreted. Second, regardless of how one defines $\nu $ on $\limfunc{span}%
(Z)$ and whether one interpretes $\limfunc{sgn}(0)$ as $0$, $1$, or $-1$,
the quantity $\nu $ is not invariant in general as can be seen from simple
examples.}

\begin{remark}
\label{MI} (i) For any test $\varphi $ invariant w.r.t. $G_{Z}$ we have $%
\varphi \left( y\right) =\varphi \left( \mathcal{I}_{Z}(y)\right) =\varphi
\left( \mathcal{I}_{Z,\zeta }(y)\right) $ for every $y\in \mathbb{R}^{n}$
and $\zeta $ as above. This is trivial for $y\in \limfunc{span}(Z)$ since $%
\varphi \left( y\right) =\varphi \left( 0\right) $ must hold by invariance.
For $y\notin \limfunc{span}(Z)$ observe that due to invariance we have 
\begin{equation*}
\varphi \left( y\right) =\varphi \left( \Pi _{\limfunc{span}(Z)^{\bot
}}y\right) =\varphi \left( \Pi _{\limfunc{span}(Z)^{\bot }}y/\Vert {\Pi _{%
\limfunc{span}(Z)^{\bot }}y}\Vert \right) =\varphi \left( \mathcal{I}%
_{Z}(y)\right) =\varphi \left( \mathcal{I}_{Z,\zeta }(y)\right) ,
\end{equation*}%
noting that $\mathcal{I}_{Z}(y)$ as well as $\mathcal{I}_{Z,\zeta }(y)$ are
proportional to $\Pi _{\limfunc{span}(Z)^{\bot }}y/\Vert {\Pi _{\limfunc{span%
}(Z)^{\bot }}y}\Vert $ with a proportionality factor equal to $\pm 1$.

(ii) For later use we note the following: if $\varphi $ is invariant w.r.t. $%
G_{Z}$, it is also invariant w.r.t. $G_{0}$. Consequently, we have $\varphi
\left( y\right) =\varphi \left( \mathcal{I}_{0}(y)\right) =\varphi \left( 
\mathcal{I}_{0,\zeta }(y)\right) $ for every $y\in \mathbb{R}^{n}$ and $%
\zeta $ as above.
\end{remark}

\begin{remark}
\label{invariance} If one assumes that the distribution of $\sigma ^{-1}%
\mathbf{u}$ does not depend on $\beta $ and $\sigma $ (as is, e.g., done in 
\cite{Mart10}, cf. Remark \ref{modelM10}(ii) above), the power function of
any $G_{X}$-invariant test $\varphi $ is then independent of $\beta $ and $%
\sigma $; that is, for every $\rho \in \lbrack 0,a)$ we have 
\begin{equation*}
E_{\beta ,\sigma ,\rho }(\varphi )=E_{0,1,\rho }(\varphi )\text{ for every }%
\beta \in \mathbb{R}^{k},0<\sigma <\infty .
\end{equation*}%
If, additionally, all the parameters of the model are identifiable, the test
problem (\ref{testproblem}) is then in fact a $G_{X}$-invariant test problem
in the sense of \cite{LR05}, Chapter 6.
\end{remark}

\begin{remark}
\label{Othgrps}In Sections \ref{IP} and \ref{idsp} as well as in Remark \ref%
{rem_gen_2} we shall also consider invariance w.r.t. the subgroups $%
G_{Z}^{+}=\left\{ g_{\gamma ,\theta }:\gamma >0,\theta \in \mathbb{R}%
^{m}\right\} $ and $G_{Z}^{1}=\left\{ g_{1,\theta }:\theta \in \mathbb{R}%
^{m}\right\} $ with associated maximal invariants 
\begin{equation*}
\mathcal{I}_{Z}^{+}(y)=%
\begin{cases}
\Pi _{\limfunc{span}(Z)^{\bot }}y/\Vert {\Pi _{\limfunc{span}(Z)^{\bot }}y}%
\Vert & \text{ if }y\notin \limfunc{span}(Z), \\ 
0 & \text{ else,}%
\end{cases}%
\end{equation*}%
and $\mathcal{I}_{Z}^{1}(y)=\Pi _{\limfunc{span}(Z)^{\bot }}y$, respectively.
\end{remark}

\subsection{Main results\label{main}}

We now set out to study the behavior of the power function of invariant
tests for the testing problem (\ref{testproblem}) when the parameter $\rho $
is `far away' from $0$, the value of $\rho $ under the null hypothesis,
i.e., when $\rho $ is close to its upper limit $a$. In particular, we are
interested in the \textit{limiting power} of such tests $\varphi $ as $\rho
\rightarrow a$, i.e., in $\lim_{\rho \rightarrow a}E_{\beta ,\sigma ,\rho
}(\varphi )$. For these limits as well as for all other limits where $\rho
\rightarrow a$ it is always implicitly understood that $\rho \in \lbrack
0,a) $, i.e., that one is considering left-hand side limits. [To avoid
confusion, we stress that throughout we consider a finite-sample situation,
i.e., sample size $n$ is fixed, and hence the notion of limiting power just
introduced has nothing to do with asymptotic power properties where sample
size increases to infinity.] To motivate our interest in this problem we
consider the following two examples.

\begin{example}
\emph{(Testing for positive autocorrelation)} Assume that the disturbances
in the regression model (\ref{linmod}) follow a Gaussian stationary
autoregressive process of order one with autoregressive parameter $\rho $.
Then the $(i,j)$-th element of $\Sigma (\rho )$ is given by $\left( 1-\rho
^{2}\right) ^{-1}\rho ^{\left\vert i-j\right\vert }$ for $\rho \in \lbrack
0,1)$, i.e., $a=1$. Unguided intuition may suggest that the power of
standard tests like the Durbin-Watson test for testing $\rho =0$ versus $%
\rho >0$ is large if $\rho $ is sufficiently different from zero, and, in
particular, if $\rho $ is close to $a=1$. In fact, this intuition may even
suggest that the power of the Durbin-Watson test should approach $1$ as $%
\rho \rightarrow a=1$. However, as already mentioned in the introduction,
this intuition is wrong: The limiting power of the Durbin-Watson test can be
zero (or one, or a number in $(0,1)$) depending on the design matrix and the
significance level employed (see \cite{kramer1985}, \cite{zeisel1989}, \cite%
{KramerZeisel1990}, and \cite{lobus2000}). $\square $
\end{example}

\begin{example}
\emph{(Testing for spatial autocorrelation)} Assume now that the
disturbances in the regression model (\ref{linmod}) are Gaussian spatial
autoregressive errors of order one. Then under typical assumptions on the
spatial weights matrix $W$ we have%
\begin{equation*}
\Sigma (\rho )=\left( I_{n}-\rho W\right) ^{-1}\left( I_{n}-\rho W^{\prime
}\right) ^{-1}
\end{equation*}%
for $\rho \in \lbrack 0,\lambda _{\max }^{-1})$, i.e., $a=\lambda _{\max
}^{-1}$. Here $\lambda _{\max }>0$ is a dominant eigenvalue of $W$. As in
the preceding example, unguided intuition may suggest that the limiting
power of standard tests like the Cliff-Ord test for $\rho \rightarrow
a=\lambda _{\max }^{-1}$ is large (e.g., is equal to $1$). However, this
intuition is again incorrect and the limiting power of the Cliff-Ord test
can be zero (or one, or a number in $(0,1)$) depending on the design matrix,
the weights matrix, and the significance level employed (see \cite%
{kramer2005}). $\square $
\end{example}

Our goal is now to develop a coherent theory for deriving the limiting power
of invariant tests for the testing problem (\ref{testproblem}), which allows
for more general correlation structures than the ones figuring in the
preceding examples and which allows for non-Gaussian distributions. As
mentioned in the introduction, an attempt at such a theory has been made in 
\cite{Mart10} and it is thus appropriate as a starting point to revisit and
discuss the main result in that paper: A large part of \cite{Mart10} is
devoted to determining the limiting power of non-randomized tests $\mathbf{1}%
_{\Phi }$ as $\rho \rightarrow a$, i.e., $\lim_{\rho \rightarrow a}P_{\beta
,\sigma ,\rho }(\Phi )$. Not surprisingly, the limiting behavior of the
power function crucially depends on the behavior of the function $\Sigma $
close to $a$. \cite{Mart10} concentrates on situations where $\Sigma
^{-1}(a-):=\lim_{\rho \rightarrow a}\Sigma ^{-1}(\rho )$ exists in $\mathbb{R%
}^{n\times n}$, and, in particular, on the case where the rank of $\Sigma
^{-1}(a-)$ equals $n-1$.\footnote{%
The case where $\Sigma ^{-1}(a-)$ exists and is positive definite
(equivalently, where the left-hand side limit $\Sigma (a-)$ of $\Sigma
(\cdot )$ exists and is positive definite) is not the focus of \cite{Mart10}
since then the model is typically also well-defined for $\rho =a$ and the
limiting power for $\rho \rightarrow a$ typically coincides with the power
for $\rho =a$.}$^{\text{,}}$\footnote{%
In \cite{Mart10}, p.~159, it is claimed that the following three cases are
exhaustive: (i) $\lim_{\rho \rightarrow a}\Sigma (\rho )$ exists and is
positive definite; (ii) $\lim_{\rho \rightarrow a}\Sigma (\rho )$ exists and
is singular and (iii) $\lim_{\rho \rightarrow a}\Sigma ^{-1}(\rho )$ exists
and is singular. This does not provide an exhaustive description of possible
cases, as there exist functions $\Sigma $ such that neither $\lim_{\rho
\rightarrow a}\Sigma (\rho )$ nor $\lim_{\rho \rightarrow a}\Sigma
^{-1}(\rho )$ exist. Let $n=2$ and define $\Sigma (\rho )$ as a diagonal
matrix with diagonal $(1-\rho ,(1-\rho )^{-1})$ for $\rho \in \lbrack 0,1)$.
Clearly, both $\Sigma (\rho )$ and its inverse do not converge as $\rho
\rightarrow 1$.} It should be observed that this condition on the function $%
\Sigma $ is satisfied in the two examples discussed above. In the following
we quote Theorem 1, the main result of \cite{Mart10}, which is set in the
framework described in Section \ref{Framework} augmented by the additional
distributional assumptions of \cite{Mart10}, discussed above in Remark \ref%
{modelM10}:

\begin{quote}
"Consider an invariant critical region $\Phi $ for testing $\rho =0$ against 
$\rho >0$ in model (1). Assume that $\Sigma (\rho )$ is positive definite as 
$\rho \rightarrow a$ \footnote{%
Positive definiteness is always assumed in \cite{Mart10} for $\rho \in
\lbrack 0,a)$, hence this assumption seems to be superfluous.}, and that $%
\limfunc{rank}(\Sigma ^{-1}(a))=n-1$. The limiting power of $\Phi $ as $\rho
\rightarrow a$ is:
\end{quote}

\begin{itemize}
\item 1 if $f_{1}(\Sigma ^{-1}(a))\in \limfunc{int}(\Phi )$;

\item in $(0,1)$ if $f_{1}(\Sigma ^{-1}(a))\in \limfunc{bd}(\Phi )$; or

\item 0 if $f_{1}(\Sigma ^{-1}(a))\notin \limfunc{cl}(\Phi )$."
\end{itemize}

From now on we shall refer to this theorem of \cite{Mart10} as MT1. A few
comments are in order: First, the notion of invariance used in the quote is
invariance w.r.t.$~G_{X}$. Second, observe that even if $\Sigma \left( \rho
\right) $ is well-defined for $\rho =a$ (which need not be the case in
general), the statement $\limfunc{rank}(\Sigma ^{-1}(a))=n-1$ as given in
the formulation of MT1 can obviously never be satisfied. To give meaning to
the above quote, the symbol $\Sigma ^{-1}(a)$ needs to be interpreted as $%
\Sigma ^{-1}(a-)$ throughout; this also becomes transparent from the proof
in \cite{Mart10}. Third, the symbol $f_{1}(A)$ in the above quote denotes a
normalized eigenvector of a symmetric matrix $A$ pertinent to $\lambda
_{1}(A)$, the smallest eigenvalue of $A$. Note that $\lambda _{1}(\Sigma
^{-1}(a-))=0$ due to the rank assumption in the quote. Furthermore, by the
rank assumption $f_{1}(\Sigma ^{-1}(a-))$ is uniquely determined up to a
sign change; because $\Phi $ is $G_{X}$-invariant, the validity of
conditions like $f_{1}(\Sigma ^{-1}(a-))\in \limfunc{int}(\Phi )$ therefore
does not depend on the choice of sign. Fourth, if $\Phi $ or its complement
is a (non-empty) $\mu _{\mathbb{R}^{n}}$-null set, then the second claim of
MT1 can obviously not hold. While these cases are unfortunately not ruled
out explicitly in the statement of MT1 (which may lead to confusion among
some readers), it should be noted that such cases are implicitly excluded in 
\cite{Mart10}, as the author considers only $G_{X}$-invariant rejection
regions $\Phi $ that have size strictly between zero and one, cf. \cite%
{Mart10}, p. 157. [Note that under the distributional assumptions in \cite%
{Mart10}, cf. Remark \ref{modelM10} above, $G_{X}$-invariance of $\Phi $
implies that the size of $\Phi $ is given by $P_{0,1,0}\left( \Phi \right) $
and that this is $0$ (or $1$) precisely if $\Phi $ (or its complement) is a $%
\mu _{\mathbb{R}^{n}}$-null set.]

Even with the just discussed appropriate interpretations, the second claim
in MT1 is incorrect (cf. also \cite{mynbaev2012}), and the proofs of the
correct parts (i.e., claims 1 and 3) are in error. Counterexamples to the
second claim in MT1 are provided in Examples \ref{EX1} and \ref{EX2} in
Appendix \ref{Newapp}. A discussion of the mistakes in the proof of the
correct parts of MT1 is also given in Appendix \ref{Newapp}. The following
section provides a generalization of the (correct) claims 1 and 3 in MT1,
whereas correct versions of the (incorrect) second claim in MT1 are provided
in Section \ref{2claim}.

\subsubsection{A generalization of the first and third claim in Theorem 1 in 
\protect\cite{Mart10}\label{1&3claim}}

The proof of MT1 given in \cite{Mart10} rests on a "concentration effect" to
occur in the distributions $P_{\beta ,\sigma ,\rho }$ as $\rho \rightarrow a$%
, namely that these distributions (for fixed $\beta $ and $\sigma $)
converge (in an appropriate sense) as $\rho \rightarrow a$ to a distribution
concentrated on a one-dimensional subspace. However, as discussed in
Appendix \ref{Newapp}, this concentration effect simply does not occur in
the way as claimed in \cite{Mart10} (cf. also \cite{mynbaev2012}). In fact,
the direct opposite happens: the distributions $P_{\beta ,\sigma ,\rho }$
stretch out, i.e., all of the mass "escapes to infinity". As we shall now
show, the problem can, however, be fixed: The crucial observation is that,
while rescaling the data has no effect on the rejection probability of $%
G_{X} $-invariant tests, an appropriate rescaling can enforce the desired
concentration effect. Formalizing this observation will lead us to Theorem %
\ref{MT}, which provides a generalization of the first and third claim of
MT1 under even weaker distributional assumptions than the ones used in MT1;
in addition, this theorem will also cover randomized tests. For a discussion
and some intuition regarding the concentration effect in a different setting
see \cite{PP13}, Section 5.2.\footnote{%
In the setting of \cite{PP13} no rescaling is needed to achieve the
concentration effect.} We shall make use of the following assumption on the
function $\Sigma $ which is weaker than the rank assumption in MT1.

\begin{assumption}
\label{ASC} $\lambda _{n}^{-1}(\Sigma (\rho ))\Sigma (\rho )\rightarrow
ee^{\prime }$ as $\rho \rightarrow a$ for some $e\in \mathbb{R}^{n}$.
\end{assumption}

Note that the vector $e$ in Assumption \ref{ASC} is necessarily normalized
and will be called \emph{concentration direction} of the underlying
covariance model. That this assumption is indeed weaker than the assumption
of a one-dimensional kernel of $\Sigma ^{-1}(a-)$ made in MT1 is shown in
the following lemma.\footnote{%
The proof idea is also used in the proof of Lemma E.4 in \cite{Mart10} in
the special case of a SAR(1) model. See also Lemma 3.3 in \cite{Mart11SPL}
and its proof.} Recall that when writing $\Sigma ^{-1}(a-)$ we always
implicitly assume that this limit exists in $\mathbb{R}^{n\times n}$.

\begin{lemma}
\label{convCM} If the normalized vector $e$ spans the kernel of $\Sigma
^{-1}(a-)$, then $\lambda _{n}^{-1}(\Sigma (\rho ))\Sigma (\rho )\rightarrow
ee^{\prime }$ as $\rho \rightarrow a$.
\end{lemma}

The converse is not true as shown in the subsequent example. This shows that
Assumption \ref{ASC} underlying Theorem \ref{MT} given below is strictly
weaker than the assumption of a one-dimensional kernel of $\Sigma ^{-1}(a-)$
underlying MT1.

\begin{example}
\label{stretch} For $\rho \in \lbrack 0,1)$ let $\Sigma (\rho )$ be a $%
2\times 2$ diagonal matrix with diagonal entries $1$ and $1-\rho $. Then the
largest eigenvalue of $\Sigma (\rho )$ equals one and $\lambda
_{n}^{-1}(\Sigma (\rho ))\Sigma (\rho )$ converges to $ee^{\prime }$ as $%
\rho \rightarrow 1$, where $e=(1,0)^{\prime }$. But the limit of $\Sigma
^{-1}\left( \rho \right) $ for $\rho \rightarrow 1$ does obviously not
exist. Another example, where $\lambda _{n}^{-1}(\Sigma (\rho ))\Sigma (\rho
)\rightarrow ee^{\prime }$ and for which the limit of $\Sigma ^{-1}\left(
\rho \right) $ for $\rho \rightarrow 1$ exists, but does not have a
one-dimensional kernel, is provided by the $2\times 2$ diagonal matrix with
diagonal entries $\left( 1-\rho \right) ^{-1}$ and $\left( 1-\rho \right)
^{-1/2}$. In this case the limit $\Sigma ^{-1}\left( 1-\right) $ exists and
equals the zero matrix. $\square $
\end{example}

For $\xi \in \mathbb{R}^{n}$ and $\delta \in \mathbb{R}\backslash \left\{
0\right\} $ let $M_{\xi ,\delta }$ denote the mapping $y\mapsto \delta
^{-1}(y-\xi )$ from $\mathbb{R}^{n}$ to $\mathbb{R}^{n}$. We now introduce
the following high-level assumption on $\mathfrak{P}$ which will be seen to
be satisfied under the assumptions in \cite{Mart10} underlying MT1. Simple
sufficient conditions for this assumption that are frequently satisfied are
discussed below.

\begin{assumption}
\label{ASD} For every $\beta \in \mathbb{R}^{k}$, $0<\sigma <\infty $, and
every sequence $\rho _{m}\in \lbrack 0,a)$ converging to $a$, every weak
accumulation point $P$ of 
\begin{equation}
P_{\beta ,\sigma ,\rho _{m}}\circ M_{X\beta ,\lambda _{n}^{1/2}(\Sigma (\rho
_{m}))\sigma }  \label{pushforward}
\end{equation}%
satisfies $P(\left\{ 0\right\} )=0$.
\end{assumption}

The measure in (\ref{pushforward}) will in general \emph{not} coincide with $%
P_{0,\lambda _{n}^{-1/2}(\Sigma (\rho _{m})),\rho _{m}}$. However, in the
important special case, where the distribution of $\sigma ^{-1}\mathbf{u}$
does not depend on $\beta $ and $\sigma $ (cf. Remark \ref{invariance}),
these two measures will indeed coincide. We furthermore note that in view of
Lemma \ref{convPM} in Appendix \ref{app_proofs} the sequence in (\ref%
{pushforward}) is automatically tight whenever Assumption \ref{ASC} is
satisfied.

Before we present our generalizations of the first and third claim in MT1 we
provide simple sufficient conditions for the high-level Assumption \ref{ASD}%
. To this end we introduce the following assumption on $\mathfrak{P}$ that
is clearly satisfied in many examples.

\begin{assumption}
\label{ASDR} There exists an $n\times 1$ random vector $\mathbf{z}$ with
mean zero and covariance matrix $I_{n}$ such that for every $\beta \in 
\mathbb{R}^{k}$, every $0<\sigma <\infty $, and every $\rho \in \lbrack 0,a)$
the distribution $P_{\beta ,\sigma ,\rho }$ is induced by model (\ref{linmod}%
) with $\mathbf{u}$ having the same distribution as $\sigma L(\rho )\mathbf{z%
}$ and where the matrices $L(\rho )$ satisfy $L(\rho )L^{\prime }(\rho
)=\Sigma (\rho )$.\footnote{%
Note, in particular, that the distribution of $\mathbf{z}$ is independent of 
$\beta $, $\sigma $, and $\rho $.}
\end{assumption}

Important examples of families $\mathfrak{P}$ satisfying Assumption \ref%
{ASDR} are provided by elliptically symmetric families. Here $\mathfrak{P}$
is said to be an \emph{elliptically symmetric family} if it satisfies
Assumption \ref{ASDR} and $\mathbf{z}$ is spherically symmetric, i.e., the
distributions of $U\mathbf{z}$ and $\mathbf{z}$ are the same for every
orthogonal matrix $U$.\footnote{%
The notion of an elliptically symmetric family implies elliptical symmetry
of its elements, but is stronger (as the distribution of $\mathbf{z}$ in
Assumption \ref{ASDR} is not allowed to vary with the parameters).}
Obviously, if $\mathfrak{P}$ is an elliptically symmetric family, we can
assume without loss of generality that $L(\rho )=\Sigma ^{1/2}(\rho )$ in
Assumption \ref{ASDR} (because any $L(\rho )$ satisfies $L(\rho )=\Sigma
^{1/2}(\rho )U(\rho )$ for some orthogonal matrix $U(\rho )$). Furthermore,
recall from Remark \ref{modelM10} that \cite{Mart10} implicitly imposes
Assumption \ref{ASDR} (with $L(\rho )=\Sigma ^{1/2}(\rho )$) and more.
Sufficient conditions for Assumption \ref{ASD} are now as follows.

\begin{proposition}
\footnote{%
Inspection of the proof shows that, more generally, Assumptions \ref{ASC}
and \ref{ASDR} imply Assumption \ref{ASD} as soon as $\Pr (e^{\prime }U%
\mathbf{z}=0)=0$ holds for any orthogonal matrix that arises as an
accumulation point of \ $\Sigma ^{-1/2}(\rho )L(\rho )$ for $\rho
\rightarrow a$.}\label{AsPexII} Suppose Assumptions \ref{ASC} and \ref{ASDR}
are satisfied.

\begin{enumerate}
\item If $L(\rho )=\Sigma ^{1/2}(\rho )$ and $\Pr (e^{\prime }\mathbf{z}%
=0)=0 $ hold for $e$ as in Assumption \ref{ASC} and for $L(\cdot )$ and $%
\mathbf{z} $ as in Assumption \ref{ASDR}, then $\mathfrak{P}$ satisfies
Assumption \ref{ASD}.

\item If the distribution of $\mathbf{z}$ is absolutely continuous w.r.t. $%
\mu _{\mathbb{R}^{n}}$, then $\mathfrak{P}$ satisfies Assumption \ref{ASD}.
More generally, if $\Pr (\mathbf{z}=0)=0$ and the distribution of $\mathbf{z}%
/\left\Vert \mathbf{z}\right\Vert $ is absolutely continuous w.r.t. the
uniform distribution $\upsilon _{S^{n-1}}$ on the unit sphere $S^{n-1}$,
then $\mathfrak{P}$ satisfies Assumption \ref{ASD}.
\end{enumerate}
\end{proposition}

We note that Part 1 of the preceding proposition shows that Assumption \ref%
{ASD} also allows for families of discrete distributions. In some contexts
(e.g., spatial regression models) it is convenient to avoid the assumption $%
L(\rho )=\Sigma ^{1/2}(\rho )$ made in Part 1. Part 2 shows that this
assumption can indeed be avoided at the cost of introducing additional
conditions on the distribution of $\mathbf{z}$. That the assumptions for the
second statement in Part 2 are indeed weaker than the assumptions for the
first statement in Part 2 follows from Lemma \ref{Proj_1} in Appendix \ref%
{auxil}.

We are now ready to present and prove a generalization of the first and
third claim in MT1. The result is stated for possibly randomized tests.

\begin{theorem}
\label{MT} Suppose Assumptions \ref{ASC} and \ref{ASD} are satisfied and let 
$\varphi $ be a test that is invariant w.r.t. $G_{X}$ and is continuous at $%
e $, where $e$ is as in Assumption \ref{ASC}. Then for every $\beta \in 
\mathbb{R}^{k}$ and $0<\sigma <\infty $ we have that $E_{\beta ,\sigma ,\rho
}(\varphi )\rightarrow \varphi (e)$ as $\rho \rightarrow a$.
\end{theorem}

In the next remark we discuss why Theorem \ref{MT} contains the first and
the third claim of MT1 as special cases.

\begin{remark}
\label{RMT} (i) First observe that in light of Lemma \ref{convCM},
Proposition \ref{AsPexII}, and Remark \ref{modelM10} the assumptions of
Theorem \ref{MT} are weaker than the assumptions in MT1. Second, under the
assumptions of MT1 the vector $e$ coincides with $f_{1}(\Sigma ^{-1}(a-))$
in MT1 (possibly up to an irrelevant sign). Third, if $\varphi $ in Theorem %
\ref{MT} is specialized to the indicator function of a rejection region $%
\Phi $ that is invariant w.r.t. $G_{X}$, the above theorem reduces to:

$\bullet $ If $e\in \limfunc{int}(\Phi )$, then for every $\beta \in \mathbb{%
R}^{k}$ and $0<\sigma <\infty $ we have $\lim\limits_{\rho \rightarrow
a}P_{\beta ,\sigma ,\rho }(\Phi )=1$, and

$\bullet $ if $e\notin \limfunc{cl}(\Phi )$, then for every $\beta \in 
\mathbb{R}^{k}$ and $0<\sigma <\infty $ we have $\lim\limits_{\rho
\rightarrow a}P_{\beta ,\sigma ,\rho }(\Phi )=0$.

To see this simply observe that $\varphi =\mathbf{1}_{\Phi }$ is continuous
at $e$ if and only if $e\notin \limfunc{bd}\left( \Phi \right) $.

(ii) If $e\in \limfunc{bd}(\Phi )$, Theorem \ref{MT} is not applicable as it
stands because $\mathbf{1}_{\Phi }$ is then not continuous at $e$. However,
in some cases Theorem \ref{MT} can be used in an indirect way as follows:
suppose the rejection region $\Phi $ can be modified into an `equivalent'
rejection region $\Phi ^{\ast }$ (in the sense that $\Phi $ and $\Phi ^{\ast
}$ differ only by a $\mathfrak{P}$-null set) such that now $e\notin \limfunc{%
bd}\left( \Phi ^{\ast }\right) $ holds. As $\Phi $ and $\Phi ^{\ast }$ give
rise to the same rejection probabilities, we can therefore obtain the limits
of the rejection probabilities of $\Phi $ by applying Theorem \ref{MT} to $%
\Phi ^{\ast }$. More generally, suppose $\varphi $ is a test that is equal
to a test $\varphi ^{\ast }$ outside of a $\mathfrak{P}$-null set and
suppose that $\varphi ^{\ast }$ satisfies the assumptions of Theorem \ref{MT}%
. As $\varphi $ and $\varphi ^{\ast }$ have the same rejection
probabilities, we can conclude that $E_{\beta ,\sigma ,\rho }(\varphi
)\rightarrow \varphi ^{\ast }(e)$ as $\rho \rightarrow a$. [Of course, a
simple sufficient condition for $\mathfrak{P}$-almost everywhere equality of 
$\varphi =\varphi ^{\ast }$ is that $\mathfrak{P}$ is dominated by a measure 
$\nu $, say, and $\varphi =\varphi ^{\ast }$ holds $\nu $-almost everywhere.]
\end{remark}

\begin{remark}
\label{RASDR} Theorem \ref{MT} applies to $G_{X}$-invariant tests. Such
tests have a natural justification if the underlying test problem is
invariant under $G_{X}$ itself (which is not in general required in Theorem %
\ref{MT}). Recall from Remark \ref{invariance} that the test problem (\ref%
{testproblem}) is invariant under $G_{X}$ provided the distribution of $%
\sigma ^{-1}\mathbf{u}$ does not depend on $\beta $ and $\sigma $ (which is,
e.g., the case under Assumption \ref{ASDR}) and the parameters of the model
are identified.
\end{remark}

\subsubsection{Correct versions of the second claim in Theorem 1 in 
\protect\cite{Mart10} \label{2claim}}

As noted before, the second claim in MT1 is incorrect in general and
counterexamples to this claim are provided in Examples \ref{EX1} and \ref%
{EX2} in Appendix \ref{Newapp}. In this section we now aim at establishing
correct versions of this result under appropriate assumptions. Theorem \ref%
{BT} below will, in particular, provide an explicit expression for the
limiting power in case $e\in \limfunc{span}(X)$. Since $\limfunc{span}(X)$
turns out to always be a subset of the boundary of any critical region $\Phi 
$ ($\neq \emptyset ,\mathbb{R}^{n}$) that is invariant under $G_{X}$ (cf.
Proposition \ref{prop_bound} below), Theorem \ref{BT} can thus be seen as a
partial substitute for the second claim in MT1 (recall that under the
assumptions in \cite{Mart10} $e$ reduces to $f_{1}(\Sigma ^{-1}(a-))$).
Furthermore, in the important special case where the critical region is of
the form $\Phi =\left\{ y\in \mathbb{R}^{n}:T(y)>\kappa \right\} $, with $T$
invariant under $G_{X}$ and satisfying some regularity conditions, Theorem %
\ref{BT2} below will provide explicit expressions for the limiting power in
case $T\left( e\right) =\kappa $. For an important subclass of $G_{X}$%
-invariant test statistics $T$ (including certain ratios of quadratic forms
in $y$), Proposition \ref{prop_bound} below will show that $\limfunc{bd}%
(\Phi )=\limfunc{span}(X)\cup \left\{ y\in \mathbb{R}^{n}:T(y)=\kappa
\right\} $ holds (provided $\emptyset \neq \Phi \neq \mathbb{R}^{n}$).
Hence, for this subclass of tests, an application of Theorems \ref{BT} and %
\ref{BT2} together provides a substitute for the second claim in MT1
(because then $e\in \limfunc{bd}(\Phi )$ amounts to $e\in \limfunc{span}(X)$
or $T\left( e\right) =\kappa $).\footnote{%
For the discussion in this paragraph we have implicitly assumed that the
vector $e$ in Assumption \ref{ASC} and Assumption \ref{ASCII} is the same;
cf. Remark \ref{RASCII} below.} Before we can give these results we need to
study the structure of $\limfunc{bd}(\Phi )$ for $\Phi $ a $G_{X}$-invariant
rejection region.

\paragraph{On the structure of the boundary of $G_{X}$-invariant rejection
regions.}

\cite{Mart10}, Footnote 9, points out that a $G_{X}$-invariant rejection
region $\Phi $ always satisfies $\limfunc{span}(X)\subseteq \limfunc{bd}%
(\Phi )$ provided its size is neither zero nor one. Even if the rejection
region is assumed to be of the form $\Phi =\left\{ y\in \mathbb{R}%
^{n}:T(y)>\kappa \right\} $, then -- contrary to claims in \cite{Mart10} --
not much more can be said about the boundary $\limfunc{bd}(\Phi )$ in
general. This is discussed in the subsequent remark. In the proposition
following the remark we show how \cite{Mart10}'s claims, which are incorrect
in general, can be saved if additional assumptions are imposed on $T$.

\begin{remark}
\label{boundary} In \cite{Mart10}, p.~162 after Equation (9) and 2nd
paragraph on p. 167, it is incorrectly claimed (without providing an
argument) that for any critical region of the form $\left\{ y\in \mathbb{R}%
^{n}:T(y)>\kappa \right\} $, where $T$ is a $G_{X}$-invariant statistic, one
has 
\begin{equation}
\limfunc{bd}\left( \left\{ y\in \mathbb{R}^{n}:T(y)>\kappa \right\} \right) =%
\text{$\limfunc{span}$}(X)\cup \left\{ y\in \mathbb{R}^{n}:T(y)=\kappa
\right\} .\footnote{%
This equality is trivially violated if the critical region is empty or is
the entire space. However, such regions are ruled out in \cite{Mart10} as we
have already noted earlier.}  \label{incorr_claim}
\end{equation}%
To see that this claim is incorrect, consider the same setting as in Example %
\ref{EX2} in Appendix \ref{Newapp} and let $T=\mathbf{1}_{\Phi }$. Observe
that $\Phi $ can be written as $\left\{ y\in \mathbb{R}^{2}:T(y)>1/2\right\} 
$ and recall that $\Phi $ has rejection probability $1/2$ under the null.
Obviously, $\left\{ y\in \mathbb{R}^{2}:T(y)=1/2\right\} =\emptyset $ and $%
\limfunc{span}(X)=\left\{ 0\right\} $ hold, but 
\begin{equation*}
\limfunc{bd}\left( \left\{ y\in \mathbb{R}^{2}:T(y)>1/2\right\} \right)
=\left\{ y\in \mathbb{R}^{2}:y_{1}y_{2}=0\right\}
\end{equation*}%
which clearly is not equal to the set $\left\{ 0\right\} $.\footnote{%
Similar examples can be given when regressors are present and $n>k+1$.}
\end{remark}

Most rejection regions considered in practice (and in \cite{Mart10}, see,
e.g., p.~157) are of the form 
\begin{equation}
\left\{ y\in \mathbb{R}^{n}:y^{\prime }C_{X}^{\prime }BC_{X}y/\Vert {C_{X}y}%
\Vert ^{2}>\kappa \right\} ,  \label{quadratic_ill-defined}
\end{equation}%
where $B\in \mathbb{R}^{(n-k)\times (n-k)}$ is a given symmetric matrix,
which may depend on $X$ and/or the function $\Sigma $, and where ${C_{X}}$
satisfies ${C_{X}C}^{\prime }{_{X}=I}_{n-k}$ and ${C}^{\prime }{%
_{X}C_{X}=\Pi }_{\text{$\limfunc{span}$}(X)^{\bot }}$. First of all, this
rejection region is strictly speaking not well-defined, as the denominator
of the test statistic can take the value zero (namely, if and only if $y\in 
\limfunc{span}(X)$). However, whenever $\limfunc{span}(X)$ is a $\mathfrak{P}
$-null set, we then can pass to the well-defined rejection region%
\begin{equation}
\Phi _{B,\kappa }=\Phi _{B,C_{X},\kappa }=\left\{ y\in \mathbb{R}%
^{n}:T_{B}\left( y\right) >\kappa \right\} ,  \label{quadratic}
\end{equation}%
where%
\begin{equation}
T_{B}\left( y\right) =T_{B,C_{X}}\left( y\right) =\left\{ 
\begin{array}{cc}
y^{\prime }C_{X}^{\prime }BC_{X}y/\Vert {C_{X}y}\Vert ^{2} & \text{if }%
y\notin \limfunc{span}(X) \\ 
\lambda _{1}(B) & \text{if }y\in \limfunc{span}(X),%
\end{array}%
\right.  \label{T_quadratic}
\end{equation}%
without affecting the rejection probabilities. The condition that $\limfunc{%
span}(X)$ is a $\mathfrak{P}$-null set is certainly satisfied if (i) the
family $\mathfrak{P}$ is absolutely continuous w.r.t. Lebesgue measure $\mu
_{\mathbb{R}^{n}}$ (since $\limfunc{span}(X)$ is a $\mu _{\mathbb{R}^{n}}$%
-null set in view of our assumption $k<n$), or if (ii) $\mathfrak{P}$ is an
elliptically symmetric family with $\Pr \left( \mathbf{z}=0\right) =0$ where 
$\mathbf{z}$ is as in Assumption \ref{ASDR} (cf. Remark \ref{E1}(iii) in
Appendix \ref{auxil}). [Note that property (i) is always maintained in \cite%
{Mart10}.] We shall adopt the definitions in (\ref{quadratic}) and (\ref%
{T_quadratic}) regardless of whether or not $\limfunc{span}(X)$ is a $%
\mathfrak{P}$-null set. While assigning the value $\lambda _{1}(B)$ to $%
T_{B} $ on $\limfunc{span}(X)$ turns out to be convenient, it is of course
completely arbitrary. However, assigning to $T_{B}$ any other value on $%
\limfunc{span}(X)$ would, of course, have no effect on the rejection
probabilities provided $\limfunc{span}(X)$ is a $\mathfrak{P}$-null set, but
it could have an effect otherwise (in which case the original definition (%
\ref{quadratic_ill-defined}) does not lead to a test at all). At any rate,
an alternative assignment on $\limfunc{span}(X)$ has an easy to understand
effect on the rejection region itself and on its boundary, see Remark \ref%
{effect_on_boundary} below. The test statistic $T_{B}$ also depends on the
choice of ${C_{X}}$, a dependence which is typically suppressed in the
notation. Note that any other choice for $C_{X}$ is necessarily of the form $%
UC_{X}$ with $U$ an orthogonal matrix, and thus only has the simple effect
of "rotating" the matrix $B$ as $T_{B,C_{X}}=T_{UBU^{\prime },UC_{X}}$
holds. Clearly, $T_{B}$ is $G_{X}$-invariant.

Furthermore, observe that in case $\lambda _{1}(B)=\lambda _{n-k}(B)$ the
test statistic $T_{B}$ is constant equal to $\lambda _{1}(B)$, and hence the
resulting test is trivial in that the rejection region is either empty or
equal to the entire sample space (depending on the choice of $\kappa $).
While this case is trivial in the sense that the power properties of the
test are then obvious, it should be noted that this case may actually arise
for commonly used tests and for certain design matrices.

The third part of the subsequent proposition now shows that for rejection
regions of the form $\Phi _{B,\kappa }$ the claim (\ref{incorr_claim})
regarding the boundary is indeed correct (provided $\Phi _{B,\kappa }$ and
its complement are not empty). The first part of the proposition is just a
slight generalization of the observation in \cite{Mart10}, Footnote 9,
mentioned above. Regarding the second part we note that simple examples can
be given which show that in general the inclusion can be strict (even if $T$
is $G_{X}$-invariant).

\begin{proposition}
\label{prop_bound}

\begin{enumerate}
\item If $\Phi $ is a $G_{X}$-invariant rejection region satisfying $%
\emptyset \neq \Phi \neq \mathbb{R}^{n}$, then $\limfunc{span}(X)\subseteq 
\limfunc{bd}(\Phi )$.

\item If $T$ is a test statistic which is continuous on $\mathbb{R}%
^{n}\backslash \limfunc{span}$$(X)$, then 
\begin{eqnarray*}
\limfunc{bd}\left( \left\{ y\in \mathbb{R}^{n}:T(y)>\kappa \right\} \right)
&\subseteq &\text{$\limfunc{span}$}(X)\cup \left\{ y\in \mathbb{R}%
^{n}:T(y)=\kappa \right\} \\
&=&\text{$\limfunc{span}$}(X)\cup \left\{ y\in \mathbb{R}^{n}\backslash 
\text{$\limfunc{span}$}(X):T(y)=\kappa \right\} .
\end{eqnarray*}

\item If $\Phi _{B,\kappa }$ is as in (\ref{quadratic}), then 
\begin{eqnarray}
\limfunc{bd}(\Phi _{B,\kappa }) &=&\text{$\limfunc{span}$}(X)\cup \left\{
y\in \mathbb{R}^{n}:T_{B}\left( y\right) =\kappa \right\}  \notag \\
&=&\text{$\limfunc{span}$}(X)\cup \left\{ y\in \mathbb{R}^{n}\backslash 
\text{$\limfunc{span}$}(X):T_{B}\left( y\right) =\kappa \right\}
\label{char_bd}
\end{eqnarray}%
provided $\emptyset \neq \Phi _{B,\kappa }\neq \mathbb{R}^{n}$.
\end{enumerate}
\end{proposition}

\begin{remark}
\label{range_for_kappa}For $\kappa <\lambda _{1}(B)$ we have $\Phi
_{B,\kappa }=\mathbb{R}^{n}$, whereas for $\kappa \geq \lambda _{n-k}(B)$ we
have $\Phi _{B,\kappa }=\emptyset $. Hence, the non-trivial cases are when $%
\kappa $ belongs to the interval $[\lambda _{1}(B),\lambda _{n-k}(B))$ (and $%
\lambda _{1}(B)<\lambda _{n-k}(B)$ holds). Note that in case $\kappa
=\lambda _{1}(B)<\lambda _{n-k}(B)$ the rejection region is the complement
of a non-empty $\mu _{\mathbb{R}^{n}}$-null set (which automatically leads
to the rejection probabilities being identically equal to $1$ in case $%
\mathfrak{P}$ is dominated by $\mu _{\mathbb{R}^{n}}$, or $\mathfrak{P}$ is
an elliptically symmetric family with $\Pr \left( \mathbf{z}=0\right) =0$
where $\mathbf{z}$ is as in Assumption \ref{ASDR} (cf. Remark \ref{E1}(iii)
in Appendix \ref{auxil})), whereas for $\kappa \in \left( \lambda
_{1}(B),\lambda _{n-k}(B)\right) $ the rejection region as well as its
complement have positive $\mu _{\mathbb{R}^{n}}$-measure.
\end{remark}

\begin{remark}
\label{effect_on_boundary}As explained above assigning another value $c$,
say, to $T_{B}$ on $\limfunc{span}(X)$, resulting in a test statistic $%
T_{B}^{\prime }$, has no effect on the rejection probabilities provided $%
\limfunc{span}(X)$ is a $\mathfrak{P}$-null set. However, it can have an
effect on the resulting rejection region $\Phi _{B,\kappa }^{\prime }$, say,
and its boundary as follows: first, such a redefinition of $T_{B}$ on $%
\limfunc{span}(X)$ can obviously only add $\limfunc{span}(X)$ to $\Phi
_{B,\kappa }$ or remove it from $\Phi _{B,\kappa }$. Second, inspection of
the proof of Part 3 of Proposition \ref{prop_bound} shows that this result
continues to hold for $\Phi _{B,\kappa }^{\prime }$ provided $\emptyset \neq
\Phi _{B,\kappa }^{\prime }\neq \mathbb{R}^{n}$ and $\left\{ y\in \mathbb{R}%
^{n}\backslash \text{$\limfunc{span}$}(X):T_{B}\left( y\right) >\kappa
\right\} \neq \emptyset $. In case the latter set is empty, we necessarily
have $\Phi _{B,\kappa }^{\prime }=\emptyset $ or $\Phi _{B,\kappa }^{\prime
}=\limfunc{span}$$(X)$ (in which case Part 3 of Proposition \ref{prop_bound}
need not hold). But these are rather uninteresting cases as then the
rejection probability is always zero provided $\limfunc{span}(X)$ is a $%
\mathfrak{P}$-null set. [Also note that in these cases $\Phi _{B,\kappa
}=\emptyset $ always holds.] In particular, in the interesting case $\kappa
\in \lbrack \lambda _{1}(B),\lambda _{n-k}(B))$ with $\lambda
_{1}(B)<\lambda _{n-k}(B)$ we have $\Phi _{B,\kappa }^{\prime }=\Phi
_{B,\kappa }$ if $c\leq \kappa $ and $\Phi _{B,\kappa }^{\prime }=\Phi
_{B,\kappa }\cup \limfunc{span}(X)$ if $c>\kappa $; in both cases we have $%
\limfunc{bd}(\Phi _{B,\kappa }^{\prime })=\limfunc{bd}(\Phi _{B,\kappa })$
and (\ref{char_bd}) also holds for $T_{B}^{\prime }$.
\end{remark}

\paragraph{Correct versions of the second claim in MT1.}

We next provide an assumption on the function $\Sigma $ that will allow us
to establish results which, in particular, imply a version of the second
claim in MT1. The assumption may look somewhat intransparent at first sight.
However, it turns out to be satisfied for commonly used correlation
structures such as the ones generated by autoregressive models of order $1$
or spatial autoregressions, see Sections \ref{SEM} and \ref{AR1C}.

\begin{assumption}
\label{ASCII} There exists a function $c:[0,a)\rightarrow \left( 0,\infty
\right) $, a normalized vector $e\in \mathbb{R}^{n}$, and a square root $%
L_{\ast }(\cdot )$ of $\Sigma (\cdot )$ such that 
\begin{equation*}
\Lambda :=\lim_{\rho \rightarrow a}c(\rho )\Pi _{\text{$\limfunc{span}$}%
(e)^{\bot }}L_{\ast }(\rho )
\end{equation*}%
exists in $\mathbb{R}^{n\times n}$ and such that the linear map $\Lambda $
is injective when restricted to $\limfunc{span}(e)^{\bot }$.
\end{assumption}

We note that then the image of $\Lambda $ necessarily is $\limfunc{span}%
(e)^{\bot }$ and $\Lambda $ is a bijection from $\limfunc{span}(e)^{\bot }$
to itself. As we shall see in later sections, this assumption can be
verified for typical spatial models. For other types of models the
equivalent condition given in the subsequent lemma is easier to verify.

\begin{lemma}
\label{AR1} Let $c:[0,a)\rightarrow \left( 0,\infty \right) $ and a
normalized vector $e\in \mathbb{R}^{n}$ be given. Then the function $\Sigma
(\cdot )$ satisfies Assumption \ref{ASCII} for the given $c(\cdot )$, $e$,
and some square root $L_{\ast }(\cdot )$ of $\Sigma (\cdot )$ if and only if 
\begin{equation}
V:=\lim_{\rho \rightarrow a}c^{2}(\rho )\Pi _{\text{$\limfunc{span}$}%
(e)^{\bot }}\Sigma (\rho )\Pi _{\text{$\limfunc{span}$}(e)^{\bot }}
\label{scaled_limit_2}
\end{equation}%
exists in $\mathbb{R}^{n\times n}$ and the linear map $V$ is injective when
restricted to $\limfunc{span}(e)^{\bot }$. [Necessarily the image of $V$ is $%
\limfunc{span}(e)^{\bot }$ and $V$ is a bijection from $\limfunc{span}%
(e)^{\bot }$ to itself.]
\end{lemma}

\begin{remark}
\label{RASCII} Although Assumption \ref{ASCII} can hold independently of
Assumption \ref{ASC}, the relevant case for our theory is the case where $%
\Sigma $ satisfies both assumptions. If Assumptions \ref{ASC} and \ref{ASCII}
hold with $e$ and $e^{\ast }$, respectively, then we claim that $\limfunc{%
span}(e)=\limfunc{span}(e^{\ast })$ must hold whenever $n>2$. Since both
conditions only depend on the span of the respective vector, we can then
always choose $e^{\ast }=e$. To establish this claim write 
\begin{equation*}
c^{2}(\rho )\Pi _{\text{$\limfunc{span}$}(e^{\ast })^{\bot }}\Sigma (\rho
)\Pi _{\text{$\limfunc{span}$}(e^{\ast })^{\bot }}=c^{2}(\rho )\lambda
_{n}(\Sigma (\rho ))\Pi _{\text{$\limfunc{span}$}(e^{\ast })^{\bot }}\lambda
_{n}^{-1}(\Sigma (\rho ))\Sigma (\rho )\Pi _{\text{$\limfunc{span}$}(e^{\ast
})^{\bot }}
\end{equation*}%
and note that%
\begin{equation*}
\Pi _{\text{$\limfunc{span}$}(e^{\ast })^{\bot }}\lambda _{n}^{-1}(\Sigma
(\rho ))\Sigma (\rho )\Pi _{\text{$\limfunc{span}$}(e^{\ast })^{\bot
}}\rightarrow \Pi _{\text{$\limfunc{span}$}(e^{\ast })^{\bot }}ee^{\prime
}\Pi _{\text{$\limfunc{span}$}(e^{\ast })^{\bot }}
\end{equation*}%
as $\rho \rightarrow a$ by Assumption \ref{ASC}. Suppose $\limfunc{span}%
(e)\neq \limfunc{span}(e^{\ast })$ holds. We can then find $z\in \limfunc{%
span}$$(e^{\ast })^{\bot }$ with $z^{\prime }e\neq 0$. But then $z^{\prime
}\Pi _{\text{$\limfunc{span}$}(e^{\ast })^{\bot }}ee^{\prime }\Pi _{\text{$%
\limfunc{span}$}(e^{\ast })^{\bot }}z=\left( z^{\prime }e\right) ^{2}>0$
follows. Also note that $z^{\prime }Vz>0$ where $V$ is defined in Lemma \ref%
{AR1}. Together with the two preceding displays these observations imply
that $c^{2}(\rho )\lambda _{n}(\Sigma (\rho ))$ converges to a finite and
positive limit $b$, say. As a consequence, $V=b\Pi _{\text{$\limfunc{span}$}%
(e^{\ast })^{\bot }}ee^{\prime }\Pi _{\text{$\limfunc{span}$}(e^{\ast
})^{\bot }}$ must hold, i.e., $V$ would have to be a matrix of rank $1$.
However, $V$ is a matrix of rank $n-1$, a contradiction as $n>2$.
\end{remark}

The first result is now as follows. Note that under the assumptions of the
subsequent theorem the rejection probabilities actually do neither depend on 
$\beta $ nor $\sigma $, i.e., $E_{\beta ,\sigma ,\rho }(\varphi
)=E_{0,1,\rho }(\varphi )$, cf. Remark \ref{invariance}. For the sake of
readability the subsequent two theorems are not presented in their utmost
general form; possible extensions are discussed in Section \ref{gen}.

\begin{theorem}
\label{BT} Suppose Assumptions \ref{ASDR} and \ref{ASCII} hold. Let $\varphi 
$ be a test that is invariant w.r.t. $G_{X}$ and additionally satisfies the
invariance property 
\begin{equation}
\varphi (y)=\varphi (y+e)  \label{INV}
\end{equation}%
for every $y\in \mathbb{R}^{n}$ where $e$ is as in Assumption \ref{ASCII}.
Let $\mathcal{U}\left( L_{\ast }^{-1}L\right) $ denote the set of all
accumulation points of the orthogonal matrices $L_{\ast }^{-1}(\rho )L(\rho
) $ for $\rho \rightarrow a$, where $L(\rho )$ and $L_{\ast }(\rho )$ are
given in Assumptions \ref{ASDR} and \ref{ASCII}, respectively. Furthermore,
let $\beta \in \mathbb{R}^{k}$ and $0<\sigma <\infty $ be arbitrary but
given.

\begin{description}
\item[A.] Suppose the distribution of $\mathbf{z}$ (figuring in Assumption %
\ref{ASDR}) possesses a density $p$ w.r.t. Lebesgue measure $\mu _{\mathbb{R}%
^{n}}$ that is $\mu _{\mathbb{R}^{n}}$-almost everywhere continuous. Then:
\end{description}

\begin{enumerate}
\item Every accumulation point of $E_{\beta ,\sigma ,\rho }(\varphi )$ for $%
\rho \rightarrow a$ has the form $E_{Q_{\Lambda ,U}}\left( \varphi \right) $
with $U\in \mathcal{U}\left( L_{\ast }^{-1}L\right) $, where $Q_{\Lambda ,U}$
denotes the distribution of $\Lambda U\mathbf{z}$ and $\Lambda $ is given in
Assumption \ref{ASCII}. Conversely, every element $E_{Q_{\Lambda ,U}}\left(
\varphi \right) $ with $U\in \mathcal{U}\left( L_{\ast }^{-1}L\right) $ is
an accumulation point of $E_{\beta ,\sigma ,\rho }(\varphi )$ for $\rho
\rightarrow a$.

\item A sufficient condition for the set of accumulation points of $E_{\beta
,\sigma ,\rho }(\varphi )$ for $\rho \rightarrow a$ to be a singleton is
that $Q_{\Lambda ,U}$ is the same for all $U\in \mathcal{U}\left( L_{\ast
}^{-1}L\right) $ (which, e.g., is the case if $\mathcal{U}\left( L_{\ast
}^{-1}L\right) $ is a singleton). In this case $\lim_{\rho \rightarrow
a}E_{\beta ,\sigma ,\rho }(\varphi )$ exists and equals $E_{Q_{\Lambda
,U}}\left( \varphi \right) $.

\item Suppose the density $p$ is such that for $\upsilon _{S^{n-1}}$-almost
all elements $s\in S^{n-1}$ the function $p_{s}:(0,\infty )\rightarrow 
\mathbb{R}$ given by $p_{s}\left( r\right) =p\left( rs\right) $ does not
vanish $\mu _{(0,\infty )}$-almost everywhere. If $\varphi $ is neither $\mu
_{\mathbb{R}^{n}}$-almost everywhere equal to zero nor $\mu _{\mathbb{R}%
^{n}} $-almost everywhere equal to one, then the set of accumulation points,
i.e., $\left\{ E_{Q_{\Lambda ,U}}\left( \varphi \right) :U\in \mathcal{U}%
\left( L_{\ast }^{-1}L\right) \right\} $, is bounded away from zero and one.
\end{enumerate}

\begin{description}
\item[B.] Suppose $\mathfrak{P}$ is an elliptically symmetric family with
the distribution of $\mathbf{z}$ (figuring in Assumption \ref{ASDR})
satisfying $\Pr \left( \mathbf{z}=0\right) =0$. Then $E_{\beta ,\sigma ,\rho
}(\varphi )$ converges to $E_{Q_{\Lambda ,I_{n}}}\left( \varphi \right) $
for $\rho \rightarrow a$ and $E_{Q_{\Lambda ,I_{n}}}\left( \varphi \right) $
equals $E\left( \varphi \left( \Lambda \mathbf{G}\right) \right) $ where $%
\mathbf{G}$ is a multivariate Gaussian random vector with mean zero and
covariance matrix $I_{n}$. Furthermore, if $\varphi $ is neither $\mu _{%
\mathbb{R}^{n}}$-almost everywhere equal to zero nor $\mu _{\mathbb{R}^{n}}$%
-almost everywhere equal to one, then $0<E_{Q_{\Lambda ,I_{n}}}\left(
\varphi \right) <1$ holds.
\end{description}
\end{theorem}

\begin{remark}
\label{RBT}(i) The condition on the density $p$ in Part A.3 is quite weak.
It is, in particular, satisfied whenever $p$ is positive on an open
neighborhood of the origin except possibly for a $\mu _{\mathbb{R}^{n}}$%
-null set, but is much weaker. In fact, given the assumption that $p$
exists, the condition on the density $p$ in Part A.3 is equivalent to the
assumption that the density of $\mathbf{z}/\left\Vert \mathbf{z}\right\Vert $
w.r.t. the uniform distribution $\upsilon _{S^{n-1}}$ on the unit sphere is $%
\upsilon _{S^{n-1}}$-almost everywhere positive; see Lemma \ref{Proj_1} in
Appendix \ref{auxil}. Hence, it is automatically satisfied under elliptical
symmetry of $\mathfrak{P}$ provided a density is assumed to exist.

(ii) All the conditions on the density $p$ in Parts A.1-A.3 are certainly
satisfied under the conditions used in \cite{Mart10}.

(iii) Part B furthermore shows that under elliptical symmetry of $\mathfrak{P%
}$ the existence of a density is in fact not required at all.

(iii) If Assumptions \ref{ASDR} and \ref{ASCII} hold with the same square
root of $\Sigma \left( \cdot \right) $ (i.e., if $L(\cdot )=L_{\ast }(\cdot
) $ can be chosen in these assumptions) as is sometimes the case, then the
above theorem simplifies as $\mathcal{U}\left( L_{\ast }^{-1}L\right) $
reduces to the singleton $\left\{ I_{n}\right\} $.

(iv) Under the distributional assumptions for Part A of the preceding
theorem, if $\varphi =0$ (or $=1$) $\mu _{\mathbb{R}^{n}}$-almost everywhere
then trivially $E_{\beta ,\sigma ,\rho }(\varphi )=0$ (or $=1$) holds for
all $\beta $, $\sigma $, and $\rho $, and hence the same holds a fortiori
for the accumulation points. That the same is true under the distributional
assumptions for Part B can be seen as follows: By $G_{X}$-invariance of $%
\varphi $ and the assumptions for Part B we have that $E_{\beta ,\sigma
,\rho }(\varphi )=E\left( \varphi \left( L\left( \rho \right) \mathbf{G}%
\right) \right) $ where $\mathbf{G}$ is standard multivariate normal and $%
L\left( \rho \right) $ is nonsingular, cf. (\ref{reduction}) in Appendix \ref%
{app_proofs}. But then $E_{\beta ,\sigma ,\rho }(\varphi )=0$ (or $=1$)
follows (and the same is then a fortiori true for the limits).

(v) Similar as in Remark \ref{RMT}(ii) we make the trivial but sometimes
useful observation that the limiting power of a test $\varphi ^{\ast }$
which does not satisfy the assumptions of Theorem \ref{BT}\ can nevertheless
be computed from that theorem in an indirect way, if one can find another
test $\varphi $ that satisfies the assumptions of the theorem and differs
from $\varphi ^{\ast }$ only on a $\mathfrak{P}$-null set. This remark
obviously applies also to all other results in the paper and will not be
repeated.

(vi) For ways of extending the results in Part B of the preceding theorem to
the case where $\Pr \left( \mathbf{z}=0\right) >0$ see Remark \ref{rem_gen_1}%
(vi) in Section \ref{gen}. In a similar way Theorem \ref{BT2} and several
other results given further below can be extended to this case. We shall not
mention this again.
\end{remark}

The relationship of the preceding theorem to the second claim in MT1 is now
as follows: The additional invariance property (\ref{INV}) in the theorem is
automatically satisfied if $e\in \limfunc{span}(X)$ (by $G_{X}$-invariance
of $\varphi $). Furthermore, under the assumptions for MT1 and if $n>2$ the
vector $e$ in the preceding theorem coincides with $f_{1}(\Sigma ^{-1}(a-))$
considered in \cite{Mart10}, cf. Lemma \ref{convCM} and Remark \ref{RASCII}.
Hence, under Assumption \ref{ASC} (which is weaker than the corresponding
assumption in MT1) and if $n>2$, the preceding theorem provides a substitute
for the (incorrect) second claim in MT1 for the case where $e\in \limfunc{%
span}(X)$ if we specialize to $\varphi =\boldsymbol{1}_{\Phi }$. Recall from
Proposition \ref{prop_bound} that $\limfunc{span}(X)$ forms a part of $%
\limfunc{bd}(\Phi )$ for $G_{X}$-invariant rejection regions $\Phi $
satisfying $\emptyset \neq \Phi \neq \mathbb{R}^{n}$. We furthermore note
that the preceding theorem does not only deliver a qualitative statement
like that the limiting power is strictly between $0$ and $1$, but provides
an explicit formula for the limiting power (or the set of accumulation
points). We also point out that Theorem 2.11 in \cite{mynbaev2012} provides
related, but only qualitative, results for a certain class of rejection
regions.

As just discussed, the preceding theorem provides a substitute for the
second claim in MT1 in case $e$ belongs to that part of $\limfunc{bd}(\Phi )$
which is represented by $\limfunc{span}(X)$. If $e\in \limfunc{bd}(\Phi
)\backslash \limfunc{span}(X)$ then, for rejection regions $\Phi $ of the
form $\left\{ y\in \mathbb{R}^{n}:T(y)>\kappa \right\} $ with $T$ satisfying
a mild continuity property, Part 2 of Proposition \ref{prop_bound} shows
that $\kappa =T\left( e\right) $ must hold. [Part 3 of the same proposition
even shows that for the frequently used rejection regions $\Phi _{B,\kappa }$
the conditions $T_{B}\left( e\right) =\kappa $ and $e\notin \limfunc{span}%
(X) $ conversely imply $e\in \limfunc{bd}(\Phi _{B,\kappa })\backslash 
\limfunc{span}(X)$ provided $\emptyset \neq \Phi _{B,\kappa }\neq \mathbb{R}%
^{n}$.] Hence, if we can determine the limiting behavior of $P_{\beta
,\sigma ,\rho }\left( \left\{ y\in \mathbb{R}^{n}:T(y)>\kappa \right\}
\right) $ as $\rho \rightarrow a$ for the case where $\kappa =T(e)$, this
can then be used to obtain a substitute for the second claim in MT1 in case $%
e\in \limfunc{bd}(\Phi )\backslash \limfunc{span}(X)$, see the discussion
following the subsequent theorem. This theorem now provides such a limiting
result.\footnote{%
It is worth noting that the assumptions of this theorem per se do not imply
the assumption of Part 2 or Part 3 of Proposition \ref{prop_bound}.} Like in
the preceding theorem the rejection probabilities actually do neither depend
on $\beta $ nor $\sigma $.

\begin{theorem}
\label{BT2} Suppose Assumptions \ref{ASC} and \ref{ASCII} hold with the same
vector $e$, and Assumption \ref{ASDR} holds. Let $T$ be a test statistic
that is invariant w.r.t. $G_{X}$. Suppose there exists a positive integer $q$
and a homogeneous multivariate polynomial $D:\mathbb{R}^{n}\rightarrow 
\mathbb{R}$ of degree $q$, which does not vanish on all of $\limfunc{span}$$%
(e)^{\bot }$, such that for every $h\in \mathbb{R}^{n}$%
\begin{equation}
T(e+h)=T(e)+D(h)+R(h)  \label{ATaylor}
\end{equation}%
holds where $R(h)/\Vert {h}\Vert ^{q}\rightarrow 0$ as $h\rightarrow 0$, $%
h\neq 0$. Let $\mathcal{U}\left( L_{\ast }^{-1}L,\Sigma ^{-1/2}L\right) $
denote the set of all accumulation points of $\left( L_{\ast }^{-1}\left(
\rho \right) L\left( \rho \right) ,\Sigma ^{-1/2}\left( \rho \right) L\left(
\rho \right) \right) $ for $\rho \rightarrow a$. Furthermore, let $\beta \in 
\mathbb{R}^{k}$ and $0<\sigma <\infty $ be arbitrary but given.

\begin{enumerate}
\item Suppose the distribution of $\mathbf{z}$ (figuring in Assumption \ref%
{ASDR}) possesses a density $p$ w.r.t. Lebesgue measure $\mu _{\mathbb{R}%
^{n}}$. Then the accumulation points of 
\begin{equation}
P_{\beta ,\sigma ,\rho }\left( \left\{ y\in \mathbb{R}^{n}:T(y)>T(e)\right\}
\right)  \label{rej_prob}
\end{equation}%
for $\rho \rightarrow a$ are, in case $q$ is even, precisely given by 
\begin{equation}
\Pr \left( D(\Lambda U\mathbf{z})>0\right)  \label{even}
\end{equation}%
with $U\in \mathcal{U}\left( L_{\ast }^{-1}L\right) $ and where $\Lambda $
is as in Assumption \ref{ASCII}; for $q$ odd, they are precisely given by%
\begin{equation}
\Pr \left( D(\Lambda U\mathbf{z})>0,e^{\prime }U_{0}\mathbf{z}>0\right) +\Pr
\left( D(\Lambda U\mathbf{z})<0,e^{\prime }U_{0}\mathbf{z}<0\right)
\label{odd}
\end{equation}%
with $\left( U,U_{0}\right) \in \mathcal{U}\left( L_{\ast }^{-1}L,\Sigma
^{-1/2}L\right) $. Thus a sufficient condition for the limit of (\ref%
{rej_prob}) for $\rho \rightarrow a$ to exist for even $q$ is that $\mathcal{%
U}\left( L_{\ast }^{-1}L\right) $ is a singleton, whereas for odd $q$ it is
that $\mathcal{U}\left( L_{\ast }^{-1}L,\Sigma ^{-1/2}L\right) $ is a
singleton.

\item Suppose $\mathfrak{P}$ is an elliptically symmetric family with the
distribution of $\mathbf{z}$ satisfying $\Pr \left( \mathbf{z}=0\right) =0$.
Then, if $q$ is even, 
\begin{equation}
\lim_{\rho \rightarrow a}P_{\beta ,\sigma ,\rho }\left( \left\{ y\in \mathbb{%
R}^{n}:T(y)>T(e)\right\} \right) =\Pr (D(\Lambda \mathbf{z})>0)=\Pr
(D(\Lambda \mathbf{G})>0)  \label{even2}
\end{equation}%
holds where $\mathbf{G}$ is a multivariate Gaussian random vector with mean
zero and covariance matrix $I_{n}$. If $q$ is odd, the accumulation points
of $P_{\beta ,\sigma ,\rho }\left( \left\{ y\in \mathbb{R}%
^{n}:T(y)>T(e)\right\} \right) $ for $\rho \rightarrow a$ are precisely
given by%
\begin{eqnarray}
&&\Pr \left( D(\Lambda \mathbf{z})>0,e^{\prime }U_{0}\mathbf{z}>0\right)
+\Pr \left( D(\Lambda \mathbf{z})<0,e^{\prime }U_{0}\mathbf{z}<0\right) 
\notag \\
&=&\Pr \left( D(\Lambda \mathbf{G})>0,e^{\prime }U_{0}\mathbf{G}>0\right)
+\Pr \left( D(\Lambda \mathbf{G})<0,e^{\prime }U_{0}\mathbf{G}<0\right)
\label{odd2}
\end{eqnarray}%
with $U_{0}\in \mathcal{U}\left( \Sigma ^{-1/2}L_{\ast }\right) $ (and hence
the limit of the rejection probabilities for $\rho \rightarrow a$
necessarily exists if $\mathcal{U}\left( \Sigma ^{-1/2}L_{\ast }\right) $ is
a singleton). If $\Lambda U_{0}^{\prime }e=0$ holds for some $U_{0}\in 
\mathcal{U}\left( \Sigma ^{-1/2}L_{\ast }\right) $, the expression in (\ref%
{odd2}) with this $U_{0}$ then equals $1/2$. [A sufficient condition for $%
\Lambda U_{0}^{\prime }e=0$ to hold is that $\Lambda U_{0}^{\prime }$ is
symmetric.]

\item Suppose $\mathfrak{P}$ is an elliptically symmetric family with the
distribution of $\mathbf{z}$ satisfying $\Pr \left( \mathbf{z}=0\right) =0$.
If $q$ is odd and if, additionally, 
\begin{equation}
\lim_{\rho \rightarrow a}\lambda _{n}^{-1/2}\left( \Sigma \left( \rho
\right) \right) c\left( \rho \right) \Pi _{\limfunc{span}\left( e\right)
^{\bot }}\Sigma \left( \rho \right) \Pi _{\limfunc{span}\left( e\right) }=0
\label{offdiag}
\end{equation}%
holds, where $c\left( \rho \right) $ is as in Assumption \ref{ASCII}, then%
\begin{equation*}
\lim_{\rho \rightarrow a}P_{\beta ,\sigma ,\rho }\left( \left\{ y\in \mathbb{%
R}^{n}:T(y)>T(e)\right\} \right) =1/2\text{.}
\end{equation*}
\end{enumerate}
\end{theorem}

We recall from Remark \ref{RASCII} that assuming that the vector $e$ is the
same in Assumptions \ref{ASC} and \ref{ASCII} entails no loss of generality
provided $n>2$. Condition (\ref{offdiag}) ensures that $\Lambda
U_{0}^{\prime }e=0$ holds for every $U_{0}\in \mathcal{U}\left( \Sigma
^{-1/2}L_{\ast }\right) $, cf. Lemma \ref{off-diag} in Appendix \ref%
{app_proofs}, which can also be used to formulate conditions equivalent to (%
\ref{offdiag}). This can be useful if one of these equivalent formulations
is easier to verify in a particular application. It will turn out that
condition (\ref{offdiag}) holds for autoregressive models of order $1$ and
certain classes of spatial error models, see Sections \ref{SEM} and \ref%
{AR1C}. Furthermore note that under the assumption that $\mathfrak{P}$ is an
elliptically symmetric family the existence of a density is not required in
the preceding theorem.

Observe that, under the assumptions of MT1, the vector $e$ in the preceding
theorem coincides with $f_{1}(\Sigma ^{-1}(a-))$ considered in \cite{Mart10}%
, cf. Lemma \ref{convCM}. Hence, for rejection regions of the form $\left\{
y\in \mathbb{R}^{n}:T(y)>\kappa \right\} $ with $T$ satisfying the
assumptions of Proposition \ref{prop_bound} as well as of the preceding
theorem, this theorem provides a substitute for the (incorrect) second claim
in MT1 in case $e\in \limfunc{bd}(\Phi )\backslash \limfunc{span}(X)$ in
that it determines the limit (or the set of accumulation points) of the
power function as $\rho \rightarrow a$. Note that the theorem itself does
not in general make a statement about the limiting expressions always being
strictly between $0$ and $1$; however, given the explicit expressions for
the accumulation points of the power function in the preceding theorem, this
can then be decided on a case by case basis. [We note that cases exist where
the above theorem applies and the limiting power is zero or one, see, e.g.,
Corollary \ref{Lem_Illust_3}, Part 1, in case $\lambda =\lambda _{1}\left(
B\right) $.]

\begin{remark}
\label{Rem_BT2}\emph{(Comments on the assumption on }$T$\emph{) }(i)\ We
note that under the assumptions of Theorem \ref{BT2} the polynomial $D$ in (%
\ref{ATaylor}) necessarily vanishes everywhere on $\limfunc{span}(e)$. More
generally, $D\left( h\right) =D\left( \Pi _{\limfunc{span}(e)^{\bot
}}h\right) $ holds for every $h\in \mathbb{R}^{n}$; see Lemma \ref{T} in
Appendix \ref{app_proofs}.

(ii) If $T$ is a test statistic that is totally differentiable at $e$, it
satisfies relation (\ref{ATaylor}) with $q=1$ and $D(h)=d^{\prime }h$, $d$ a 
$n\times 1$ vector. If $d\notin \limfunc{span}(e)$, then $D$ satisfies the
assumption of the theorem. In case $d\in \limfunc{span}(e)$ this is not so,
since $D$ then vanishes identically on $\limfunc{span}$$(e)^{\bot }$ (in
fact, $d=0$ must then hold provided $T$ is $G_{X}$-invariant). In this case
one can try to resort to higher order Taylor expansions: For example, if $T$
is twice continuously partially differentiable in a neighborhood of $e$,
then $D$ can be chosen as $1/2$ times the quadratic form corresponding to
the Hessian matrix of $T$ at the point $e$, provided that $D$ does not
vanish identically on $\limfunc{span}(e)^{\bot }$.

(iii) In Theorem \ref{BT2} the element $e$ does not belong to the rejection
region by construction. In case $q$ is odd, $e$ always belongs to the
boundary of that region in view of homogeneity of $D$. The same is true in
case $q$ is even provided $D\left( h\right) >0$ holds for some $h\in \mathbb{%
R}^{n}$ (which by (i) above is equivalent to $D\left( h\right) >0$ for some $%
h\in \limfunc{span}$$(e)^{\bot }$ with $h\neq 0$). If $q$ is even and $%
D\left( h\right) <0$ holds for all $h\notin \limfunc{span}$$(e)$ (which by
(i) above is equivalent to $D\left( h\right) <0$ for all $h\in \limfunc{span}
$$(e)^{\bot }$ with $h\neq 0$), Lemma \ref{T} in Appendix \ref{app_proofs}
shows that then $e$ is not an element of the boundary, but is an element of
the exterior (i.e., of the complement of the closure) of the rejection
region. In the remaining case, i.e., $q$ even and $D\left( h\right) \leq 0$
for all $h\notin \limfunc{span}$$(e)$ but $D\left( h\right) =0$ for some $%
h\notin \limfunc{span}$$(e)$, no conclusion can be drawn in general.
\end{remark}

In the following example we illustrate how the assumptions on $T$ in the
preceding theorem can be verified for the important class of test statistics 
$T_{B}$.

\begin{example}
\label{EX_quad}We consider the test statistic $T_{B}=T_{B,C_{X}}$ given by (%
\ref{T_quadratic}). We assume that $B$ is not a multiple of $I_{n-k}$, since
otherwise $T_{B}$ is constant which is a trivial case. If $e\in \limfunc{span%
}(X)$ holds, then $T_{B}$ is not even continuous at $e$, showing that
condition (\ref{ATaylor}) can not be satisfied. We hence assume $e\notin 
\limfunc{span}(X)$. Elementary calculations show that then (\ref{ATaylor})
with $q=1$ and 
\begin{equation}
D\left( h\right) =2\left\Vert C_{X}e\right\Vert ^{-2}\left( e^{\prime
}C_{X}^{\prime }BC_{X}-\left\Vert C_{X}e\right\Vert ^{-2}\left( e^{\prime
}C_{X}^{\prime }BC_{X}e\right) e^{\prime }C_{X}^{\prime }C_{X}\right) h
\label{D_1}
\end{equation}%
holds. In view of $D\left( e\right) =0$ and surjectivity of $C_{X}$ we see
that $D$ does not vanish on all of $\limfunc{span}$$(e)^{\bot }$ if and only
if $e^{\prime }C_{X}^{\prime }B\neq \left\Vert C_{X}e\right\Vert ^{-2}\left(
e^{\prime }C_{X}^{\prime }BC_{X}e\right) e^{\prime }C_{X}^{\prime }$, or in
other words if and only if $C_{X}e$ is not an eigenvector of $B$, a
condition that can easily be checked. As a point of interest we note that in
this case $T_{B}(e)\in (\lambda _{1}(B),\lambda _{n-k}(B))$ must hold,
entailing that $\emptyset \neq \Phi _{B,T_{B}(e)}\neq \mathbb{R}^{n}$ (in
fact, neither $\Phi _{B,T_{B}(e)}$ nor its complement are $\mu _{\mathbb{R}%
^{n}}$-null sets, cf. Remark \ref{range_for_kappa}). It then follows from
Proposition \ref{prop_bound} that $e$ is an element of the boundary of $\Phi
_{B,T_{B}(e)}$. [This can alternatively be deduced from Remark \ref{Rem_BT2}%
(iii).] Next consider the case where $C_{X}e$ is an eigenvector of $B$ with
eigenvalue $\lambda $. Applying now (\ref{ATaylor}) with $q=2$ leads to 
\begin{equation}
D\left( h\right) =\left\Vert C_{X}e\right\Vert ^{-2}h^{\prime }\left(
C_{X}^{\prime }BC_{X}-\lambda C_{X}^{\prime }C_{X}\right) h  \label{D_2}
\end{equation}%
which is homogeneous of degree $q=2$ and which does not vanish on all of $%
\limfunc{span}$$(e)^{\bot }$ (except if $B=\lambda I_{n-k}$, a case we have
ruled out). We note that now $T_{B}(e)=\lambda \in \lbrack \lambda
_{1}(B),\lambda _{n-k}(B)]$ must hold. Recall from Remark \ref%
{range_for_kappa} that in case $\lambda $ is not the largest eigenvalue of $%
B $, we know that $\emptyset \neq \Phi _{B,T_{B}(e)}\neq \mathbb{R}^{n}$ (in
fact, neither $\Phi _{B,T_{B}(e)}$ nor its complement are $\mu _{\mathbb{R}%
^{n}}$-null sets if additionally $\lambda >\lambda _{1}(B)$ holds, whereas $%
\Phi _{B,T_{B}(e)}$ is the complement of a non-empty $\mu _{\mathbb{R}^{n}}$%
-null set if $\lambda =\lambda _{1}(B)$). Hence Proposition \ref{prop_bound}
shows that $e$ then belongs to the boundary of $\Phi _{B,T_{B}(e)}$. [Since $%
D\left( h\right) >0$ holds for some $h\in \mathbb{R}^{n}$ if $\lambda $ is
not the largest eigenvalue of $B$, this can alternatively be deduced from
Remark \ref{Rem_BT2}(iii).] In case $\lambda $ is the largest eigenvalue of $%
B$, then $\Phi _{B,T_{B}(e)}$ is empty. The last case shows that, although
Theorem \ref{BT2} is geared to the case where $e$ belongs to the boundary of
the rejection region, its assumptions do not rule out other cases.
Furthermore, the case where $C_{X}e$ is an eigenvector of $B$ with
eigenvalue $\lambda $ satisfying $\lambda =\lambda _{1}(B)$ shows that
Theorem \ref{BT2} also applies to cases where, although $e$ belongs to the
boundary of the rejection region, the limiting rejection probabilities are
not necessarily in $(0,1)$. $\square $
\end{example}

\begin{remark}
\label{rem_accu}\emph{(Comments on the set of accumulation points) }(i) If
one can choose $L_{\ast }=\Sigma ^{1/2}$ in the second part of the preceding
theorem then $\mathcal{U}\left( \Sigma ^{-1/2}L_{\ast }\right) $ reduces to
the singleton $\left\{ I_{n}\right\} $ and the statement in Part 2
simplifies accordingly. A similar remark applies to the first part of the
theorem in case $L_{\ast }=L$ and/or $L=\Sigma ^{1/2}$ can be chosen.

(ii) It is not difficult to see that the accumulation points as given in (%
\ref{even}) and (\ref{odd}) depend continuously on $U$ and $U_{0}$. [This
follows from the portmanteau theorem observing that $D\left( \Lambda U%
\mathbf{z}\right) $ as well as $e^{\prime }U_{0}\mathbf{z}$ depend
continuously on $U$ and $U_{0}$, respectively, and that both expressions are
nonzero almost surely as shown in the proof of Theorem \ref{BT2}.] Since $%
\mathcal{U}\left( L_{\ast }^{-1}L\right) $ as well as $\mathcal{U}\left(
L_{\ast }^{-1}L,\Sigma ^{-1/2}L\right) $ are compact, the question of
whether or not the set of accumulation points is bounded away from $0$ (or $%
1 $, respectively) then just reduces to the question as to whether every
accumulation point is larger than $0$ (smaller than $1$, respectively). The
latter question can often easily be answered by examining the explicit
expressions provided by (\ref{even}) and (\ref{odd}). For an example see the
remark immediately below.

(iii) Suppose $q=1$ in the second part of the theorem. Observe that then $%
D(h)=d^{\prime }h$ with $d\notin \limfunc{span}(e)$ by Remark \ref{Rem_BT2}%
(ii). Hence, in case $d^{\prime }\Lambda $ and $e^{\prime }U_{0}$ are not
collinear, the accumulation point given by (\ref{odd2}) is in the open
interval $(0,1)$. If $d^{\prime }\Lambda $ and $e^{\prime }U_{0}$ are
collinear, then the accumulation point is either $0$ or $1$.
\end{remark}

\subsubsection{An illustration for tests based on $T_{B}$\label{T_B}}

We now illustrate the results obtained so far by applying them to tests
based on the statistic $T_{B}=T_{B,C_{X}}$ defined in (\ref{T_quadratic}).
We note that, under regularity conditions (including appropriate
distributional assumptions) and excluding degenerate cases, point-optimal
invariant tests and locally best invariant tests are of this form with $%
B=-\left( C_{X}\Sigma (\bar{\rho})C_{X}^{\prime }\right) ^{-1}$ and $B=C_{X}%
\dot{\Sigma}(0)C_{X}^{\prime }$, respectively, with $\dot{\Sigma}(0)$
denoting the derivative at $\rho =0$ (ensured to exist under the
aforementioned regularity conditions), see, e.g., \cite{KingHillier1985}.%
\footnote{%
These tests are point-optimal (locally best) in the class of all $G_{X}^{+}$%
-invariant tests. As they are also $G_{X}$-invariant, they are a fortiori
also point-optimal (locally best) tests in the class of $G_{X}$-invariant
tests.}

Recall that under the assumptions in \cite{Mart10} the vector $e$ given by
Assumption \ref{ASC} corresponds to the eigenvector $f_{1}(\Sigma ^{-1}(a-))$
in MT1, possibly up to a sign change. For that reason we impose Assumption %
\ref{ASC} in all of the three corollaries that follow, although this
assumption would not be needed for the second one of the corollaries (but
note that then $e$ would be determined by Assumption \ref{ASCII} only).
Furthermore, recall from Remark \ref{range_for_kappa} that $\emptyset \neq
\Phi _{B,\kappa }\neq \mathbb{R}^{n}$ occurs if and only if $\kappa \in
\lbrack \lambda _{1}(B),\lambda _{n-k}(B))$ (the interval being non-empty if
and only if $\lambda _{1}(B)<\lambda _{n-k}(B)$). We shall in the following
corollaries hence always assume that $\kappa $ is in that range and thus
shall exclude the trivial cases where $\Phi _{B,\kappa }=\emptyset $ or $%
\Phi _{B,\kappa }=\mathbb{R}^{n}$ from the formulation of the corollaries.

The first corollary is based on Theorem \ref{MT}. Recall that the conditions
in this corollary are weaker than the conditions used in MT1 (cf. Remark \ref%
{RMT}) and that sufficient conditions for the high-level Assumption \ref{ASD}
have been given in Proposition \ref{AsPexII} (under which the rejection
probabilities actually do neither depend on $\beta $ nor $\sigma $).

\begin{corollary}
\label{Lem_Illust_1}Suppose Assumptions \ref{ASC} and \ref{ASD} are
satisfied. Assume that $\kappa \in \lbrack \lambda _{1}(B),\lambda
_{n-k}(B)) $ with $\lambda _{1}(B)<\lambda _{n-k}(B)$. Then we have:

\begin{enumerate}
\item $T_{B}\left( e\right) >\kappa $ (i.e., $e\in \limfunc{int}\left( \Phi
_{B,\kappa }\right) $) implies $\lim_{\rho \rightarrow a}P_{\beta ,\sigma
,\rho }\left( \Phi _{B,\kappa }\right) =1$ for every $\beta \in \mathbb{R}%
^{k}$ and $0<\sigma <\infty $.\footnote{%
Note that $T_{B}\left( e\right) >\kappa $ entails $e\notin \limfunc{span}%
\left( X\right) $ in view of (\ref{T_quadratic}) and $\kappa \geq \lambda
_{1}\left( B\right) $.}

\item $T_{B}\left( e\right) <\kappa $ and $e\notin \limfunc{span}(X)$ (i.e., 
$e\notin \limfunc{cl}\left( \Phi _{B,\kappa }\right) $) implies $\lim_{\rho
\rightarrow a}P_{\beta ,\sigma ,\rho }\left( \Phi _{B,\kappa }\right) =0$
for every $\beta \in \mathbb{R}^{k}$ and $0<\sigma <\infty $.
\end{enumerate}
\end{corollary}

It is worth pointing out here that the second case, i.e., the zero-power
trap, can occur even for point-optimal invariant or locally best invariant
tests as has been documented in the literature cited in the introduction.
The next two corollaries now deal with the case where $e$ belongs to the
boundary of the rejection region. They are based on Theorems \ref{BT} and %
\ref{BT2}, respectively. For simplicity of presentation we concentrate only
on the case of elliptically symmetric families. We remind the reader that in
the two subsequent corollaries the rejection probabilities actually neither
depend on $\beta $ nor $\sigma $, i.e., $P_{\beta ,\sigma ,\rho }\left( \Phi
_{B,\kappa }\right) =P_{0,1,\rho }\left( \Phi _{B,\kappa }\right) $ holds.

\begin{corollary}
\label{Lem_Illust_2}Suppose Assumptions \ref{ASC} and \ref{ASCII} are
satisfied with the same vector $e$.\footnote{%
For reasons of conformity we have here included the condition that the
vector $e$ is the same in both assumptions, although this does not impose a
restriction here. This is so because of Remark \ref{RASCII} and since $n>2$
must hold in this corollary: Suppose $n=2$ would hold. Then $k=1$ would
follow in view of $0\leq k<n$ and the assumption $e\in \limfunc{span}(X)$.
But this would be in conflict with $\lambda _{1}(B)<\lambda _{n-k}(B)$.}
Furthermore, assume that $\mathfrak{P}$ is an elliptically symmetric family
(i.e., Assumption \ref{ASDR} holds with a spherically distributed $\mathbf{z}
$) and $\Pr \left( \mathbf{z}=0\right) =0$. Assume that $\kappa \in \lbrack
\lambda _{1}(B),\lambda _{n-k}(B))$ with $\lambda _{1}(B)<\lambda _{n-k}(B)$%
. Suppose $e\in \limfunc{span}(X)$ holds. Then $\lim_{\rho \rightarrow
a}P_{\beta ,\sigma ,\rho }\left( \Phi _{B,\kappa }\right) $ exists and
equals $\Pr \left( T_{B}\left( \Lambda \mathbf{G}\right) >\kappa \right) $
where $\mathbf{G}$ is a multivariate Gaussian random vector with mean zero
and covariance matrix $I_{n}$. Furthermore, the limit satisfies%
\begin{equation*}
0<\lim_{\rho \rightarrow a}P_{\beta ,\sigma ,\rho }\left( \Phi _{B,\kappa
}\right) <1
\end{equation*}%
provided $\kappa >\lambda _{1}(B)$, whereas it equals $1$ in case $\kappa
=\lambda _{1}(B)$.\footnote{\label{footnote_b}In case $\kappa =\lambda
_{1}(B)$ the rejection region is the complement of a $\mu _{\mathbb{R}^{n}}$%
-null set. As discussed in Remark \ref{RBT}(iv), we then even have $P_{\beta
,\sigma ,\rho }\left( \Phi _{B,\kappa }\right) =1$ for every $\beta $, $%
\sigma $, and $\rho $ (although we do not require $\mathbf{z}$ to possess a
density).}
\end{corollary}

The next result covers the case where $e\in \limfunc{bd}(\Phi _{B,\kappa
})\backslash \limfunc{span}(X)$. Recall from Proposition \ref{prop_bound}
and Remark \ref{range_for_kappa} that this is equivalent to $e\notin 
\limfunc{span}(X)$ and $\kappa =T_{B}\left( e\right) \in \lbrack \lambda
_{1}(B),\lambda _{n-k}(B))$ with $\lambda _{1}(B)<\lambda _{n-k}(B)$. Note
that $T_{B}\left( e\right) \in \lbrack \lambda _{1}(B),\lambda _{n-k}(B)]$
always holds by definition of $T_{B}$.

\begin{corollary}
\label{Lem_Illust_3}Suppose Assumptions \ref{ASC} and \ref{ASCII} are
satisfied with the same vector $e$. Furthermore, assume that $\mathfrak{P}$
is an elliptically symmetric family (i.e., Assumption \ref{ASDR} holds with
a spherically distributed $\mathbf{z}$) and $\Pr \left( \mathbf{z}=0\right)
=0$. Assume $e\in \limfunc{bd}(\Phi _{B,\kappa })\backslash \limfunc{span}%
(X) $ (i.e., $e\notin \limfunc{span}(X)$ and $\kappa =T_{B}\left( e\right)
\in \lbrack \lambda _{1}(B),\lambda _{n-k}(B))$ with $\lambda
_{1}(B)<\lambda _{n-k}(B)$ hold).

\begin{enumerate}
\item Suppose $C_{X}e$ is an eigenvector of $B$ with eigenvalue $\lambda $,
say. Then $\lambda =T_{B}\left( e\right) =\kappa $ and%
\begin{equation}
\lim_{\rho \rightarrow a}P_{\beta ,\sigma ,\rho }\left( \Phi _{B,\kappa
}\right) =\Pr \left( \mathbf{G}^{\prime }\Lambda ^{\prime }\left(
C_{X}^{\prime }BC_{X}-\lambda C_{X}^{\prime }C_{X}\right) \Lambda \mathbf{G}%
>0\right) ,  \label{even_special}
\end{equation}%
where $\mathbf{G}$ is a multivariate Gaussian random vector with mean zero
and covariance matrix $I_{n}$. Furthermore, the limit belongs to the open
interval $(0,1)$ if $\lambda >\lambda _{1}(B)$, whereas it equals $1$ in
case $\lambda =\lambda _{1}(B)$.\footnote{%
Cf. Footnote \ref{footnote_b}.}

\item Suppose $C_{X}e$ is not an eigenvector of $B$. Then the accumulation
points of $P_{\beta ,\sigma ,\rho }\left( \Phi _{B,\kappa }\right) $ for $%
\rho \rightarrow a$ are given by%
\begin{eqnarray}
&&\Pr \left( \left( e^{\prime }C_{X}^{\prime }BC_{X}-\left\Vert
C_{X}e\right\Vert ^{-2}\left( e^{\prime }C_{X}^{\prime }BC_{X}e\right)
e^{\prime }C_{X}^{\prime }C_{X}\right) \Lambda \mathbf{G}>0,e^{\prime }U_{0}%
\mathbf{G}>0\right) +  \notag \\
&&\Pr \left( \left( e^{\prime }C_{X}^{\prime }BC_{X}-\left\Vert
C_{X}e\right\Vert ^{-2}\left( e^{\prime }C_{X}^{\prime }BC_{X}e\right)
e^{\prime }C_{X}^{\prime }C_{X}\right) \Lambda \mathbf{G}<0,e^{\prime }U_{0}%
\mathbf{G}<0\right)  \label{odd_special}
\end{eqnarray}%
with $U_{0}\in \mathcal{U}\left( \Sigma ^{-1/2}L_{\ast }\right) $. The
expression in (\ref{odd_special}) is in the open interval $(0,1)$ for every $%
U_{0}\in \mathcal{U}\left( \Sigma ^{-1/2}L_{\ast }\right) $ which has the
property that $\left( e^{\prime }C_{X}^{\prime }BC_{X}-\left\Vert
C_{X}e\right\Vert ^{-2}\left( e^{\prime }C_{X}^{\prime }BC_{X}e\right)
e^{\prime }C_{X}^{\prime }C_{X}\right) \Lambda $ and $e^{\prime }U_{0}$ are
not collinear.\footnote{%
If these two vectors are collinear, then the expression in (\ref{odd_special}%
) is $0$ or $1$.} [This non-collinearity is, in particular, the case if $%
\Lambda U_{0}^{\prime }e=0$ holds, in which case the expression in (\ref%
{odd_special}) equals $1/2$.] Furthermore, the set of all accumulation
points is bounded away from $0$ and $1$ provided $\left( e^{\prime
}C_{X}^{\prime }BC_{X}-\left\Vert C_{X}e\right\Vert ^{-2}\left( e^{\prime
}C_{X}^{\prime }BC_{X}e\right) e^{\prime }C_{X}^{\prime }C_{X}\right)
\Lambda $ and $e^{\prime }U_{0}$ are not collinear for every $U_{0}\in 
\mathcal{U}\left( \Sigma ^{-1/2}L_{\ast }\right) $.

\item Suppose $C_{X}e$ is not an eigenvector of $B$. If, additionally, 
\begin{equation*}
\lim_{\rho \rightarrow a}\lambda _{n}^{-1/2}\left( \Sigma \left( \rho
\right) \right) c\left( \rho \right) \Pi _{\limfunc{span}\left( e\right)
^{\bot }}\Sigma \left( \rho \right) \Pi _{\limfunc{span}\left( e\right) }=0
\end{equation*}%
holds, where $c\left( \rho \right) $ is as in Assumption \ref{ASCII}, then 
\begin{equation*}
\lim_{\rho \rightarrow a}P_{\beta ,\sigma ,\rho }\left( \Phi _{B,\kappa
}\right) =1/2.
\end{equation*}
\end{enumerate}
\end{corollary}

In the preceding corollary we have excluded the case where $\kappa
=T_{B}\left( e\right) =\lambda _{n-k}(B)>\lambda _{1}(B)$. While we already
know that this is a trivial case as then $\Phi _{B,\kappa }$ is empty, it is
interesting to note that even in this case the proof of the above corollary,
which is based on Theorem \ref{BT2}, would still go through and would
deliver (\ref{even_special}), which -- as it should -- would then reduce to
zero since the matrix $\Lambda ^{\prime }\left( C_{X}^{\prime
}BC_{X}-\lambda C_{X}^{\prime }C_{X}\right) \Lambda $ is non-positive
definite in this case.\footnote{%
The proof of the corollary makes use of Example \ref{EX_quad} which assumes $%
e\notin \limfunc{span}(X)$. Note that $\kappa =T_{B}\left( e\right) =\lambda
_{n-k}(B)$ implies $e\notin \limfunc{span}(X)$ if $\lambda _{1}(B)<\lambda
_{n-k}(B)$, allowing one to directly extend the proof of the corollary to
this case.}

\begin{remark}
(i) Appropriate versions of Corollaries \ref{Lem_Illust_1}-\ref{Lem_Illust_3}
can also be given for a test statistic $T_{B}^{\prime }$ that takes a value $%
c\neq \lambda _{1}(B)$ on all of $\limfunc{span}\left( X\right) $ and
coincides with $T_{B}$ on the complement of $\limfunc{span}\left( X\right) $%
. For example, in such a version of Corollary \ref{Lem_Illust_1} one needs
to add the assumption $e\notin \limfunc{span}\left( X\right) $ in Part 1 of
that corollary, because there is then no guarantee that the condition $%
T_{B}^{\prime }\left( e\right) >\kappa $ is equivalent to $e\in \limfunc{int}%
\left( \Phi _{B,\kappa }^{\prime }\right) $.

(ii) In case $\limfunc{span}\left( X\right) $ is a $\mathfrak{P}$-null set
(which is, e.g., the case under the assumptions of Corollaries \ref%
{Lem_Illust_2} and \ref{Lem_Illust_3}, cf Remark \ref{E1}(iii)) we have $%
P_{\beta ,\sigma ,\rho }\left( \Phi _{B,\kappa }^{\prime }\right) =P_{\beta
,\sigma ,\rho }\left( \Phi _{B,\kappa }\right) $. Applying the above
corollaries as they stand to $T_{B}$ thus immediately provides information
on $P_{\beta ,\sigma ,\rho }\left( \Phi _{B,\kappa }^{\prime }\right) $
without the need of obtaining appropriate versions of the above corollaries
for $T_{B}^{\prime }$.
\end{remark}

\subsubsection{On the relationship between the size of a test and the
zero-power trap\label{alpha_star}}

Given a $G_{X}$-invariant test statistic $T$, we have seen in previous
sections that the limiting power of the test with rejection region $\Phi
_{\kappa }:=\left\{ y\in \mathbb{R}^{n}:T(y)>\kappa \right\} $ can be zero
(zero-power trap). Of course, an important question to ask is for which
critical values $\kappa $ this occurs. An (essentially) equivalent
formulation is to ask for which values of the sizes of the rejection regions 
$\Phi _{\kappa }$ the limiting power is zero; i.e., for which values of the
sizes the zero-power trap arises (at least along a subsequence of values of $%
\rho $). To this end we define 
\begin{equation}
\alpha ^{\ast }(T)=\inf \left\{ P_{0,1,0}(\Phi _{\kappa }):\kappa \in 
\mathbb{R}\text{ and }\liminf\limits_{\rho \rightarrow a}P_{0,1,\rho }(\Phi
_{\kappa })>0\right\} .  \label{astar}
\end{equation}%
We note that whenever the rejection probabilities $P_{\beta ,\sigma ,\rho
}(\Phi _{\kappa })$ are independent of $\beta $ and $\sigma $, which is
often the case (e.g., under Assumption \ref{ASDR}, see Remark \ref%
{invariance}), the quantity $\alpha ^{\ast }(T)$ is the infimum of the sizes
of all rejection regions $\Phi _{\kappa }$, the limiting power of which does
not vanish. Thus $\alpha ^{\ast }(T)$ describes the size where a phase
transition occurs: for sizes above $\alpha ^{\ast }(T)$ the zero-power trap
does not occur, while it occurs for sizes below $\alpha ^{\ast }(T)$ (at
least along a subsequence).\footnote{%
While $\alpha >\alpha ^{\ast }(T)$ implies that the zero-power trap does not
occur, it may in general still be the case that the limiting power is very
low.} We investigate properties of $\alpha ^{\ast }(T)$ in this section.

Before proceeding we note that in the more narrow context of spatial
regression models \cite{Mart10} also discusses the quantity $\alpha ^{\ast
}(T)$ in his Lemmata D.2 and D.3, which provide the basis for a large part
of the results beyond Theorem 1 in that reference.\footnote{%
That $\alpha ^{\ast }(T)$ defined above is indeed equivalent to the quantity 
$\alpha ^{\ast }$ described in \cite{Mart10}, p. 165, is discussed in
Appendix \ref{A1}.} Unfortunately, these lemmata are inappropriately stated
and the proofs contain several errors. We discuss this in detail in Appendix %
\ref{A1}. In the present section we provide correct versions of these two
lemmata, simultaneously freeing them from the spatial context, thus making
them applicable to much more general covariance structures.

The subsequent lemma can now be seen as a general version of Lemma D.2 in 
\cite{Mart10}. It gives an expression for $\alpha ^{\ast }(T)$ and shows
that -- under the assumptions of the lemma -- for every $\kappa $ with $%
P_{0,1,0}(\Phi _{\kappa })>\alpha ^{\ast }(T)$ the limiting power is not
only positive but in fact equals one.

\begin{lemma}
\label{D2new} Suppose Assumptions \ref{ASC} and \ref{ASD} are satisfied and
let $T:\mathbb{R}\rightarrow \mathbb{R}$ be a test statistic that is
invariant w.r.t. $G_{X}$. Consider the family of rejection regions 
\begin{equation*}
\Phi _{\kappa }=\left\{ y\in \mathbb{R}^{n}:T(y)>\kappa \right\}
\end{equation*}%
for $\kappa \in \mathbb{R}$. Suppose there exists a $\delta >0$ such that $%
e\notin \limfunc{bd}(\Phi _{\kappa })$ holds for every $0<|\kappa
-T(e)|<\delta $ where $e$ is the vector figuring in Assumption \ref{ASC}.
[This is, in particular, satisfied if $e\notin \limfunc{span}$$(X)$ and $T$
is continuous on $\mathbb{R}^{n}\backslash \limfunc{span}$$(X)$.] If the
cumulative distribution function of $P_{0,1,0}\circ T$ is continuous at $%
T(e) $, then 
\begin{equation*}
\alpha ^{\ast }(T)=P_{0,1,0}(\Phi _{T(e)}).
\end{equation*}%
Furthermore, if for some $\kappa $ we have $P_{0,1,0}(\Phi _{\kappa
})>\alpha ^{\ast }(T)$ ($<\alpha ^{\ast }(T)$, respectively), then $\kappa
<T(e)$ ($\kappa >T(e)$, respectively) and $\lim_{\rho \rightarrow
a}P_{0,1,\rho }(\Phi _{\kappa })=1$ ($=0$, respectively) hold.
\end{lemma}

The next result, which is based on Lemma \ref{D2new} above, considers the
test statistic $T_{B}=T_{B,C_{X}}$ and, in particular, characterizes
situations when the zero-power trap occurs or does not occur at all
significance levels. Restricted to regression models with spatial
autoregressive errors of order one, the subsequent lemma contains a correct
and improved version of Lemma D.3 in \cite{Mart10} as a special case, the
improvement relating amongst others to the fact that we do not only
characterize when $\alpha ^{\ast }\left( T_{B}\right) $ equals $0$ or $1$,
but that we also determine the limiting power in each case. Before
presenting the result, we note that Lemma D.3 in \cite{Mart10} is stated for
tests obtained from $T_{B}$ by rejecting for small values of the test
statistic while we state our result for tests that reject for large values
of $T_{B}$. However, this is immaterial as Lemma D.3 in \cite{Mart10} can
trivially be rephrased in our setting by simply passing from $B$ to $-B$. In
the subsequent two propositions we exclude the trivial case where $\lambda
_{1}(B)=\lambda _{n-k}(B)$ holds, in which case $\alpha ^{\ast }\left(
T_{B}\right) =1$. [To see this note that then $T_{B}$ is constant equal to $%
\lambda _{1}(B)$ and thus all rejection probabilities are zero or one
depending on whether $\kappa \geq \lambda _{1}(B)$ or $\kappa <\lambda
_{1}(B)$.]

\begin{proposition}
\label{D3new} Suppose Assumptions \ref{ASC} and \ref{ASD} hold. Furthermore,
assume that $P_{0,1,0}$ is $\mu _{\mathbb{R}^{n}}$-absolutely continuous
with a density that is positive on an open neighborhood of the origin except
possibly for a $\mu _{\mathbb{R}^{n}}$-null set. Suppose $e\notin \limfunc{%
span}(X)$ where $e$ is the vector figuring in Assumption \ref{ASC} and
suppose that $\lambda _{1}(B)<\lambda _{n-k}(B)$ holds. Then:

\begin{enumerate}
\item $\alpha ^{\ast }\left( T_{B}\right) =0$ if and only if $C_{X}e\in 
\limfunc{Eig}\left( B,\lambda _{n-k}(B)\right) $. If $C_{X}e\in \limfunc{Eig}%
\left( B,\lambda _{n-k}(B)\right) $ holds, then $\lim_{\rho \rightarrow
a}P_{0,1,\rho }(\Phi _{B,\kappa })=1$ for every $\kappa \in (-\infty
,\lambda _{n-k}(B))$. [For $\kappa \geq \lambda _{n-k}(B)$ we trivially
always have $\Phi _{B,\kappa }=\emptyset $.]

\item $\alpha ^{\ast }\left( T_{B}\right) =1$ if and only if $C_{X}e\in 
\limfunc{Eig}\left( B,\lambda _{1}(B)\right) $. If $C_{X}e\in \limfunc{Eig}%
\left( B,\lambda _{1}(B)\right) $ holds, then $\lim_{\rho \rightarrow
a}P_{0,1,\rho }(\Phi _{B,\kappa })=0$ for every $\kappa \in (\lambda
_{1}(B),\infty )$. [For $\kappa <\lambda _{1}(B)$ we trivially always have $%
\Phi _{B,\kappa }=\mathbb{R}^{n}$, whereas $\Phi _{B,\kappa }$ is the
complement of a $\mu _{\mathbb{R}^{n}}$-null set in case $\kappa =\lambda
_{1}(B)<\lambda _{n-k}(B)$.\footnote{%
Hence, $P_{0,1,0}(\Phi _{B,\kappa })=1$ holds in case $\kappa =\lambda
_{1}(B)<\lambda _{n-k}(B)$. Furthermore, $\lim_{\rho \rightarrow
a}P_{0,1,\rho }(\Phi _{B,\kappa })=1$ will then also hold provided, e.g.,
the measures $P_{0,1,\rho }$ are $\mu _{\mathbb{R}^{n}}$-absolutely
continuous.}]

\item $0<\alpha ^{\ast }\left( T_{B}\right) <1$ if and only if $C_{X}e$
neither belongs to $\limfunc{Eig}\left( B,\lambda _{1}(B)\right) $ nor $%
\limfunc{Eig}\left( B,\lambda _{n-k}(B)\right) $. If $C_{X}e$ neither
belongs to $\limfunc{Eig}\left( B,\lambda _{1}(B)\right) $ nor $\limfunc{Eig}%
\left( B,\lambda _{n-k}(B)\right) $, there exists a unique $\kappa ^{\ast
}\in \left( \lambda _{1}(B),\lambda _{n-k}(B)\right) $ such that $%
P_{0,1,0}(\Phi _{B,\kappa ^{\ast }})=\alpha ^{\ast }\left( T_{B}\right) $;
furthermore, $\kappa ^{\ast }=T_{B}(e)$ holds, and for $\kappa <\kappa
^{\ast }$ ($\kappa >\kappa ^{\ast }$, respectively) we have $\lim_{\rho
\rightarrow a}P_{0,1,\rho }(\Phi _{B,\kappa })=1$ ($=0$, respectively).
\end{enumerate}
\end{proposition}

Part 3 is silent on the limiting power in case $\kappa =\kappa ^{\ast }$.
Under additional assumptions, information on the limiting power in this case
has been provided in Corollary \ref{Lem_Illust_3}; we do not repeat the
results. Furthermore, note that in view of Lemma \ref{DT} in Appendix \ref%
{app_proofs} the analogon to $\kappa ^{\ast }$ in Part 1 is $\lambda
_{n-k}(B)$, whereas in Part 2 it is $\lambda _{1}(B)$.

In the case of a pure correlation model, i.e., $k=0$, the condition $e\notin 
\limfunc{span}(X)$ is always satisfied and the preceding lemma tells us that
the limiting power of the test based on $T_{B}$ is then always $1$ (for
every choice of $\kappa $ $<\lambda _{n-k}(B)$), and thus the power-trap
never arises, if and only if $e$ is an eigenvector of $B$ to the eigenvalue $%
\lambda _{n-k}(B)$.

As already noted in the discussion following Corollary \ref{Lem_Illust_1},
point-optimal invariant as well as locally best invariant tests are in
general not guaranteed to be immune to the zero-power trap phenomenon, i.e.,
they can fall under the wrath of Case 2 or 3 of the preceding proposition.
However, under its assumptions, Proposition \ref{D3new} also tells us how we
may construct -- for a given covariance model $\Sigma \left( \cdot \right) $
and a given design matrix $X$ -- a test that avoids the zero-power trap and
even has limiting power equal to $1$: All that needs to be done is to choose 
$B$ such that $C_{X}e\in \limfunc{Eig}\left( B,\lambda _{n-k}(B)\right) $
holds; one such choice is given by $B=C_{X}ee^{\prime }C_{X}^{\prime }$, but
there are many other choices. However, this result does not tell us anything
about whether or not such a test has good power properties for $\rho $ not
close to $a$. For more on ways to overcome the zero-power trap see \cite%
{Prein2014}.

\begin{remark}
\label{rem_D3new}Suppose the rejection probabilities $P_{\beta ,\sigma ,\rho
}(\Phi _{B,\kappa })$ are independent of $\beta $ and $\sigma $ (which is,
e.g., the case under Assumption \ref{ASDR} (see Remark \ref{invariance}))
and suppose $P_{0,1,0}$ is absolutely continuous w.r.t. $\mu _{\mathbb{R}%
^{n}}$. Furthermore assume that $\lambda _{1}(B)<\lambda _{n-k}(B)$ holds.
Then it follows from our Lemma \ref{DT} in Appendix \ref{app_proofs} that
for every $\alpha \in \left( 0,1\right) $ one can find a $\kappa \left(
\alpha \right) \in \left( \lambda _{1}(B),\lambda _{n-k}(B)\right) $ such
that $\Phi _{B,\kappa \left( \alpha \right) }$ has size $\alpha $. If,
additionally, $P_{0,1,0}$ has a density that is positive on an open
neighborhood of the origin except possibly for an $\mu _{\mathbb{R}^{n}}$%
-null set, then $\kappa \left( \alpha \right) $ is unique and satisfies $%
\kappa \left( \alpha \right) \rightarrow \lambda _{1}(B)$ ($\rightarrow
\lambda _{n-k}(B)$) as $\alpha \rightarrow 1$ ($\rightarrow 0$).
\end{remark}

While Proposition \ref{D3new} concerns the case $e\notin \limfunc{span}(X)$,
we have, as a simple consequence of Theorem \ref{BT}, the following result
in case $e\in \limfunc{span}(X)$. In contrast to the cases discussed in the
preceding proposition only the case $\alpha ^{\ast }\left( T_{B}\right) =0$
can occur. Recall from Remark \ref{RASCII} that whenever Assumptions \ref%
{ASC} and \ref{ASCII} both hold, the vector $e$ in the subsequent
proposition is the same as the vector $e$ in Proposition \ref{D3new} above
(since $n>2$ must hold in the subsequent proposition). Also recall that the
rejection probabilities do neither depend on $\beta $ nor $\sigma $ under
the assumption of the subsequent proposition, hence the results could be
rephrased for $P_{\beta ,\sigma ,\rho }(\Phi _{B,\kappa })$ where $\beta $
and $\sigma $ are arbitrary.

\begin{proposition}
\label{D4new} Suppose Assumptions \ref{ASDR} and \ref{ASCII} hold.
Furthermore, assume that the distribution of $\mathbf{z}$ (figuring in
Assumption \ref{ASDR}) possesses a density $p$ w.r.t. Lebesgue measure $\mu
_{\mathbb{R}^{n}}$, which is $\mu _{\mathbb{R}^{n}}$-almost everywhere
continuous and which is positive on an open neighborhood of the origin
except possibly for a $\mu _{\mathbb{R}^{n}}$-null set. Suppose $e\in 
\limfunc{span}(X)$, where $e$ is the vector figuring in Assumption \ref%
{ASCII} and suppose that $\lambda _{1}(B)<\lambda _{n-k}(B)$ holds. Then $%
\alpha ^{\ast }\left( T_{B}\right) =0$ always holds. Furthermore,%
\begin{equation*}
0<\liminf_{\rho \rightarrow a}P_{0,1,\rho }(\Phi _{B,\kappa })\leq
\limsup_{\rho \rightarrow a}P_{0,1,\rho }(\Phi _{B,\kappa })<1
\end{equation*}%
holds for every $\kappa \in (\lambda _{1}(B),\lambda _{n-k}(B))$, whereas%
\begin{equation*}
\liminf_{\rho \rightarrow a}P_{0,1,\rho }(\Phi _{B,\kappa })=1
\end{equation*}%
holds for $\kappa \leq \lambda _{1}(B)$. [For $\kappa \geq \lambda _{n-k}(B)$
we trivially always have $\Phi _{B,\kappa }=\emptyset $.]
\end{proposition}

\begin{remark}
The assumption in the preceding proposition, that $p$ is positive on an open
neighborhood of the origin except possibly for a $\mu _{\mathbb{R}^{n}}$%
-null set, can be replaced by the weaker assumption used in Part A.3 in
Theorem \ref{BT}. Furthermore, the assumption that a density $p$ exists can
be completely removed if $\mathfrak{P}$ is assumed to be an elliptically
symmetric family with $\Pr \left( \mathbf{z}=0\right) =0$.
\end{remark}

\subsection{On indistinguishability by invariant tests\label{IP}}

The discussion so far has been concerned with evaluating the power function
of a $G_{X}$-invariant test for values of $\rho $ close to $a$, the upper
bound of the range of $\rho $. In particular, we have identified conditions
under which the power function approaches zero for $\rho \rightarrow a$
(zero-power trap). These conditions, of course, depend on the test
considered as well as on the underlying model. In this section we now
isolate conditions on the model alone under which the null and alternative
hypotheses are indistinguishable by \emph{any} $G_{X}$-invariant test (in
fact, by any $G_{X}^{1}$-invariant test) whatsoever. These results, given in
Theorem \ref{ID} and Corollary \ref{IDD} below, contain a number of results
in the literature as special cases: (i) The univariate case of Theorem 5 in 
\cite{Arnold79} concerning flatness of the power function of invariant tests
in a linear model with intercept and exchangeably distributed errors, (ii)
Theorem 5 in \cite{Kadiyala}, (iii) those parts of Propositions 3-5 in \cite%
{Mart10} regarding flatness of the power function of the tests considered
there (see Section \ref{idsp} for further discussion), (iv) the first half
of the theorem proved in \cite{Mart11} (see also Section \ref{idsp}), and
(v) the result on the likelihood ratio test in \cite{Kariya80JASA}.

\begin{theorem}
\label{ID}Suppose that for some $0<\rho ^{\ast }<a$ the matrix $C_{X}\Sigma
(\rho ^{\ast })C_{X}^{\prime }$ is a multiple of $I_{n-k}$, i.e., $%
C_{X}\Sigma (\rho ^{\ast })C_{X}^{\prime }=\delta \left( \rho ^{\ast
}\right) I_{n-k}$.

\begin{enumerate}
\item Then for every $n\times n$ matrix $K(\rho ^{\ast })$ satisfying $%
K(\rho ^{\ast })K^{\prime }(\rho ^{\ast })=\Sigma (\rho ^{\ast })$ there
exists an orthogonal $n\times n$ matrix $U(\rho ^{\ast })$ such that for
every $\beta \in \mathbb{R}^{k}$ and every $0<\sigma <\infty $,%
\begin{eqnarray}
\mathcal{I}_{X}(X\beta +\sigma K(\rho ^{\ast })z) &=&\mathcal{I}_{X}(U(\rho
^{\ast })z),  \notag \\
\mathcal{I}_{X}^{+}(X\beta +\sigma K(\rho ^{\ast })z) &=&\mathcal{I}%
_{X}^{+}(U(\rho ^{\ast })z),  \label{eq0ip} \\
\mathcal{I}_{X}^{1}(X\beta +\sigma K(\rho ^{\ast })z) &=&\mathcal{I}%
_{X}^{1}(\sigma \delta ^{1/2}\left( \rho ^{\ast }\right) U(\rho ^{\ast })z) 
\notag
\end{eqnarray}%
hold for every $z\in \mathbb{R}^{n}$, where $\mathcal{I}_{X}$, $\mathcal{I}%
_{X}^{+}$, and $\mathcal{I}_{X}^{1}$ have been defined in Section \ref%
{groups}.

\item Suppose, furthermore, that $\mathfrak{P}$ is an elliptically symmetric
family. Then for every $\beta \in \mathbb{R}^{k}$ and every $0<\sigma
<\infty $ 
\begin{equation*}
P_{\beta ,\sigma ,\rho ^{\ast }}\circ \mathcal{I}_{X}=P_{\beta ,\sigma
,0}\circ \mathcal{I}_{X}=P_{0,1,0}\circ \mathcal{I}_{X},
\end{equation*}%
and 
\begin{equation*}
P_{\beta ,\sigma ,\rho ^{\ast }}\circ \mathcal{I}_{X}^{+}=P_{\beta ,\sigma
,0}\circ \mathcal{I}_{X}^{+}=P_{0,1,0}\circ \mathcal{I}_{X}^{+},
\end{equation*}%
whereas 
\begin{equation*}
P_{\beta ,\sigma ,\rho ^{\ast }}\circ \mathcal{I}_{X}^{1}=P_{0,\sigma \delta
^{1/2}\left( \rho ^{\ast }\right) ,0}\circ \mathcal{I}_{X}^{1}\text{ \ and \ 
}P_{\beta ,\sigma ,0}\circ \mathcal{I}_{X}^{1}=P_{0,\sigma ,0}\circ \mathcal{%
I}_{X}^{1}.
\end{equation*}%
The same relations hold with $\mathcal{I}_{X}$, $\mathcal{I}_{X}^{+}$, and $%
\mathcal{I}_{X}^{1}$, respectively, replaced by arbitrary $G_{X}$-, $%
G_{X}^{+}$-, or $G_{X}^{1}$-invariant statistics, meaning that no $G_{X}^{1}$%
-invariant test (and a fortiori no $G_{X}^{+}$-invariant or $G_{X}$%
-invariant test) can distinguish the null $H_{0}$ defined in (\ref%
{testproblem}) from the alternative $\rho =\rho ^{\ast }$, $\beta \in 
\mathbb{R}^{k}$, $0<\sigma <\infty $. In particular, the power function of
any $G_{X}^{+}$-invariant test (and a fortiori of any $G_{X}$-invariant
test) is constant on $\mathbb{R}^{k}\times (0,\infty )\times \left\{ 0,\rho
^{\ast }\right\} $, whereas for any $G_{X}^{1}$-invariant test power is
always less than or equal to size.
\end{enumerate}
\end{theorem}

\begin{corollary}
\label{IDD}Suppose $C_{X}\Sigma (\rho ^{\ast })C_{X}^{\prime }$ is a
multiple of $I_{n-k}$ for every $\rho ^{\ast }\in (0,a)$ and $\mathfrak{P}$
is an elliptically symmetric family. Then no $G_{X}^{1}$-invariant test (and
a fortiori no $G_{X}^{+}$-invariant or $G_{X}$-invariant test) can
distinguish $H_{0}$ from the alternative $H_{1}$ defined in (\ref%
{testproblem}). In particular, the power function of any $G_{X}^{+}$%
-invariant test (and a fortiori of any $G_{X}$-invariant test) is constant
on $\mathbb{R}^{k}\times (0,\infty )\times \lbrack 0,a)$, whereas for any $%
G_{X}^{1}$-invariant test power is always less than or equal to size.
\end{corollary}

\begin{remark}
(i) The condition that $C_{X}\Sigma (\rho ^{\ast })C_{X}^{\prime }$ is a
multiple of $I_{n-k}$ does not depend on the particular choice of $C_{X}$ as
any two such choices differ only by premultiplication with an orthogonal
matrix. Furthermore, note that the condition $C_{X}\Sigma (\rho ^{\ast
})C_{X}^{\prime }=\delta \left( \rho ^{\ast }\right) I_{n-k}$ is equivalent
to $\Pi _{\limfunc{span}(X)^{\bot }}\Sigma (\rho ^{\ast })\Pi _{\limfunc{span%
}(X)^{\bot }}=\delta \left( \rho ^{\ast }\right) \Pi _{\limfunc{span}%
(X)^{\bot }}$, see Lemma \ref{auxid} in Appendix \ref{app_proofs}.

(ii) Suppose that for some $0<\rho ^{\ast }<a$ the matrix $C_{X}\Sigma (\rho
^{\ast })C_{X}^{\prime }$ is not a multiple of $I_{n-k}$ and that $\mathfrak{%
P}$ is an elliptically symmetric family. Then it can be shown that for every 
$\alpha \in (0,1)$ there exists a $G_{X}$-invariant size $\alpha $ test with
power at $(\beta ,\sigma ,\rho ^{\ast })$ strictly larger than $\alpha $ for
every $\beta \in \mathbb{R}^{k}$ and every $0<\sigma <\infty $. As a
consequence of this result and Theorem \ref{ID} we see that the hypothesis $%
\rho =0$ and the alternative $\rho =\rho ^{\ast }$ are distinguishable by a $%
G_{X}$-invariant ($G_{X}^{+}$-invariant, $G_{X}^{1}$-invariant) test if and
only if $C_{X}\Sigma (\rho ^{\ast })C_{X}^{\prime }$ is not a multiple of $%
I_{n-k}$. [If Assumption \ref{ASDR} is satisfied but $\mathfrak{P}$ is not
an elliptically symmetric family, the hypothesis $\rho =0$ and the
alternative $\rho =\rho ^{\ast }$ may still be distinguishable by a $G_{X}$%
-invariant test even in the case where $C_{X}\Sigma (\rho ^{\ast
})C_{X}^{\prime }$ is a multiple of $I_{n-k}$, provided $L(\rho ^{\ast })$
from Assumption \ref{ASDR} gives rise to a $U(\rho ^{\ast })\neq L\left(
0\right) $ when it is used for $K(\rho ^{\ast })$ in Part 1 of Theorem \ref%
{ID}.]
\end{remark}

\begin{remark}
\label{IDD2}\emph{(Generalization of Theorem \ref{ID} and Corollary \ref{IDD}%
) }Part 2 of Theorem \ref{ID} is true more generally if $\mathfrak{P}$
satisfies Assumption \ref{ASDR} and if $\Pi _{\limfunc{span}(X)^{\bot
}}L(\rho ^{\ast })\mathbf{z}$ has the same distribution as a positive
multiple of $\Pi _{\limfunc{span}(X)^{\bot }}L(0)\mathbf{z}$, where $L\left(
\cdot \right) $ and $\mathbf{z}$ are as in Assumption \ref{ASDR} (the
multiple then being necessarily equal to $\delta ^{1/2}\left( \rho ^{\ast
}\right) $). A sufficient condition for this clearly is that $\Pi _{\limfunc{%
span}(X)^{\bot }}L(\rho ^{\ast })$ is a positive multiple of $\Pi _{\limfunc{%
span}(X)^{\bot }}L(0)$, for which in turn a sufficient condition is that
both of these two matrices are a multiple of $\Pi _{\limfunc{span}(X)^{\bot
}}$ with the multiples being non-zero and having the same sign. Similarly,
Corollary \ref{IDD} holds if $\mathfrak{P}$ satisfies Assumption \ref{ASDR}
and the distributions of $\delta ^{-1/2}\left( \rho ^{\ast }\right) \Pi _{%
\limfunc{span}(X)^{\bot }}L(\rho ^{\ast })\mathbf{z}$ for $\rho ^{\ast }\in
\lbrack 0,a)$ do not depend on $\rho ^{\ast }$ (a sufficient condition for
this being that $\Pi _{\limfunc{span}(X)^{\bot }}L(\rho ^{\ast })$ is a
positive multiple of $\Pi _{\limfunc{span}(X)^{\bot }}L(0)$ for every $\rho
^{\ast }\in \lbrack 0,a)$). Such cases arise naturally in the context of
spatial models, see Section \ref{idsp}.
\end{remark}

\begin{remark}
Suppose Assumption \ref{ASC} holds and $C_{X}\Sigma (\rho ^{\ast
})C_{X}^{\prime }=\delta \left( \rho ^{\ast }\right) I_{n-k}$ for all $\rho
^{\ast }\in (0,a)$ (or at least for a sequence $\rho _{m}^{\ast }$
converging to $a$). It is then not difficult to see that then either $e\in 
\limfunc{span}(X)$ or $n=k+1$ must hold.\footnote{%
Cf. Footnote \ref{FnA}.}
\end{remark}

Theorem \ref{ID} explains the flatness of power functions of $G_{X}$- (or $%
G_{X}^{+}$-) invariant tests observed in the literature cited above in terms
of an identification problem in the "reduced" experiment, where the
reduction is effected by the action of the group $G_{X}$ (or $G_{X}^{+}$)
(i.e., the parameters are not identifiable from the distribution of the
corresponding maximal invariant statistic); cf. Remark 2 in \cite{Mart11}
for a special case. In our framework this shows that what has been dubbed 
\textit{non-identifiability as a hypothesis} in \cite{Kariya80JASA} is
simply an identification problem in the distribution of the maximal
invariant statistic.

\section{Some generalizations\label{gen}}

\begin{remark}
\label{rem_gen_1} \emph{(Generalizations of the distributional assumptions)}
(i) We start with the following simple observation: Suppose Assumption \ref%
{ASDR} holds with $\Pr \left( \mathbf{z}=0\right) =0$. Let $\mathbf{z}^{\dag
}$ be another random vector of the same dimension as $\mathbf{z}$ (possibly
defined on another probability space) with $\Pr \left( \mathbf{z}^{\dag
}=0\right) =0$ and such that $\mathbf{z}^{\dag }/\left\Vert \mathbf{z}^{\dag
}\right\Vert $ has the same distribution as $\mathbf{z}/\left\Vert \mathbf{z}%
\right\Vert $\emph{. }It is then easy to see that the rejection
probabilities of any $G_{X}$-invariant (or $G_{X}^{+}$-invariant) test are
the same whether they are computed under $P_{\beta ,\sigma ,\rho }$ or under 
$P_{\beta ,\sigma ,\rho }^{\dag }$, where $P_{\beta ,\sigma ,\rho }^{\dag }$
is the distribution of $\mathbf{y}^{\dag }$ which is obtained from $\mathbf{z%
}^{\dag }$ via Assumption \ref{ASDR} in the same way as $\mathbf{y}$ is
obtained from $\mathbf{z}$. Hence, any result that holds for rejection
probabilities of a $G_{X}$-invariant (or $G_{X}^{+}$-invariant) test
obtained under model $\mathfrak{P}^{\dag }$ automatically carries over to
the rejection probabilities of the same test obtained under model $\mathfrak{%
P}$.

(ii) An immediate consequence of the preceding observation is, for example,
that Part 1 of Theorem \ref{BT2} continues to hold if the requirement that $%
\mathbf{z}$ has a density is replaced by the following weaker condition
(just apply Part 1 of Theorem \ref{BT2} to $\mathfrak{P}^{\dag }$):

\emph{Condition (*): }$\Pr \left( \mathbf{z}=0\right) =0$\emph{\ and there
exists a random vector }$\mathbf{z}^{\dag }$\emph{, which possesses a
density }$p^{\dag }$\emph{\ w.r.t. Lebesgue measure, such that }$\mathbf{z}%
/\left\Vert \mathbf{z}\right\Vert $\emph{\ and }$\mathbf{z}^{\dag
}/\left\Vert \mathbf{z}^{\dag }\right\Vert $\emph{\ have the same
distribution. }

This condition can be shown to be equivalent to the more explicit condition
that $\Pr \left( \mathbf{z}=0\right) =0$ and that $\mathbf{z}/\left\Vert 
\mathbf{z}\right\Vert $ possesses a density with respect to the uniform
probability measure $\upsilon _{S^{n-1}}$ on $S^{n-1}$, see Lemmata \ref%
{Proj_1} and \ref{Proj_2} in Appendix \ref{auxil}. As a consequence, Part 1
of Theorem \ref{BT2} could have been stated more generally under the
assumption that $\Pr \left( \mathbf{z}=0\right) =0$ and that $\mathbf{z}%
/\left\Vert \mathbf{z}\right\Vert $ possesses a density with respect to the
uniform probability measure $\upsilon _{S^{n-1}}$ on $S^{n-1}$.

(iii) The same reasoning as in (ii) shows that Part A of Theorem \ref{BT}
holds even without the assumption of absolute continuity of the distribution
of $\mathbf{z}$ under the following weaker assumptions: Parts A.1 and A.2
hold provided Condition (*) is satisfied and provided $\mathbf{z}^{\dag }$
can be chosen in such a way that the density $p^{\dag }$ is $\mu _{\mathbb{R}%
^{n}}$-almost everywhere continuous; an explicit sufficient condition for
this is that $\Pr \left( \mathbf{z}=0\right) =0$ holds and that $\mathbf{z}%
/\left\Vert \mathbf{z}\right\Vert $ possesses a $\upsilon _{S^{n-1}}$-almost
everywhere continuous density, see Lemma \ref{Proj_2} in Appendix \ref{auxil}%
. [Unfortunately, this explicit condition is not necessary, making it
difficult to give a simple equivalent condition which is in terms of the
distribution of $\mathbf{z}/\left\Vert \mathbf{z}\right\Vert $ only.]
Furthermore, Part A.3 holds, provided Condition (*) is satisfied and
provided $\mathbf{z}^{\dag }$ can be chosen in such a way that the density $%
p^{\dag }$ is $\mu _{\mathbb{R}^{n}}$-almost everywhere continuous and has
the property that for $\upsilon _{S^{n-1}}$-almost all $s\in S^{n-1}$ the
function $p^{\dag }\left( rs\right) $ does not vanish $\mu _{\left( 0,\infty
\right) }$-almost everywhere. An explicit sufficient condition for this is
that $\Pr \left( \mathbf{z}=0\right) =0$ and $\mathbf{z}/\left\Vert \mathbf{z%
}\right\Vert $ possesses a $\upsilon _{S^{n-1}}$-almost everywhere
continuous and $\upsilon _{S^{n-1}}$-almost everywhere positive density, see
Lemmata \ref{Proj_1} and \ref{Proj_2} in Appendix \ref{auxil}.

(iv) In case $\mathfrak{P}$ is an elliptically symmetric family with $\Pr
\left( \mathbf{z}=0\right) =0$ then $\mathbf{z}$ is spherically symmetric
entailing that the distribution of $\mathbf{z}/\left\Vert \mathbf{z}%
\right\Vert $ is the uniform distribution on the unit sphere $S^{n-1}$.
Hence, the explicit conditions discussed above are met, entailing that in
this case Condition (*) is always satisfied and $\mathbf{z}^{\dag }$ can be
chosen such that $p^{\dag }$ is $\mu _{\mathbb{R}^{n}}$-almost everywhere
continuous and $\mu _{\mathbb{R}^{n}}$-almost everywhere positive (in fact, $%
\mathbf{z}^{\dag }$ can be chosen to be Gaussian). This is what underlies
Parts 2 and 3 of Theorem \ref{BT2} as well as Part B of Theorem \ref{BT}.

(v) Suppose $\mathfrak{P}$ does not satisfy Assumption \ref{ASDR} but each
element $P_{\beta ,\sigma ,\rho }$ of $\mathfrak{P}$ is elliptically
symmetric and does not have an atom at $X\beta $ (that is, $\mathbf{u}$ is
now distributed as $\sigma L(\rho )\mathbf{w}$ where $\mathbf{w}$ has zero
mean, identity covariance matrix, and is spherically symmetric with $\Pr
\left( \mathbf{w}=0\right) =0$, but where the distribution of $\mathbf{w}$
now may depend on the parameters $\beta ,\sigma ,\rho $). Then it follows
from the results in Appendix \ref{auxil} and from the argument underlying
the discussion in (i) above that we may replace $\mathfrak{P}$ by an \emph{%
elliptically symmetric\ family }$\mathfrak{P}^{\dag }$ (even by a Gaussian
family) without affecting the rejection probabilities of $G_{X}$-invariant
(or $G_{X}^{+}$-invariant) tests and then apply our results. [More
generally, if $\mathbf{w}$ is not necessarily spherically symmetric, but the
distribution of $\mathbf{w}/\left\Vert \mathbf{w}\right\Vert $ does not
depend on the parameters $\beta ,\sigma ,\rho $, we may replace $\mathfrak{P}
$ by a family $\mathfrak{P}^{\dag }$ that is based on a $\mathbf{z}^{\dag }$%
, the distribution of which does not depend on the parameters, and
consequently satisfies Assumption \ref{ASDR}.]

(vi) In the above discussion we have so far not considered cases where
Assumption \ref{ASDR} holds, but $\vartheta :=\Pr \left( \mathbf{z}=0\right) 
$ is positive. These cases can be treated as follows: Observe that then $%
P_{\beta ,\sigma ,\rho }=\vartheta \delta _{X\beta }+\left( 1-\vartheta
\right) \tilde{P}_{\beta ,\sigma ,\rho }$ where now $\tilde{P}_{\beta
,\sigma ,\rho }$ satisfies Assumption \ref{ASDR} and the corresponding $%
\mathbf{\tilde{z}}$ has no mass at the origin (and is spherically symmetric
if $\mathbf{z}$ is so). Now for a $G_{X}$-invariant ($G_{X}^{+}$-invariant, $%
G_{X}^{1}$-invariant) test $\varphi $ the rejection probabilities satisfy $%
E_{\beta ,\sigma ,\rho }\left( \varphi \right) =\vartheta \varphi \left(
X\beta \right) +\left( 1-\vartheta \right) \tilde{E}_{\beta ,\sigma ,\rho
}\left( \varphi \right) $, where we observe that $\varphi \left( X\beta
\right) =\varphi \left( 0\right) $ is a constant not depending on $\beta $
(due to invariance of $\varphi $). Hence, the behavior of $E_{\beta ,\sigma
,\rho }\left( \varphi \right) $ can be deduced from the behavior of $\tilde{E%
}_{\beta ,\sigma ,\rho }\left( \varphi \right) $, to which our results are
applicable.
\end{remark}

\begin{remark}
\label{rem_gen_1.5}\emph{(Semiparametric Models) }Throughout the paper we
have taken a parametric viewpoint in that the distribution of $\mathbf{y}$
is assumed to be completely determined by the parameters $\beta $, $\sigma $%
, and $\rho $. The above discussion shows that some of the results of the
paper like Theorems \ref{BT} and \ref{BT2} also apply in broader
semiparametric settings (as only properties of the distribution of $\mathbf{z%
}/\left\Vert \mathbf{z}\right\Vert $ and $\Pr \left( \mathbf{z}=0\right) =0$
matter). To give just one example, let $\mathfrak{P}_{all}$ denote the \emph{%
semiparametric} model of \emph{all} elliptically symmetric distributions
with mean $X\beta $ and covariance matrix $\sigma ^{2}\Sigma \left( \rho
\right) $ that have no atom at $X\beta $ and where $\left( \beta ,\sigma
,\rho \right) $ varies in $\mathbb{R}^{k}\times (0,\infty )\times \lbrack
0,a)$. The preceding discussion then shows that the rejection probabilities
of a $G_{X}$-invariant test coincide with the rejection probabilities of a
corresponding parametric elliptically symmetric family $\mathfrak{P}$ (which
actually can be assumed be to Gaussian). Hence, the behavior of the
rejection probabilities corresponding to $\mathfrak{P}_{all}$ can
immediately be deduced from Theorems \ref{BT} and \ref{BT2} (applied to $%
\mathfrak{P}$).
\end{remark}

\begin{remark}
\label{rem_gen_2} \emph{(Extensions to }$G_{X}^{+}$\emph{-invariant tests) }%
The results of the present paper, apart from a few exceptions, are concerned
with properties of $G_{X}$-invariant tests. Concentrating on $G_{X}$%
-invariant tests, however, does not seem to impose a serious restriction
since most tests for the testing problem (\ref{testproblem}) available in
the literature satisfy this invariance property. If one nevertheless is
interested in the larger class of $G_{X}^{+}$-invariant tests, the following
observation is of interest as it allows one to extend our results to this
larger class of tests: Suppose Assumption \ref{ASDR} holds with the vector $%
\mathbf{z}$ having the same distribution as $-\mathbf{z}$ (which, in
particular, is the case under spherical symmetry). For a $G_{X}^{+}$-
invariant test $\varphi $ define the test $\varphi ^{\ast }$ by $\varphi
^{\ast }\left( y\right) =\left( \varphi \left( y\right) +\varphi \left(
-y\right) \right) /2$, which clearly is $G_{X}$-invariant. Furthermore, $%
E_{\beta ,\sigma ,\rho }\varphi =E_{\beta ,\sigma ,\rho }\varphi ^{\ast }$
holds for every $\beta $, $\sigma $, and $\rho $. Applying now our results
to $\varphi ^{\ast }$ then delivers corresponding results for $\varphi $.
\end{remark}

\begin{remark}
\emph{(Further Generalizations) }Our results easily extend to the case where
the covariance model is defined only on a set $R$, with $0\in R\subseteq
\lbrack 0,a)$, that has $a$ as its accumulation point. This observation, in
particular, allows one to obtain limiting power results along certain
sequences $\rho _{m}$, $\rho _{m}\rightarrow a$, when some of the
assumptions (like Assumptions \ref{ASC}, \ref{ASD}, \ref{ASDR}, or \ref%
{ASCII}) hold only along these sequences.
\end{remark}

\section{An application to spatial regression models\label{spatial}}

In this section we focus on regression models with spatial autoregressive
errors of order one, i.e., SAR(1) disturbances, and on spatial lag models.
First, we consider the case of a regression model with SAR(1) errors, i.e.,
what is sometimes also called a spatial error model. Second, we consider a
spatial lag model.

\subsection{Spatial error models\label{SEM}}

Let $n\geq 2$ and let $W$ be a given $n\times n$ matrix, the \textit{weights
matrix}. We assume that the diagonal elements of $W$ are all zero and that $%
W $ has a positive (real) eigenvalue, denoted by $\lambda _{\max }$, such
that any other real or complex zero of the characteristic polynomial of $W$
is in absolute value not larger than $\lambda _{\max }$. We also assume that 
$\lambda _{\max }$ has algebraic multiplicity (and thus also geometric
multiplicity) equal to $1$. Choose $f_{\max }$ as a normalized eigenvector
of $W$ corresponding to $\lambda _{\max }$ (which is unique up to
multiplication by $-1$). The spatial error model (SEM) is then given by the
regression model in equation (\ref{linmod}) where the disturbance vector $%
\mathbf{u}$ is SAR(1), i.e., for given $\beta \in \mathbb{R}^{k}$, $0<\sigma
<\infty $, and $\rho \in \lbrack 0,\lambda _{\max }^{-1})$ we have%
\begin{equation}
\mathbf{u}=\rho W\mathbf{u}+\sigma \mathbf{\varepsilon }  \label{SAR}
\end{equation}%
where $\mathbf{\varepsilon }$ is a mean zero random vector with covariance
matrix $I_{n}$. Observe that then clearly 
\begin{equation}
\mathbf{u}=(I_{n}-\rho W)^{-1}\sigma \mathbf{\varepsilon }
\label{SAR_explicit}
\end{equation}%
holds and that the covariance matrix of $\mathbf{u}$ is given by $\sigma
^{2}\Sigma _{SEM}(\rho )$ where $\Sigma _{SEM}(\rho )=[(I_{n}-\rho W^{\prime
})(I_{n}-\rho W)]^{-1}$ for $\rho \in \lbrack 0,a)$ where here $a=\lambda
_{\max }^{-1}$. Additionally we assume that the distribution of $\mathbf{%
\varepsilon }$ is a fixed distribution independent of $\beta $, $\sigma $,
and $\rho $.\footnote{%
It appears that it is implicitly assumed in \cite{Mart10} that $\mathbf{%
\varepsilon }$ is a random vector whose distribution is independent of $%
\beta $, $\sigma $, and $\rho $, cf. \cite{Mart10}, p. 155. As discussed in
Remark \ref{modelM10}(ii), it is also implicitly assumed in \cite{Mart10}
that the distribution of $\sigma ^{-1}\Sigma _{SEM}^{-1/2}(\rho )\mathbf{u}$
is independent of $\beta $, $\sigma $, and $\rho $. Note that the latter
random vector is connected to $\mathbf{\varepsilon }$ via multiplication by
an orthogonal matrix $U\left( \rho \right) $, say. If $W$ is symmetric, $%
U\left( \rho \right) \equiv I_{n}$ holds and hence both implicit assumptions
are equivalent. However, for nonsymmetric $W$, these two implicit
assumptions will typically be compatible only if the distribution of $%
\mathbf{\varepsilon }$ is spherically symmetric.} \emph{The above are the
maintained assumptions for the SEM considered in this section.} The
parametric family $\mathfrak{P}$ of probability measures induced by (\ref%
{linmod}) and (\ref{SAR}) under the maintained assumptions will be denoted
by $\mathfrak{P}_{SEM}$.

\begin{remark}
If $W$ is an (elementwise) nonnegative and irreducible matrix with zero
elements on the main diagonal, a frequent assumption for spatial weights
matrices, then\ the above assumptions on $W$ are satisfied by the
Perron-Frobenius theorem and $\lambda _{\max }$ is then the Perron-Frobenius
root of $W$ (see, e.g., \cite{HJ1985}, Theorem 8.4.4, p. 508). In this case
one can always choose $f_{\max }$ to be entrywise positive.
\end{remark}

The next lemma shows identifiability of the parameters in the model,
identifiability of $\beta $ being trivial. An immediate consequence is that
the two subsets of $\mathfrak{P}_{SEM}$ corresponding to the null hypothesis 
$\rho =0$ and alternative hypothesis $\rho >0$ are disjoint.\footnote{%
Lemmata \ref{ISAR1} and \ref{SARL} actually hold without the additional
assumption on the distribution of $\mathbf{\varepsilon }$ made above.}

\begin{lemma}
\label{ISAR1} If $\sigma _{1}^{2}\Sigma _{SEM}(\rho _{1})=\sigma
_{2}^{2}\Sigma _{SEM}(\rho _{2})$ holds for $\rho _{i}\in \lbrack 0,\lambda
_{\max }^{-1})$ and $0<\sigma _{i}<\infty $ ($i=1,2$) then $\rho _{1}=\rho
_{2}$ and $\sigma _{1}=\sigma _{2}$.
\end{lemma}

We next verify that the spatial error model satisfies Assumptions \ref{ASC}, %
\ref{ASDR}, and \ref{ASCII}, and that it satisfies Assumption \ref{ASD}
under a mild condition on the distribution of $\mathbf{\varepsilon }$. The
first claim in Lemma \ref{SARL} also appears in \cite{Mart11}, Lemma 3.3.

\begin{lemma}
\label{SARL} $\Sigma _{SEM}(\cdot )$ satisfies Assumption \ref{ASC} with $%
e=f_{\max }$ as well as Assumption \ref{ASCII} with $e=f_{\max }$, $c(\rho
)=1$, $L_{\ast }(\rho )=(I_{n}-\rho W)^{-1}$, and $\Lambda =\left(
I_{n}-\lambda _{\max }^{-1}\Pi _{\limfunc{span}(f_{\max })^{\bot }}W\right)
^{-1}-\Pi _{\limfunc{span}(f_{\max })}$.
\end{lemma}

\begin{lemma}
\label{SARL_2} $\mathfrak{P}_{SEM}$ satisfies Assumption \ref{ASDR} with $%
L\left( \rho \right) =(I_{n}-\rho W)^{-1}$ and $\mathbf{z}$ a random vector
distributed like $\mathbf{\varepsilon }$. Furthermore, if the distribution
of $\mathbf{\varepsilon }$ is absolutely continuous w.r.t. $\mu _{\mathbb{R}%
^{n}}$, or, more generally, if $\Pr (\mathbf{\varepsilon }=0)=0$ and the
distribution of $\mathbf{\varepsilon }/\left\Vert \mathbf{\varepsilon }%
\right\Vert $ is absolutely continuous w.r.t. the uniform distribution $%
\upsilon _{S^{n-1}}$ on the unit sphere $S^{n-1}$, then $\mathfrak{P}_{SEM}$
satisfies Assumption \ref{ASD}.
\end{lemma}

Given the preceding two lemmata the main results of Section \ref{main},
i.e., Theorems \ref{MT}, \ref{BT}, and \ref{BT2}, can be immediately applied
to obtain results for the spatial error model. Rather than spelling out
these general results, we provide the following two corollaries for the
purpose of illustration and thus do not strive for the weakest conditions.
These corollaries provide, in particular, correct versions of the claims in
Corollary 1 in \cite{Mart10}. Recall that by the assumed $G_{X}$-invariance
the rejection probabilities $E_{\beta ,\sigma ,\rho }\left( \varphi \right) $
in the subsequent results do in fact neither depend on $\beta $ nor $\sigma $%
, cf. Remark \ref{invariance}.

\begin{corollary}
\label{Cor1new} Given the maintained assumptions for the SEM suppose
furthermore that either (i) the distribution of $\mathbf{\varepsilon }$
possesses a $\mu _{\mathbb{R}^{n}}$-density $p$ that is continuous $\mu _{%
\mathbb{R}^{n}}$-almost everywhere and that is positive on an open
neighborhood of the origin except possibly for a $\mu _{\mathbb{R}^{n}}$%
-null set, or (ii) the distribution of $\mathbf{\varepsilon }$ is
spherically symmetric with no atom at the origin. Then for every $G_{X}$%
-invariant test $\varphi $ the following statements hold:

\begin{enumerate}
\item If $\varphi $ is continuous at $f_{\max }$ then for every $\beta \in 
\mathbb{R}^{k}$, $0<\sigma <\infty $, we have $E_{\beta ,\sigma ,\rho
}\left( \varphi \right) \rightarrow \varphi \left( f_{\max }\right) $ for $%
\rho \rightarrow \lambda _{\max }^{-1}$, $\rho \in \lbrack 0,\lambda _{\max
}^{-1})$.

\item Suppose $\varphi $ satisfies $\varphi \left( y\right) =\varphi \left(
y+f_{\max }\right) $ for every $y\in \mathbb{R}^{k}$ (which is certainly the
case if $f_{\max }\in \limfunc{span}\left( X\right) $). Then for every $%
\beta \in \mathbb{R}^{k}$, $0<\sigma <\infty $, we have $E_{\beta ,\sigma
,\rho }\left( \varphi \right) \rightarrow E\varphi \left( \Lambda \mathbf{%
\varepsilon }\right) $ for $\rho \rightarrow \lambda _{\max }^{-1}$, $\rho
\in \lbrack 0,\lambda _{\max }^{-1})$. The limit $E\varphi \left( \Lambda 
\mathbf{\varepsilon }\right) $ is strictly between $0$ and $1$ provided
neither $\varphi =0$ $\mu _{\mathbb{R}^{n}}$-almost everywhere nor $\varphi
=1$ $\mu _{\mathbb{R}^{n}}$-almost everywhere holds. [The matrix $\Lambda $
is defined in Lemma \ref{SARL}.]

\item If $\varphi $ is the indicator function of a critical region $\Phi $,
we have for every $\beta \in \mathbb{R}^{k}$, $0<\sigma <\infty $, and as $%
\rho \rightarrow \lambda _{\max }^{-1}$, $\rho \in \lbrack 0,\lambda _{\max
}^{-1})$:

\begin{itemize}
\item $f_{\max }\in \limfunc{int}(\Phi )$ implies $P_{\beta ,\sigma ,\rho
}(\Phi )\rightarrow 1$.

\item $f_{\max }\notin \limfunc{cl}(\Phi )$ implies $P_{\beta ,\sigma ,\rho
}(\Phi )\rightarrow 0$.

\item $f_{\max }\in \limfunc{span}\left( X\right) $ implies $P_{\beta
,\sigma ,\rho }(\Phi )\rightarrow \Pr \left( \Lambda \mathbf{\varepsilon }%
\in \Phi \right) $. The limiting probability is strictly between $0$ and $1$
provided neither $\Phi $ nor its complement are $\mu _{\mathbb{R}^{n}}$-null
sets.
\end{itemize}

\item If $\varphi $ is the indicator function of the critical region $\Phi
_{B,\kappa }$ given by (\ref{quadratic}) with $B$ satisfying $\lambda
_{1}\left( B\right) <\lambda _{n-k}\left( B\right) $ and with $\kappa \in
\lbrack \lambda _{1}\left( B\right) ,\lambda _{n-k}\left( B\right) )$, then
we have for every $\beta \in \mathbb{R}^{k}$, $0<\sigma <\infty $, and as $%
\rho \rightarrow \lambda _{\max }^{-1}$, $\rho \in \lbrack 0,\lambda _{\max
}^{-1})$:

\begin{itemize}
\item $T_{B}\left( f_{\max }\right) >\kappa $ implies $P_{\beta ,\sigma
,\rho }(\Phi _{B,\kappa })\rightarrow 1$.\footnote{%
Note that $T_{B}\left( f_{\max }\right) >\kappa $ entails $f_{\max }\notin 
\limfunc{span}\left( X\right) $ in view of (\ref{T_quadratic}) and $\kappa
\geq \lambda _{1}\left( B\right) $.}

\item $T_{B}\left( f_{\max }\right) <\kappa $ and $f_{\max }\notin \limfunc{%
span}\left( X\right) $ implies $P_{\beta ,\sigma ,\rho }(\Phi _{B,\kappa
})\rightarrow 0$.

\item $f_{\max }\in \limfunc{span}\left( X\right) $ implies $P_{\beta
,\sigma ,\rho }(\Phi _{B,\kappa })\rightarrow \Pr \left( \Lambda \mathbf{%
\varepsilon }\in \Phi _{B,\kappa }\right) $. The limiting probability is
strictly between $0$ and $1$ provided $\kappa \in (\lambda _{1}\left(
B\right) ,\lambda _{n-k}\left( B\right) )$, while it is $1$ for $\kappa
=\lambda _{1}\left( B\right) $.
\end{itemize}
\end{enumerate}
\end{corollary}

\begin{remark}
If $\varphi =0$ ($=1$) $\mu _{\mathbb{R}^{n}}$-almost everywhere in Part 2
or in the last claim of Part 3 of the preceding corollary, then $E_{\beta
,\sigma ,\rho }(\varphi )=0$ (or $=1$) holds for all $\beta $, $\sigma $,
and $\rho $, and hence the same holds a fortiori for the accumulation
points, see Remark \ref{RBT}(iv).
\end{remark}

Parts 3 and 4 of the preceding corollary are silent on the case $f_{\max
}\in \limfunc{bd}\left( \Phi \right) \backslash \limfunc{span}\left(
X\right) $ (recall that $\limfunc{span}\left( X\right) \subseteq \limfunc{bd}%
\left( \Phi \right) $ holds provided $\emptyset \neq \Phi \neq \mathbb{R}%
^{n} $). The next corollary provides such a result for the important
critical regions $\Phi _{B,\kappa }$ under an elliptical symmetry assumption
on $\mathfrak{P}_{SEM}$ and under the assumption of a symmetric weights
matrix $W $. More general results without the symmetry assumption on $W$,
without the elliptical symmetry assumption, and for more general classes of
tests can of course be obtained from Theorem \ref{BT2}.

\begin{corollary}
\label{Cor2new} Given the maintained assumptions for the SEM suppose
furthermore that the distribution of $\mathbf{\varepsilon }$ is spherically
symmetric with no atom at the origin and that $W$ is symmetric. Let the
critical region $\Phi _{B,\kappa }$ be given by (\ref{quadratic}). Assume $%
f_{\max }\in \limfunc{bd}\left( \Phi _{B,\kappa }\right) \backslash \limfunc{%
span}\left( X\right) $ (i.e., $f_{\max }\notin \limfunc{span}\left( X\right) 
$ and $\kappa =T_{B}\left( f_{\max }\right) \in \lbrack \lambda _{1}\left(
B\right) ,\lambda _{n-k}\left( B\right) )$ with $\lambda _{1}\left( B\right)
<\lambda _{n-k}\left( B\right) $ hold).

\begin{enumerate}
\item Suppose $C_{X}f_{\max }$ is an eigenvector of $B$ with eigenvalue $%
\lambda $. Then $\lambda =T_{B}\left( f_{\max }\right) =\kappa $ and 
\begin{equation}
P_{\beta ,\sigma ,\rho }(\Phi _{B,\kappa })\rightarrow \Pr \left( \mathbf{G}%
^{\prime }\Lambda ^{\prime }\left( C_{X}^{\prime }BC_{X}-\lambda
C_{X}^{\prime }C_{X}\right) \Lambda \mathbf{G}>0\right)  \label{li}
\end{equation}%
for $\rho \rightarrow \lambda _{\max }^{-1}$, $\rho \in \lbrack 0,\lambda
_{\max }^{-1})$, and for every $\beta \in \mathbb{R}^{k}$, $0<\sigma <\infty 
$, where $\mathbf{G}$ is a multivariate Gaussian random vector with mean
zero and covariance matrix $I_{n}$. The limit in (\ref{li}) is strictly
between $0$ and $1$ if $\lambda >\lambda _{1}(B)$, whereas it equals $1$ in
case $\lambda =\lambda _{1}(B)$.

\item Suppose $C_{X}f_{\max }$ is not an eigenvector of $B$. Then $P_{\beta
,\sigma ,\rho }(\Phi _{B,\kappa })\rightarrow 1/2$ for $\rho \rightarrow
\lambda _{\max }^{-1}$, $\rho \in \lbrack 0,\lambda _{\max }^{-1})$, and for
every $\beta \in \mathbb{R}^{k}$, $0<\sigma <\infty $.
\end{enumerate}
\end{corollary}

\begin{remark}
\label{kraemer}\emph{(Some comments on \cite{kramer2005})} (i) \cite%
{kramer2005} considers "test statistics" of the form $\mathbf{u}^{\prime
}Q_{1}\mathbf{u}/\mathbf{u}^{\prime }Q_{2}\mathbf{u}$ for general matrices $%
Q_{1}$ and $Q_{2}$. However, this ratio will then in general not be
observable and thus will not be a test statistic. Fortunately, the problem
disappears in the leading cases where $Q_{1}$ and $Q_{2}$ are such that $%
\mathbf{u}^{\prime }Q_{i}\mathbf{u}=\mathbf{y}^{\prime }Q_{i}\mathbf{y}$.
The same problem also appears in \cite{KramerZeisel1990} and \cite{small1993}%
.

(ii) The proof of the last claim in Theorem 1 of \cite{kramer2005} is in
error, as -- contrary to the claim in \cite{kramer2005} -- the quantity $%
d_{1}$ need not be strictly positive. This has already be noted by \cite%
{Mart12}, Footnote 5.

(iii) Theorem 2 in \cite{kramer2005} is not a theorem in the mathematical
sense, as it is not made precise what it means that the limiting power "is
in general strictly between $0$ and $1$".
\end{remark}

As discussed earlier, point-optimal invariant and locally best invariant
tests are in general not immune to the zero-power trap. The next result,
which is a correct version of Proposition 1 in \cite{Mart10}, now provides a
necessary and sufficient condition for the Cliff-Ord test (i.e., $%
B=W+W^{\prime }$) and a point-optimal invariant test (i.e., $B=-\Sigma
_{SEM}^{-1}(\bar{\rho})$)\ in a pure SAR-model (i.e., $k=0$) to have
limiting power equal to $1$ for every choice of the critical value $\kappa $
(excluding trivial cases). For a discussion of the problems with Proposition
1 in \cite{Mart10} see Appendix \ref{Prop1}. In the subsequent proposition
we always have $\lambda _{1}\left( B\right) <\lambda _{n}\left( B\right) $
as a consequence of the assumptions. We also note that the condition $\kappa
\in \left( \lambda _{1}\left( B\right) ,\lambda _{n}\left( B\right) \right) $
in this proposition precisely corresponds to the condition that the test has
size strictly between zero and one, cf. Remark \ref{rem_D3new}. Furthermore,
observe that while the statement that the limiting power (as $\rho
\rightarrow \lambda _{\max }^{-1}$) equals $1$ for every $\kappa \in \left(
\lambda _{1}\left( B\right) ,\lambda _{n}\left( B\right) \right) $ is in
general clearly stronger than the statement that $\alpha ^{\ast }\left(
T_{B}\right) =0$, Proposition \ref{D3new} shows that these statements are in
fact equivalent in the context of the following result. Finally, recall that
in view of invariance and the maintained assumptions of this section the
rejection probabilities do neither depend on $\beta $ nor $\sigma $.

\begin{proposition}
\label{prop1new} Given the maintained assumptions for the SEM, suppose that
the distribution of $\mathbf{\varepsilon }$ is absolutely continuous w.r.t. $%
\mu _{\mathbb{R}^{n}}$ with a density that is positive on an open
neighborhood of the origin except possibly for a $\mu _{\mathbb{R}^{n}}$%
-null set. Furthermore, assume that $k=0$. Let (i) $B=W+W^{\prime }$ or (ii) 
$B=-\Sigma _{SEM}^{-1}(\bar{\rho})$ for some $0<\bar{\rho}<\lambda _{\max
}^{-1}$. Consider the rejection region $\Phi _{B,\kappa }$ given by (\ref%
{quadratic}) with $C_{0}=I_{n}$. Then for every $\beta \in \mathbb{R}^{k}$
and $0<\sigma <\infty $ we have in both cases (i) and (ii): $P_{\beta
,\sigma ,\rho }(\Phi _{B,\kappa })\rightarrow 1$ for every $\kappa \in
\left( \lambda _{1}\left( B\right) ,\lambda _{n}\left( B\right) \right) $ as 
$\rho \rightarrow \lambda _{\max }^{-1}$, $\rho \in \lbrack 0,\lambda _{\max
}^{-1})$, if and only if $f_{\max }\in \limfunc{Eig}\left( B,\lambda
_{n}(B)\right) $. In particular, if $W$ is (elementwise) nonnegative and
irreducible, then, for both choices of $B$, $f_{\max }\in \limfunc{Eig}%
\left( B,\lambda _{n}(B)\right) $ is equivalent to $f_{\max }$ being an
eigenvector of $W^{\prime }$.
\end{proposition}

The next proposition is a correct version of Lemma E.4 in \cite{Mart10}; see
Appendix \ref{Prop1} for a discussion of the shortcomings of that lemma. It
provides conditions under which the Cliff-Ord test and point-optimal
invariant tests in a SEM with exogenous variables are not subject to the
zero-power trap and even have limiting power equal to $1$.

\begin{proposition}
\label{E4new} Given the maintained assumptions for the SEM, suppose that the
distribution of $\mathbf{\varepsilon }$ is absolutely continuous w.r.t. $\mu
_{\mathbb{R}^{n}}$ with a density that is positive on an open neighborhood
of the origin except possibly for a $\mu _{\mathbb{R}^{n}}$-null set.
Suppose further that $f_{\max }\notin \limfunc{span}(X)$, that $\limfunc{Eig}%
\left( C_{X}\Sigma _{SEM}(\rho )C_{X}^{\prime },\lambda _{n-k}(C_{X}\Sigma
_{SEM}(\rho )C_{X}^{\prime })\right) $ is independent of $0<\rho <\lambda
_{\max }^{-1}$, and that $n-k>1$. Let (i) $B=C_{X}\left( W+W^{\prime
}\right) C_{X}^{\prime }$ and suppose that $\lambda _{1}\left( B\right)
<\lambda _{n-k}\left( B\right) $, or (ii) $B=-\left( C_{X}\Sigma _{SEM}(\bar{%
\rho})C_{X}^{\prime }\right) ^{-1}$ for some $0<\bar{\rho}<\lambda _{\max
}^{-1}$. Consider the rejection region $\Phi _{B,\kappa }$ given by (\ref%
{quadratic}). Then for every $\beta \in \mathbb{R}^{k}$ and $0<\sigma
<\infty $ we have in both cases (i) and (ii): $P_{\beta ,\sigma ,\rho }(\Phi
_{B,\kappa })\rightarrow 1$ for every $\kappa \in \left( \lambda _{1}\left(
B\right) ,\lambda _{n-k}\left( B\right) \right) $ as $\rho \rightarrow
\lambda _{\max }^{-1}$, $\rho \in \lbrack 0,\lambda _{\max }^{-1})$.
\end{proposition}

\begin{remark}
\label{RE4} (i) The condition that $\limfunc{Eig}\left( C_{X}\Sigma
_{SEM}(\rho )C_{X}^{\prime },\lambda _{n-k}(C_{X}\Sigma _{SEM}(\rho
)C_{X}^{\prime })\right) $ is independent of $\rho $ is easily seen to be
satisfied, e.g., if $W$ is symmetric and if $f_{\max }\in \limfunc{span}%
(X)^{\bot }$ (and thus, in particular, if $k=0$).

(ii) If in the preceding proposition $W$ is symmetric and $f_{\max }\in 
\limfunc{span}(X)^{\bot }$ holds, then the condition $\lambda _{1}\left(
B\right) <\lambda _{n-k}\left( B\right) $ in case $B=C_{X}(W+W^{\prime
})C_{X}^{\prime }=2C_{X}WC_{X}^{\prime }$ is automatically satisfied. This
can be seen as follows: Since $f_{\max }\in \limfunc{span}(X)^{\bot }$ we
can represent $f_{\max }$ as $C_{X}^{\prime }\gamma $ for some $\gamma \in 
\mathbb{R}^{n-k}$ with $\gamma ^{\prime }\gamma =1$. On the one hand, the
largest eigenvalue of $2C_{X}WC_{X}^{\prime }$, as the maximum of $2\delta
^{\prime }C_{X}WC_{X}^{\prime }\delta $ over all normalized vectors $\delta
\in \mathbb{R}^{n-k}$, is therefore not less than $2\gamma ^{\prime
}C_{X}WC_{X}^{\prime }\gamma =2f_{\max }^{\prime }Wf_{\max }=2\lambda _{\max
}$. On the other hand, noting that $\left\Vert C_{X}^{\prime }\delta
\right\Vert =\left\Vert \delta \right\Vert $, the maximum of $2\delta
^{\prime }C_{X}WC_{X}^{\prime }\delta $ over all normalized vectors $\delta
\in \mathbb{R}^{n-k}$ is not larger than the maximum of $v^{\prime }Wv$ over
all normalized vectors $v\in \mathbb{R}^{n}$, which shows that the largest
eigenvalue of $2C_{X}WC_{X}^{\prime }$ is equal to $2\lambda _{\max }$.
Because $\lambda _{\max }$ as the largest eigenvalue of $W$ has algebraic
multiplicity $1$ by the assumptions in this section and since $%
2C_{X}WC_{X}^{\prime }$ is symmetric, we see that the algebraic multiplicity
of $2\lambda _{\max }$ as an eigenvalue of $2C_{X}WC_{X}^{\prime }$ must
also be $1$. But then $\lambda _{1}\left( B\right) <\lambda _{n-k}\left(
B\right) $ follows since $n-k>1$ has been assumed in the proposition.

(iii) If $n-k=1$ or if $n-k>1$, but $\lambda _{1}\left( B\right) =\lambda
_{n-k}\left( B\right) $ holds for $B=C_{X}\left( W+W^{\prime }\right)
C_{X}^{\prime }$, then the test statistic degenerates to a constant (and the
proposition trivially holds as $\left( \lambda _{1}\left( B\right) ,\lambda
_{n-k}\left( B\right) \right) $ is then empty).
\end{remark}

\subsection{Spatial lag models\label{SPLM}}

Let $X$ be as in Section \ref{Framework}, let $W$ be as in Section \ref{SEM}%
, and consider the spatial lag model (SLM) of the form%
\begin{equation}
\mathbf{y}=\rho W\mathbf{y}+X\beta +\sigma \mathbf{\varepsilon },
\label{SL1}
\end{equation}%
where $\beta \in \mathbb{R}^{k}$, $\rho \in \lbrack 0,\lambda _{\max }^{-1})$%
, and $0<\sigma <\infty $, and where $\mathbf{\varepsilon }$ is a mean zero
random vector with covariance matrix $I_{n}$. As in Section \ref{SEM}, we
assume that the distribution of $\mathbf{\varepsilon }$ is a fixed
distribution independent of $\beta \in \mathbb{R}^{k}$, $\sigma \in \left(
0,\infty \right) $, and $\rho \in \lbrack 0,\lambda _{\max }^{-1})$. \emph{%
The above are the maintained assumptions for the SLM considered in this
section. }Because the SLM and the SEM have the same covariance structure, a
simple consequence of Lemma \ref{ISAR1} is that also the parameters of the
SLM are identifiable. For $\rho \in \lbrack 0,\lambda _{\max }^{-1})$ we can
rewrite the above equation as 
\begin{equation}
\mathbf{y}=(I_{n}-\rho W)^{-1}(X\beta +\sigma \mathbf{\varepsilon }).
\label{SL2}
\end{equation}%
Obviously, in case $k=0$ the spatial lag model of order one coincides with
the SAR(1) model. For $k>0$, however, the SLM does \emph{not }fit into the
general framework of Section \ref{main} of the present paper. In particular,
while the problem of testing $\rho =0$ versus $\rho >0$ is still invariant
under the group $G_{0}$, it is typically no longer invariant under the
larger group $G_{X}$. Nevertheless we can establish the following result
which is similar in spirit to Theorem \ref{MT}. In the following result, $%
P_{\beta ,\sigma ,\rho }$ denotes the distribution of $\mathbf{y}$ given by (%
\ref{SL2}) under the parameters $\beta \in \mathbb{R}^{k}$, $0<\sigma
<\infty $, and $\rho \in \lbrack 0,\lambda _{\max }^{-1})$ and $E_{\beta
,\sigma ,\rho }$ denotes the corresponding expectation operator.

\begin{theorem}
\label{theorem_SPLM}Given the maintained assumptions for the SLM, assume
furthermore that the distribution of $\mathbf{\varepsilon }$ does not put
positive mass on a proper affine subspace of $\mathbb{R}^{n}$. Let $\varphi $
be a $G_{0}$-invariant test. If $\varphi $ is continuous at $f_{\max }$ then
for every $\beta \in \mathbb{R}^{k}$ and $0<\sigma <\infty $ we have $%
E_{\beta ,\sigma ,\rho }\left( \varphi \right) \rightarrow \varphi \left(
f_{\max }\right) $ for $\rho \rightarrow \lambda _{\max }^{-1}$, $\rho \in
\lbrack 0,\lambda _{\max }^{-1})$. In particular, if $\varphi $ is the
indicator function of a critical region $\Phi $ we have for every $\beta \in 
\mathbb{R}^{k}$, $0<\sigma <\infty $, and as $\rho \rightarrow \lambda
_{\max }^{-1}$, $\rho \in \lbrack 0,\lambda _{\max }^{-1})$:

\begin{itemize}
\item $f_{\max }\in \limfunc{int}(\Phi )$ implies $P_{\beta ,\sigma ,\rho
}(\Phi )\rightarrow 1$.

\item $f_{\max }\notin \limfunc{cl}(\Phi )$ implies $P_{\beta ,\sigma ,\rho
}(\Phi )\rightarrow 0$.
\end{itemize}
\end{theorem}

The above result provides a correct version of the first and third claim in
Proposition 2 in \cite{Mart10}, the proofs of which in \cite{Mart10} suffer
from the same problems as the proofs of the corresponding parts of MT1. The
second claim in Proposition 2 in \cite{Mart10} is incorrect for the same
reasons as is the second part of MT1. While Theorems \ref{BT} and \ref{BT2}
provide correct versions of the second claim of MT1, these results can not
directly be used in the context of the SLM as this model does not fit into
the framework of Section \ref{main} as noted above. We do not investigate
this issue any further here.\footnote{%
Under additional restrictive assumptions (such as $\limfunc{span}(\left(
I_{n}-\rho W\right) ^{-1}X)\subseteq \limfunc{span}\left( X\right) $ for
every $\rho \in \lbrack 0,\lambda _{\max }^{-1})$) invariance w.r.t. $G_{X}$
can again become an appropriate assumption on a test statistic and a version
of Theorem \ref{BT} can then be produced. We abstain from pursuing this any
further.}

\subsection{Indistinguishability by invariant tests in spatial regression
models \label{idsp}}

We close our discussion of spatial regression models by applying the results
on indistinguishability developed in Section \ref{IP} to these models. It
turns out that a number of results in \cite{Mart10} (namely all parts of
Proposition 3, 4, and 5 that are based on degeneracy of the test statistic)
as well as the first part of the theorem in \cite{Mart11} are consequences
of an identification problem in the distribution of the maximal invariant
statistic (more precisely, an identification problem in the "reduced"
experiment). Theorem \ref{ID} and Corollary \ref{IDD} thus provide a simple
and systematic way to recognize when this identification problem occurs.

Consider the SEM with the maintained assumptions of Section \ref{SEM} and
additionally assume \emph{for this paragraph only} that the distribution of
the error $\mathbf{\varepsilon }$ is spherically symmetric. As shown in
Section \ref{IP}, the condition for the identification problem in the
reduced experiment to occur, entailing a constant power function for any $%
G_{X}$-invariant (even for any $G_{X}^{+}$-invariant) test, is then that $%
C_{X}\Sigma _{SEM}(\rho )C_{X}^{\prime }$ is a multiple of $I_{n-k}$ for
every $\rho \in (0,\lambda _{\max }^{-1})$. As can be seen from Lemma \ref%
{auxid} in Appendix \ref{app_proofs}, a sufficient condition for this is
that $\limfunc{span}(X)^{\bot }$ is contained in an eigenspace of $\Sigma
_{SEM}(\rho )$ for every $\rho \in (0,\lambda _{\max }^{-1})$, a condition
that appears in Proposition 3 of \cite{Mart10}, which is a statement about
point-optimal invariant and locally best invariant tests. Thus the
corresponding part of this proposition is an immediate consequence of
Corollary \ref{IDD}; moreover, and in contrast to this proposition in \cite%
{Mart10}, it now follows that this result holds more generally for \emph{any}
$G_{X}$-invariant (even any $G_{X}^{+}$-invariant) test and that the
Gaussianity assumption in this proposition can be weakened to elliptical
symmetry. In a similar way, Propositions 4 and 5 in \cite{Mart10} make use
of the conditions that $W$ is symmetric and $\limfunc{span}$$(X)^{\bot }$ is
contained in an eigenspace of $W$. In the subsequent lemma we show that the
condition that $\limfunc{span}$$(X)^{\bot }$ is contained in an eigenspace
of $W^{\prime }$ is already sufficient for $C_{X}\Sigma _{SEM}(\rho
)C_{X}^{\prime }$ to be a multiple of $I_{n-k}$ for every $\rho \in
(0,\lambda _{\max }^{-1})$. Thus the subsequent lemma combined with
Corollary \ref{IDD} establishes, in particular, the respective parts of
Propositions 4 and 5 in \cite{Mart10}. The preceding comments are of some
importance as there are several problems with Propositions 3, 4, and 5 in 
\cite{Mart10} which are discussed in Appendix \ref{P3}.

\begin{lemma}
\label{SEMID} Let $W$ be a weights matrix as in Section \ref{SEM} and let $X$
be an $n\times k$ matrix ($n>k$) such that 
\begin{equation}
\text{$\limfunc{span}$}(X)^{\bot }\subseteq \limfunc{Eig}\left( W^{\prime
},\lambda \right)  \label{eig}
\end{equation}%
is satisfied for some eigenvalue $\lambda \in \mathbb{R}$ of $W^{\prime }$.
Then%
\begin{equation*}
\Pi _{\limfunc{span}(X)^{\bot }}\left( I_{n}-\rho W\right) ^{-1}=\left(
1-\rho \lambda \right) ^{-1}\Pi _{\limfunc{span}(X)^{\bot }}
\end{equation*}%
and%
\begin{equation*}
C_{X}\Sigma _{SEM}(\rho )C_{X}^{\prime }=\left( 1-\rho \lambda \right)
^{-2}I_{n-k}
\end{equation*}%
hold for every $0\leq \rho <\lambda _{\max }^{-1}$.
\end{lemma}

In the following example we show that the first half of the theorem in \cite%
{Mart11} is a special case of Remark \ref{IDD2} following Corollary \ref{IDD}
combined with the preceding lemma.

\begin{example}
(i) Consider the SEM with the maintained assumptions of Section \ref{SEM}.
Suppose that $W$ is an $n\times n$ ($n\geq 2$) equal weights matrix, i.e., $%
w_{ij}$ is constant for $i\neq j$ and zero else and that $\limfunc{span}(X)$
contains the intercept. Without loss of generality we assume $w_{ij}=1$ for $%
i\neq j$. Clearly, $W$ is symmetric and has the eigenvalues $\lambda
_{1}\left( W\right) =\ldots =\lambda _{n-1}\left( W\right) =-1$ and $\lambda
_{\max }=\lambda _{n}\left( W\right) =n-1$. The eigenspace corresponding to $%
\lambda _{\max }$ is spanned by the eigenvector $f_{\max }=n^{-1/2}(1,\ldots
,1)^{\prime }$ and the other eigenspace consists of all vectors orthogonal
to $f_{\max }$. Since every element of $\limfunc{span}(X)^{\bot }$ is
orthogonal to $(1,\ldots ,1)^{\prime }$ we have 
\begin{equation*}
\text{$\limfunc{span}$}(X)^{\bot }\subseteq \limfunc{Eig}(W,-1)=\limfunc{Eig}%
(W^{\prime },-1).
\end{equation*}%
Therefore, by Lemma \ref{SEMID} together with Remark \ref{IDD2}, the power
function of every $G_{X}^{+}$-invariant test must be constant.

(ii)\ Consider next the SLM with the maintained assumptions of Section \ref%
{SPLM} with the same weights matrix and the same design matrix as in (i).
Observe that $W$ can be written as $W=nf_{\max }f_{\max }^{\prime }-I_{n}$,
a matrix which obviously maps $\limfunc{span}(X)$ into $\limfunc{span}(X)$
as the intercept has been assumed to be an element of $\limfunc{span}(X)$.
Consequently, also $I_{n}-\rho W$ maps $\limfunc{span}(X)$ into $\limfunc{%
span}(X)$ for every $\rho \in \lbrack 0,\lambda _{\max }^{-1})$. Because $%
I_{n}-\rho W$ is nonsingular for $\rho $ in that range, it follows that this
mapping is onto and furthermore that also $\left( I_{n}-\rho W\right) ^{-1}$
maps $\limfunc{span}(X)$ into $\limfunc{span}(X)$ in a bijective way. As a
consequence, the mean of $\mathbf{y}$, which equals $(I_{n}-\rho
W)^{-1}X\beta $, is an element $X\gamma \left( \beta ,\rho \right) $, say,
of $\limfunc{span}(X)$ for every $\rho \in \lbrack 0,\lambda _{\max }^{-1})$
and $\beta \in \mathbb{R}^{k}$.\footnote{%
Compare Footnote 2 in \cite{Mart11}, where the author attempts to justify
invariance in case of a spatial lag model. The argument given there to show
that $(I-\rho W)^{-1}$ for $W$ an equal weights matrix maps $\limfunc{span}%
(X)$ into itself, however, does not make sense as it is based on an
incorrect expression for $E(\mathbf{y})$, which is incorrectly given as $%
\Delta _{\rho }X$.} Let $\varphi $ be any $G_{X}^{+}$-invariant test. Then
by $G_{X}^{+}$-invariance we have 
\begin{equation*}
E\varphi \left( \mathbf{y}\right) =E\varphi \left( X\gamma \left( \beta
,\rho \right) +\sigma (I_{n}-\rho W)^{-1}\mathbf{\varepsilon }\right)
=E\varphi \left( \sigma (I_{n}-\rho W)^{-1}\mathbf{\varepsilon }\right)
=E\varphi \left( X\beta +\sigma (I_{n}-\rho W)^{-1}\mathbf{\varepsilon }%
\right) ,
\end{equation*}%
which coincides with the power function in a SEM as in (i) above and thus is
independent of $\beta $, $\sigma $, and $\rho $, showing that the power
function of any $G_{X}^{+}$-invariant test in the SLM considered here is
constant. $\square $
\end{example}

\section{An application to time-series regression models\label{AR1C}}

In this section we briefly comment on the case where the error vector $%
\mathbf{u}$ in (\ref{linmod}) has covariance matrix $\Sigma (\rho )$ for $%
\rho \in \lbrack 0,1)$ with the $(i,j)$-th element of $\Sigma (\rho )$ given
by $\rho ^{|i-j|}$ (Case I) or $\left( -\rho \right) ^{|i-j|}$ (Case II).
Clearly, Case I corresponds to testing against positive autocorrelation,
while Case II corresponds to testing against negative autocorrelation. More
precisely, in both cases we assume that $\mathbf{u}$ is distributed as $%
\sigma \Sigma ^{1/2}(\rho )\mathbf{\varepsilon }$, where $\mathbf{%
\varepsilon }$ has mean zero, has covariance matrix $I_{n}$, and has a fixed
distribution that is spherically symmetric (and hence does not depend on any
parameters); in particular, Assumption \ref{ASDR} is maintained.
Furthermore, assume that $\Pr (\mathbf{\varepsilon }=0)=0$. \emph{We shall
refer to these assumptions as the maintained assumptions of this section}.
This framework clearly covers the case where the vector $\mathbf{u}$ is a
segment of a Gaussian stationary autoregressive process of order $1$. In
Case I it is now readily verified that Assumption \ref{ASC} holds with $%
e=n^{-1/2}(1,\ldots ,1)^{\prime }$, while in Case II this assumption is
satisfied with $e=n^{-1/2}(-1,1,\ldots ,\left( -1\right) ^{n})^{\prime }$.
The validity of Assumption \ref{ASD} then follows from Proposition \ref%
{AsPexII}. Furthermore, Assumption \ref{ASCII} (more precisely, the
equivalent condition given in Lemma \ref{AR1}) has been shown to be
satisfied in Case I as well as in Case II in Lemma G.1 of \cite{PP13}, where
the form of the matrix $V$ (denoted by $D$ in that reference) is also given;
this lemma also establishes condition (\ref{offdiag}) in view of the fact
that obviously $\lambda _{n}\left( \Sigma (\rho )\right) \rightarrow n$ for $%
\rho \rightarrow 1$ (in Case I as well as in Case II). We thus immediately
get the following result as a special case of the results in Section \ref%
{main}:

\begin{corollary}
\label{TS}Suppose the maintained assumptions hold. Let $e$ denote $%
n^{-1/2}(1,\ldots ,1)^{\prime }$ in Case I while it denotes $%
n^{-1/2}(-1,1,\ldots ,\left( -1\right) ^{n})^{\prime }$ in Case II.

\begin{enumerate}
\item Then every $G_{X}$-invariant test $\varphi $ satisfies the conclusions
1.-4. of Corollary \ref{Cor1new} subject to replacing $\lambda _{\max }^{-1}$
by $1$, $f_{\max }$ by $e$, and where $\Lambda $ now represents a
square-root of the matrix $D$ given in Lemma G.1 of \cite{PP13}.

\item Let the critical region $\Phi _{B,\kappa }$ be given by (\ref%
{quadratic}). Assume $e\in \limfunc{bd}\left( \Phi _{B,\kappa }\right)
\backslash \limfunc{span}\left( X\right) $ (i.e., $e\notin \limfunc{span}%
\left( X\right) $ and $\kappa =T_{B}\left( e\right) \in \lbrack \lambda
_{1}\left( B\right) ,\lambda _{n-k}\left( B\right) )$ with $\lambda
_{1}\left( B\right) <\lambda _{n-k}\left( B\right) $ hold). Then:

(i) Suppose $C_{X}e$ is an eigenvector of $B$ with eigenvalue $\lambda $.
Then $\lambda =T_{B}\left( e\right) =\kappa $ and 
\begin{equation}
P_{\beta ,\sigma ,\rho }(\Phi _{B,\kappa })\rightarrow \Pr \left( \mathbf{G}%
^{\prime }\Lambda ^{\prime }\left( C_{X}^{\prime }BC_{X}-\lambda
C_{X}^{\prime }C_{X}\right) \Lambda \mathbf{G}>0\right)  \label{li_1}
\end{equation}%
for $\rho \rightarrow 1$, $\rho \in \lbrack 0,1)$, and for every $\beta \in 
\mathbb{R}^{k}$, $0<\sigma <\infty $, where $\mathbf{G}$ is a multivariate
Gaussian random vector with mean zero and covariance matrix $I_{n}$. The
limit in (\ref{li_1}) is strictly between $0$ and $1$ if $\lambda >\lambda
_{1}(B)$, whereas it equals $1$ in case $\lambda =\lambda _{1}(B)$.

(ii) Suppose $C_{X}e$ is not an eigenvector of $B$. Then $P_{\beta ,\sigma
,\rho }(\Phi _{B,\kappa })\rightarrow 1/2$ for $\rho \rightarrow 1$, $\rho
\in \lbrack 0,1)$, and for every $\beta \in \mathbb{R}^{k}$, $0<\sigma
<\infty $.
\end{enumerate}
\end{corollary}

The proof of the corollary is similar to the proof of Corollaries \ref%
{Cor1new} and \ref{Cor2new} and consists of a straightforward application of
Theorems \ref{MT}, \ref{BT}, and Corollary \ref{Lem_Illust_3}, noting that
condition (\ref{offdiag}) has been verified in Lemma G.1 of \cite{PP13}. At
the expense of arriving at a more complicated result, some of the maintained
assumptions like the spherical symmetry assumption could be weakened, while
nevertheless allowing the application of the results in Section \ref{main}.
In the literature often the alternative parameterization $\varsigma
^{2}\left( 1-\rho ^{2}\right) ^{-1}\Sigma (\rho )$ for the covariance matrix
of $\mathbf{u}$ is used, which just amounts to parametrizing $\sigma ^{2}$
as $\varsigma ^{2}\left( 1-\rho ^{2}\right) ^{-1}$. In view of Remark \ref%
{invariance} and $G_{X}$-invariance of the tests considered, such an
alternative reparameterization has no effect on the results in this section
at all.

Even after specializing to the Gaussian case, the preceding corollary
provides a substantial generalization of a number of results in the
literature in that (i) it allows for general $G_{X}$-invariant tests rather
than discussing some specific tests, and (ii) provides explicit expressions
for the limiting power also in the case where the limit is neither zero nor
one: \cite{kramer1985} appears to have been the first to notice that the
zero-power trap can arise for the Durbin-Watson test in that he showed that
the limiting power (as the autocorrelation tends to $1$) of the
Durbin-Watson test can be zero when one considers a linear regression model
without an intercept and with the errors following a Gaussian autoregressive
process of order one. More precisely, he established that in this model the
limiting power is zero (is one) if -- in our notation -- the vector $e$ is
outside the closure (is inside the interior) of the rejection region of the
Durbin-Watson test.\footnote{%
We note that the words "inside" and "outside" in the Corollary of \cite%
{kramer1985} should be interchanged.} Based on numerical results, he also
noted that the zero-power trap does not seem to arise in models that contain
an intercept. Subsequently, \cite{zeisel1989} showed that indeed in models
with an intercept the limiting power of the Durbin-Watson test (except in
degenerate cases) is always strictly between zero and one.\footnote{%
The argument in \cite{zeisel1989} tacitly makes use of the Portmanteau
theorem in deriving the formula for the limiting rejection probability
without providing the necessary justification. For a more complete proof
along the same lines as the one in \cite{zeisel1989} see \cite{lobus2000}.}
The results in \cite{kramer1985} and \cite{zeisel1989} just mentioned are
extended in \cite{KramerZeisel1990} from the Durbin-Watson test to tests
that can be expressed as ratios of quadratic forms, see also \cite{small1993}%
.\footnote{%
See Remark \ref{kraemer}(i).} [We note that \cite{kramer1985} and \cite%
{KramerZeisel1990} additionally also consider the case where the
autocorrelation tends to $-1$.] All these results can be easily read off
from Part 1 of our Corollary \ref{TS}. The analysis in \cite{kramer1985}, 
\cite{zeisel1989}, and \cite{KramerZeisel1990} always excludes a particular
case, which is treated in \cite{lobus2000} for the Durbin-Watson test. This
result is again easily seen to be a special case of Part 2 of our Corollary %
\ref{TS}. Furthermore, \cite{zeisel1989} shows that for any sample size $n$
and number of regressors $k<n$ a design matrix exists such that zero-power
trap arises. For a systematic investigation of the set of regressors for
which the zero-power trap occurs see \cite{Prein2014}.

\appendix{}

\section{Comments on and counterexamples to Theorem 1 in Martellosio (2010) 
\label{Newapp}}

As already mentioned in Section \ref{main}, the first and third claim in MT1
are correct, but the proof of these statements as given in \cite{Mart10} is
not (cf. also \cite{mynbaev2012}). To explain the mistake, we assume for
simplicity that $\mathbf{u}$ is Gaussian. With this additional assumption
the model satisfies all the requirements imposed in \cite{Mart10}, page 154
(cf. Remark \ref{modelM10} above). The proof of MT1 in \cite{Mart10} is
given for $\beta $ arbitrary and $\sigma =1$. Set $\beta =0$ for simplicity.
In the proof of MT1 it is argued that the density of $\mathbf{y}$ tends, as $%
\rho \rightarrow a$, to a degenerate "density" which is supported on a set
that simplifies to the eigenspace of $\Sigma ^{-1}(a-)$ corresponding to its
smallest eigenvalue in the case $\beta =0$ considered here. However, for $%
\rho \in \lbrack 0,a)$, the density of $\mathbf{y}$ is 
\begin{equation*}
f(y)=(2\pi )^{-n/2}\left( \det \left( \Sigma ^{-1}(\rho )\right) \right)
^{1/2}\exp \left\{ -\frac{1}{2}y^{\prime }\Sigma ^{-1}(\rho )y\right\} .
\end{equation*}%
As $\rho \rightarrow a$ we have $\det \left( \Sigma ^{-1}(\rho )\right)
\rightarrow 0$ in view of the assumption $\limfunc{rank}\left( \Sigma
^{-1}(a-)\right) =n-1$. Furthermore, $\exp \left\{ -\frac{1}{2}y^{\prime
}\Sigma ^{-1}(\rho )y\right\} \rightarrow \exp \left\{ -\frac{1}{2}y^{\prime
}\Sigma ^{-1}(a-)y\right\} $ (even uniformly on compact subsets). Therefore
the density converges to zero everywhere (and even uniformly on compact
subsets). In particular, it does not tend to a degenerate "density"
supported on the eigenspace of $\Sigma ^{-1}(a-)$ corresponding to its
smallest eigenvalue in any suitable way. Note that $P_{0,1,\rho }$ does also
not converge weakly as $\rho \rightarrow a$ as the sequence $P_{0,1,\rho
_{m}}$ for any $\rho _{m}\rightarrow a$ is obviously not tight. This shows
that the proof in \cite{Mart10} is incorrect. Furthermore, the concentration
effect discussed after the theorem in \cite{Mart10} simply does not occur in
the way as claimed. In fact, the direct opposite happens: the distributions
stretch out, i.e., all of the mass "escapes to infinity".

We next turn to the second claim in MT1 and show by two simple
counterexamples that this claim is not correct.\footnote{\cite{mynbaev2012}
also claims to provide counterexamples to the second claim in MT1. However,
strictly speaking, these examples are not counterexamples as the tests
constructed always have either size $0$ or $1$, a case ruled out in the main
body of \cite{Mart10}.} The first example below is based on the following
simple observation: Suppose $k=0$, the testing problem satisfies all
assumptions of MT1, and we can find an invariant rejection region $\tilde{%
\Phi}$ of size $\alpha $, $0<\alpha <1$, with $e\in \limfunc{int}(\tilde{\Phi%
})$. The (correct) first claim of MT1 then implies that the limiting power
of $\tilde{\Phi}$ is $1$. Now define $\Phi =\tilde{\Phi}\backslash \limfunc{%
span}(e)$ and observe that $\Phi $ is again invariant and that $\Phi $ and $%
\tilde{\Phi}$ have the same rejection probabilities as they differ only by a 
$\mu _{\mathbb{R}^{n}}$-null set and the family $\mathfrak{P}$ is dominated
by $\mu _{\mathbb{R}^{n}}$ under the assumptions in \cite{Mart10}. Now $e\in 
\limfunc{bd}(\Phi )$ holds, but the limiting power of $\Phi $ is obviously $%
1 $. A concrete counterexample is as follows:

\begin{example}
\label{EX1} Assume that the elements $P_{\beta ,\sigma ,\rho }$ of the
family $\mathfrak{P}$ are Gaussian, i.e., $\mathfrak{P}$ satisfies
Assumption \ref{ASDR} with $\mathbf{z}$ a standard normally distributed
vector (and without loss of generality we may set $L(\rho )=\Sigma
^{1/2}(\rho )$). For simplicity we consider the case without regressors
(i.e., we assume $k=0$ and thus $\beta =0$ holds by our conventions). Let 
\begin{equation*}
\Sigma (\rho )=I_{n}+(1-\rho )^{-1}\rho ee^{\prime }
\end{equation*}%
for every $\rho \in \lbrack 0,1)$ where $e$ is normalized. Clearly, $\Sigma
(\rho )$ is symmetric and positive definite for $\rho \in \lbrack 0,1)$ and $%
\Sigma (0)=I_{n}$ holds. Observe that $\Sigma ^{-1}(\rho )=I_{n}-\rho
ee^{\prime }$, and thus $\Sigma ^{-1}(1-)=I_{n}-ee^{\prime }$, which has
rank $n-1$. The family $\mathfrak{P}$ hence clearly satisfies all the
assumptions for MT1 imposed in \cite{Mart10}, cf. Remark \ref{modelM10} in
Section \ref{Framework} above. Now fix an arbitrary $\alpha \in (0,1)$ and
choose a rejection region $\tilde{\Phi}\in \mathcal{B}(\mathbb{R}^{n})$ that
is (i) invariant w.r.t. $G_{0}$, (ii) satisfies $P_{0,1,0}(\tilde{\Phi}%
)=\alpha $ (and thus $P_{0,\sigma ,0}(\tilde{\Phi})=\alpha $ for every $%
0<\sigma <\infty $ by $G_{0}$-invariance), and (iii) $e\in \limfunc{int}(%
\tilde{\Phi})$. [For example, choose $M$ equal to a spherical cap on the
unit sphere $S^{n-1}$ centered at $e$ such that $M$ has measure $\alpha /2$
under the uniform distribution on $S^{n-1}$, and set $\tilde{\Phi}=\left\{
\gamma y:\gamma \neq 0,y\in M\right\} $.] From Remark \ref{RMT}(i) we obtain
that the limiting power of $\tilde{\Phi}$ is $1$. [The assumptions of
Theorem \ref{MT} are obviously satisfied in view of Lemma \ref{convCM} and
Proposition \ref{AsPexII}.] We now define a new rejection region $\Phi =%
\tilde{\Phi}\backslash \limfunc{span}(e)$. Clearly, $\Phi $ is also $G_{0}$%
-invariant, and $\Phi $ and $\tilde{\Phi}$ have the same rejection
probabilities since $\limfunc{span}(e)$ is an $\mu _{\mathbb{R}^{n}}$-null
set (as we have assumed $n\geq 2$) and the elements of $\mathfrak{P}$ are
absolutely continuous w.r.t. $\mu _{\mathbb{R}^{n}}$. However, now $e\in 
\limfunc{bd}(\Phi )$ holds, showing that the second claim in MT1 is
incorrect. A similar example, starting with the rejection region $\tilde{\Psi%
}=\mathbb{R}^{n}\backslash \tilde{\Phi}$, where $\tilde{\Phi}$ is as before,
and then passing to $\Psi =\mathbb{R}^{n}\backslash \Phi $ provides an
example where $e\in \limfunc{bd}(\Psi )$ holds, but the limiting power is
zero. $\square $
\end{example}

The argument underlying this counterexample works more generally for any
covariance model $\Sigma (\cdot )$ that satisfies the assumptions of Theorem %
\ref{MT}, and thus, in particular, for spatial models.

While the rejection region $\Phi $ constructed in the preceding example
certainly provides a counterexample to the second claim in MT1, one could
argue that it is somewhat artificial since $\Phi $ can be modified by a $\mu
_{\mathbb{R}^{n}}$-null set into the rejection region $\tilde{\Phi}$ which
does not have $e$ on its boundary. One could therefore ask if there is a
more genuine counterexample to the second claim of MT1 in the sense that the
rejection region in such a counterexample can not be modified by a $\mu _{%
\mathbb{R}^{n}}$-null set in such a way that the modified region does not
have $e$ on its boundary. This is indeed the case as shown by the subsequent
example.

\begin{example}
\label{EX2} Consider the same model as in the previous example, except that
we now assume $n=2$ and $\Sigma (\rho )$ is given by 
\begin{equation*}
\Sigma (\rho )=I_{n}+(1-\rho )^{-1}\rho e(\rho )e^{\prime }(\rho )
\end{equation*}%
for $\rho \in \lbrack 0,1)$ where $e(\rho )=(\cos (\phi (\rho )),\sin (\phi
(\rho )))^{\prime }$ with $\phi $ a strictly monotone and continuous
function on $[0,1)$ satisfying $\phi (0)=0$ and $\phi (1-)=\pi /2$. Again $%
\Sigma (\rho )$ is symmetric and positive definite for $\rho \in \lbrack
0,1) $ and $\Sigma (0)=I_{n}$ holds. Observe that $\Sigma ^{-1}(\rho
)=I_{n}-\rho e(\rho )e^{\prime }(\rho )$ holds and thus $\Sigma
^{-1}(1-)=I_{n}-ee^{\prime }$ where $e=(0,1)^{\prime }$. Obviously, $\Sigma
^{-1}(1-)$ has rank $n-1$. Again the family $\mathfrak{P}$ satisfies all the
assumptions for MT1 imposed in \cite{Mart10}. Consider the rejection region $%
\Phi =\left\{ y\in \mathbb{R}^{2}:y_{1}y_{2}\geq 0\right\} $ which is $G_{0}$%
-invariant. The rejection probability under the null is always equal to $1/2$%
. Furthermore, $e\in \limfunc{bd}(\Phi )$ holds (and obviously there is no
modification by a $\mu _{\mathbb{R}^{n}}$-null set such that $e\notin 
\limfunc{bd}(\Phi )$). We next show that $P_{0,1,\rho }(\Phi )$ converges to 
$1$ for $\rho \rightarrow 1$ under a suitable choice of the function $\phi $%
: By $G_{0}$-invariance, 
\begin{equation}
P_{0,1,\rho }(\Phi )=Q_{0,(1-\rho )I_{n}+\rho e(\rho )e^{\prime }(\rho
)}(\Phi )  \label{help}
\end{equation}%
where $Q_{0,\Omega }$ denotes the Gaussian measure on $\mathbb{R}^{n}$ with
mean zero and variance covariance matrix $\Omega $. Now for fixed $\eta $, $%
0<\eta <1$, we have that $e(\eta )\in \limfunc{int}\left( \Phi \right) $
because of strict monotonicity of $\phi $. Furthermore,%
\begin{equation*}
Q_{0,(1-\rho )I_{n}+\rho e(\eta )e^{\prime }(\eta )}(\Phi )\rightarrow
Q_{0,e(\eta )e^{\prime }(\eta )}(\Phi )\geq Q_{0,e(\eta )e^{\prime }(\eta )}(%
\limfunc{span}(e(\eta )))=1
\end{equation*}%
because $Q_{0,(1-\rho )I_{n}+\rho e(\eta )e^{\prime }(\eta )}$ converges to $%
Q_{0,e(\eta )e^{\prime }(\eta )}$ weakly, and because 
\begin{equation*}
Q_{0,e(\eta )e^{\prime }(\eta )}(\limfunc{bd}(\Phi ))=Q_{0,e(\eta )e^{\prime
}(\eta )}(\limfunc{bd}(\Phi )\cap \limfunc{span}(e(\eta )))=Q_{0,e(\eta
)e^{\prime }(\eta )}(\left\{ 0\right\} )=0.
\end{equation*}%
It is now obvious that if $\phi (\rho )$ converges to $\pi /2$ sufficiently
slowly, we can also achieve that (\ref{help}) converges to $1$ as $\rho
\rightarrow 1$. Furthermore, we also conclude that the invariant rejection
region $\Psi =\mathbb{R}^{2}\backslash \Phi $, which also has rejection
probability $1/2$ under the null, provides an example where $e\in \limfunc{bd%
}(\Psi )$ holds, but the limiting power is zero. $\square $
\end{example}

Similar counterexamples to the second claim in MT1 can also be constructed
when regressors are present (except if $n=k+1$).\footnote{\label{FnA}The
case $n=k+1$ is somewhat trivial as we now explain: If $n=k+1$, every $G_{X}$%
-invariant test $\varphi $ is $\mu _{\mathbb{R}^{n}}$-almost everywhere
constant. [To see this observe that $\limfunc{span}(X)$ is a $\mu _{\mathbb{R%
}^{n}}$-null set and that every element of $\limfunc{span}(X)^{\bot }$ is of
the form $\lambda b$ for a fixed vector $b$ and hence $\varphi \left(
y\right) =\varphi (\Pi _{\limfunc{span}(X)^{\bot }}y)=\varphi \left( \lambda
b\right) =\varphi \left( b\right) $ holds whenever $\lambda \neq 0$, i.e.,
whenever $y\notin \limfunc{span}(X)$. Additionally, note that $\varphi $ is
constant on $\limfunc{span}(X)$.] Consequently, $\varphi $ has a constant
power function if (i) the family of probability measures in (\ref{family})
is absolutely continuous w.r.t. $\mu _{\mathbb{R}^{n}}$, or if (ii) this
family is an elliptically symmetric family (to see this in case $\Pr \left( 
\mathbf{z}=0\right) =0$ use the argument given in Remark \ref{RBT}(iv); in
case $\Pr \left( \mathbf{z}=0\right) >0$ combine this argument with Remark %
\ref{rem_gen_1}(vi)). In particular, if $\varphi $ is non-randomized, it is
then a trivial test in that its size and power are either both zero or one,
provided (i) holds or (ii) holds with $\Pr \left( \mathbf{z}=0\right) =0$.}

\section{Comments on further results in Martellosio (2010)\label{A2}}

In this section we comment on problems in some results in \cite{Mart10} that
have not been discussed so far. We also discuss if and how these problems
can be fixed.

\subsection{Comments on Lemmata D.2 and D.3 in Martellosio (2010)\label{A1}}

Here we discuss problems with Lemmata D.2 and D.3 in Martellosio (2010)
which are phrased in a spatial error model context. Correct versions of
these lemmata, which furthermore are also not restricted to spatial
regression models, have been given in Section \ref{alpha_star} above. Both
lemmata in \cite{Mart10} concern the quantity $\alpha ^{\ast }$, which is
defined on p. 165 of \cite{Mart10} as follows:

\begin{quote}
"For an exact invariant test of $\rho =0$ against $\rho >0$ in a SAR(1)
model, $\alpha ^{\ast }$ is the infimum of the set of values of $\alpha \in
(0,1]$ such that the limiting power does not vanish."
\end{quote}

In this definition $\alpha $ denotes a generic symbol for the size of the
test. Taken literally, the definition refers to one test only and hence does
not make sense (as there is then only one associated value of $\alpha $).
From later usage of this definition in \cite{Mart10}, it seems that the
author had in mind a \emph{family} of tests (rejection regions) like $\Phi
_{\kappa }=\left\{ y\in \mathbb{R}^{n}:T(y)>\kappa \right\} $, where $T$ is
a test statistic. Interpreting Martellosio's definition this way, it is
clear that under the assumptions made in \cite{Mart10} (see Remark \ref%
{modelM10}(ii) and Remark \ref{invariance} above) his $\alpha ^{\ast }$
coincides with $\alpha ^{\ast }\left( T\right) $ defined in (\ref{astar}).

Lemma D.2 of \cite{Mart10}, p. 181, then reads as follows:

\begin{quote}
"Consider a model $G(X\beta ,\sigma ^{2}[(I-\rho W^{\prime })(I-\rho
W)]^{-1})$, where $G(\mu ,\Gamma )$ denotes some multivariate distribution
with mean $\mu $ and variance matrix $\Gamma $. When an invariant critical
region for testing $\rho =0$ against $\rho >0$ is in form (9) [i.e., is of
the form $\left\{ y\in \mathbb{R}^{n}:T(y)>\kappa \right\} $ for some
univariate test statistic $T$], and is such that $f_{\max }$ is not
contained in its boundary, $\alpha ^{\ast }=\Pr (T(\boldsymbol{z})>T(f_{\max
});\boldsymbol{z}\sim G(0,I))$."
\end{quote}

The statement of this lemma as well as its proof are problematic for the
following reasons:

\begin{enumerate}
\item The lemma makes a statement about $\alpha ^{\ast }$, which is a
quantity that depends not only on one specific critical region, but on a 
\emph{family} of critical regions corresponding to a \emph{family} of
critical values $\kappa $ against which the test statistic is compared. The
critical region usually depends on $\kappa $ and so does its boundary (cf.
Proposition \ref{prop_bound}). Therefore, the assumption "... $f_{\max }$ is
not contained in its [the invariant critical region's] boundary..." has
little meaning in this context as it is not clear to which one of the many
rejection regions the statement refers to. [Alternatively, if one interprets
the statement of the lemma as requiring $f_{\max }$ not to be contained in
the boundary of \emph{every} rejection region in the family considered, this
leads to a condition that typically will never be satisfied.]

\item The proof of the lemma is based on Corollary 1 in \cite{Mart10}, the
proof of which is incorrect as it is based on the incorrect Theorem 1 of 
\cite{Mart10}.

\item The proof implicitly uses a continuity assumption on the cumulative
distribution function of the test statistic under the null at the point $%
T(f_{\max })$ which is not satisfied in general.
\end{enumerate}

Next we turn to Lemma D.3 in \cite{Mart10}, which reads:

\begin{quote}
"Consider a test that, in the context of a spatial error model with
symmetric $W$, rejects $\rho =0$ for small values of a statistic $\nu
^{\prime }B\nu $, where $B$ is an $(n-k)\times (n-k)$ known symmetric matrix
that does not depend on $\alpha $, and $\nu $ is as defined in Section 2.2.
Provided that $f_{\max }\notin \limfunc{bd}(\Phi )$, $\alpha ^{\ast }=0$ if
and only if $Cf_{\max }\in E_{1}(B)$, and $\alpha ^{\ast }=1$ if and only if 
$Cf_{\max }\in E_{n-k}(B)$."
\end{quote}

Here $\alpha $ refers to the size of the test, $\nu $ is given by $\limfunc{%
sign}(y_{i})Cy/\Vert {Cy}\Vert $ for some fixed $i\in \left\{ 1,\ldots
,n\right\} $, and $\Phi $ is not explicitly defined, but presumably denotes
a rejection region corresponding to the test statistic $\nu ^{\prime }B\nu $%
. [Although the test statistic is not defined whenever ${Cy=0}$, this does
not pose a severe problem here since \cite{Mart10} considers only absolutely
continuous distributions and since he assumes $k<n$; cf Remark \ref%
{effect_on_boundary}. Note furthermore that the factor $\limfunc{sign}%
(y_{i}) $ is irrelevant here.] Furthermore, $E_{1}(B)$ ($E_{n-k}(B)$)
denotes the eigenspace corresponding to the smallest (largest) eigenvalue of 
$B$, and $C$ in \cite{Mart10} stands for $C_{X}$. The statement of the lemma
and its content are inappropriate for the following reasons:

\begin{enumerate}
\item The proof of this lemma is based on Lemma D.2 of \cite{Mart10} which
is invalid as discussed above.

\item Again, as in the statement of Lemma D.2 of \cite{Mart10}, the author
assumes that `... $f_{\max }\notin \limfunc{bd}(\Phi )$ ...', which is not
meaningful, as the boundary typically depends on the critical value.

\item The above lemma in \cite{Mart10} requires $W$ to be symmetric
(although this is actually not used in the proof). Nevertheless, it is later
applied to nonsymmetric weights matrices in the proof of Proposition 1 in 
\cite{Mart10}.
\end{enumerate}

As a point of interest we note that naively applying Lemma D.3 in \cite%
{Mart10} to the case where $B$ is a multiple of the identity matrix $I_{n-k}$
leads to the contradictory statement $0=\alpha ^{\ast }=1$. However, in case 
$B$ is a multiple of $I_{n-k}$, the test statistic degenerates, and thus the
size of the test is $0$ or $1$, a case that is ruled out in \cite{Mart10}
from the very beginning.

\subsection{Comments on Proposition 1 and Lemma E.4 in Martellosio (2010) 
\label{Prop1}}

Proposition 1 in \cite{Mart10} considers the pure SAR(1) model, i.e., $k=0$
is assumed. This proposition reads as follows:

\begin{quote}
"Consider testing $\rho =0$ against $\rho >0$ in a pure SAR(1) model. The
limiting power of the Cliff-Ord test [cf. eq. (\ref{COT}) below] or of a
test (8) [cf. eq. (\ref{LRT}) below] is $1$ irrespective of $\alpha $ [the
size of the test] if and only if $f_{\max }$ is an eigenvector of $W^{\prime
}$."
\end{quote}

We note that, while not explicit in the above statement, it is understood in 
\cite{Mart10} that $0\leq \rho <\lambda _{\max }^{-1}$ \ is assumed.
Similarly, the case $n=1$ is not ruled out explicitly in the statement of
the proposition, but it seems to be implicitly understood in \cite{Mart10}
that $n\geq 2$ holds (note that in case $n=1$ the test statistics degenerate
and therefore the associated tests trivially have size equal to $0$ or $1$,
depending on the choice of the critical value).

The test defined in equation (8) of \cite{Mart10} rejects for small values
of 
\begin{equation}
y^{\prime }(I_{n}-\bar{\rho}W^{\prime })(I_{n}-\bar{\rho}W)y/\Vert {y}\Vert
^{2},  \label{LRT}
\end{equation}%
where $0<\bar{\rho}<\lambda _{\max }^{-1}$ is specified by the user. The
argument in the proof of the proposition in \cite{Mart10} for this class of
tests is incorrect for the following reasons:

\begin{enumerate}
\item The proof is based on Lemma D.3 in \cite{Mart10} which is incorrect as
discussed in Appendix \ref{A1}.

\item Even if Lemma D.3 in \cite{Mart10} were correct and could be used,
this lemma would only deliver the result $\alpha ^{\ast }=0$ which does 
\emph{not} imply, without a further argument, that the limiting power is
equal to one for every size $\alpha \in (0,1)$. By definition of $\alpha
^{\ast }$, $\alpha ^{\ast }=0$ only implies that the limiting power is \emph{%
nonzero} for every size $\alpha \in (0,1)$.
\end{enumerate}

For the case of the Cliff-Ord test, i.e., the test rejecting for small
values of 
\begin{equation}
-y^{\prime }Wy/\Vert {y}\Vert ^{2}=-0.5y^{\prime }(W+W^{\prime })y/\Vert {y}%
\Vert ^{2},  \label{COT}
\end{equation}%
\cite{Mart10} argues that this can be reduced to the previously considered
case, the proof of which is flawed as just shown. Apart from this, the
reduction argument, which we now quote, has its own problems:

\begin{quote}
"... By Lemma D.3 with $B=\Gamma ^{-1}(\bar{\rho})$ [which equals $(I_{n}-%
\bar{\rho}W^{\prime })(I_{n}-\bar{\rho}W)$], in order to prove that the
limiting power of test (8) [cf. eq. (\ref{LRT}) above] is $1$ for any $%
\alpha $ [the size of the test], we need to show that $W^{\prime }f_{\max
}=\lambda _{\max }f_{\max }$ is necessary and sufficient for $f_{\max }\in
E_{n}(\Gamma (\bar{\rho}))$. Clearly, if this holds for any $\bar{\rho}>0$,
it holds for $\bar{\rho}\rightarrow 0$ too, establishing also the part of
the proposition regarding the Cliff-Ord test. ..."
\end{quote}

The problem here is that it is less than clear what the precise mathematical
"approximation" argument is. If we interpret it as deriving limiting power
equal to $1$ for the Cliff-Ord test from the corresponding result for tests
of the form (8) and the fact that the Cliff-Ord test emerges as a limit of
these tests for $\bar{\rho}\rightarrow 0$, then this involves an interchange
of two limiting operations, namely $\rho \rightarrow \lambda _{\max }^{-1}$
and $\bar{\rho}\rightarrow 0$, for which no justification is provided.
Alternatively, one could try to interpret the "approximation" argument as an
argument that tries to derive $f_{\max }\in E_{n}(W+W^{\prime })$ from $%
f_{\max }\in E_{n}(\Gamma (\bar{\rho}))$ for every $\bar{\rho}>0$; of
course, such an argument would need some justification which, however, is
not provided. We note that this argument could perhaps be saved by using the
arguments we provide in the proof of Proposition \ref{E4new}, but the proof
of our correct version of Proposition 1 in \cite{Mart10}, i.e., Proposition %
\ref{prop1new} in Section \ref{SEM}, is more direct and does not need such a
reasoning. Furthermore, note that the proof of Proposition \ref{prop1new} is
based on our Proposition \ref{D3new}, which is a correct version of Lemma
D.3 in \cite{Mart10} and which delivers not only the conclusion $\alpha
^{\ast }=0$, but the stronger conclusion that the limiting power is indeed
equal to $1$ for every size in $(0,1)$.

We now turn to a discussion of Lemma E.4, which is again a statement about
the Cliff-Ord test and tests of the form (8) in \cite{Mart10}, but now in
the context of the SEM (i.e., $k>0$ is possible). The statement and the
proof of the lemma suffer from the following shortcomings (again Lemma E.4
implicitly assumes that $0\leq \rho <\lambda _{\max }^{-1}$ holds):

\begin{enumerate}
\item The proof of the lemma is based on Lemma D.3 in \cite{Mart10}, which
is incorrect (cf. the discussion in Appendix \ref{A1}).

\item The proof uses non-rigorous arguments such as arguments involving a
`limiting matrix' with an infinite eigenvalue. Additionally, continuity of
the dependence of eigenspaces on the underlying matrix is used without
providing the necessary justification.

\item For the case of the Cliff-Ord test the same unjustified reduction
argument as in the proof of Proposition 1 of \cite{Mart10} is used, cf. the
preceding discussion.
\end{enumerate}

For a correct version of Lemma E.4 of \cite{Mart10} see Proposition \ref%
{E4new}\ in Section \ref{SEM} above. As a point of interest we furthermore
note that cases where the test statistics become degenerate (e.g., the case $%
n-k=1$) are not ruled out explicitly in Lemma E.4 in \cite{Mart10}; in these
cases $\alpha ^{\ast }=1$ (and not $\alpha ^{\ast }=0$) holds.

\subsection{Comments on Propositions 3, 4, and 5 in Martellosio (2010) \label%
{P3}}

The proof of the part of Proposition 3 of \cite{Mart10} regarding
point-optimal invariant tests seems to be correct except for the case where $%
\limfunc{span}(X)^{\bot }$ is contained in one of the eigenspaces of $\Sigma
(\rho )$. In this case the test statistic of the form (8) in \cite{Mart10}
is degenerate (see Section \ref{idsp} above) and does not give the
point-optimal invariant test (except in the trivial case where the size is $%
0 $ or $1$, a case always excluded in \cite{Mart10}). However, this problem
is easily fixed by observing that the point-optimal invariant test in this
case is given by the randomized test $\varphi \equiv \alpha $, which is
trivially unbiased. Two minor issues in the proof are as follows: (i) Lemma
E.3 can only be applied as long as $\mathbf{z}_{i}^{2}>0$ for every $i\in H$%
. Fortunately, the complement of this event is a null-set allowing the
argument to go through. (ii) The expression `stochastically larger' in the
paragraph following (E.4) should read `stochastically smaller'. We also note
that the assumption of Gaussianity can easily be relaxed to elliptical
symmetry in view of $G_{X}$-invariance of the tests considered.

More importantly, the proof of the part of Proposition 3 of \cite{Mart10}
concerning locally best invariant tests is highly deficient for at least two
reasons: First, it is claimed that locally best invariant tests are of the
form (7) in \cite{Mart10} with $Q=d\Sigma (\rho )/d\rho |_{\rho =0}$. While
this is correct under regularity conditions (including a differentiability
assumption on $\Sigma (\rho )$), such conditions are, however, missing in
Proposition 3 of \cite{Mart10}. Also, the case where $\limfunc{span}%
(X)^{\bot }$ is contained in one of the eigenspaces of $\Sigma (\rho )$ has
to be treated separately, as then the locally best invariant test is given
by the randomized test $\varphi \equiv \alpha $. Second, the proof uses once
more an unjustified approximation argument in an attempt to reduce the case
of locally best invariant tests to the case of point-optimal invariant
tests. It is not clear what the precise nature of the approximation argument
is. Furthermore, even if the approximation argument could be somehow
repaired to deliver unbiasedness of locally best invariant tests, it is less
than clear that strict unbiasedness could be obtained this way as strict
inequalities are not preserved by limiting operations.

We next turn to the part of Proposition 4 of \cite{Mart10} regarding
point-optimal invariant tests.\footnote{%
While not explicit in the statement of this proposition, it is implicitly
assumed that $0\leq \rho <\lambda _{\max }^{-1}$ holds.} As in the case of
Proposition 3 discussed above, the case where $\limfunc{span}(X)^{\bot }$ is
contained in one of the eigenspaces of $\Sigma (\rho )$ has to be treated
separately, and Gaussianity can be relaxed to elliptical symmetry. We note
that the clause `if and only if' in the last but one line of p. 185 of \cite%
{Mart10} should read `if'. We also note that the verification of the first
displayed inequality on p. 186 of \cite{Mart10} could be shortened (using
Lemma E.3 (more precisely, the more general result referred to in the proof
of this lemma) with $a_{i}=\lambda _{i}(W)/\tau _{i}(\rho )$, $b_{i}=\tau
_{i}^{2}(\bar{\rho})$, and $p_{i}=\mathbf{z}_{i}^{2}/\tau _{i}^{2}(\rho )$
to conclude that the first display on p. 186 holds almost surely, and
furthermore that it holds almost surely with equality if and only if all $%
b_{i}$ or all $a_{i}$ are equal, which is equivalent to all $\lambda _{i}(W)$
for $i\in H$ being equal).

Again, the proof of the part of Proposition 4 of \cite{Mart10} concerning
locally best invariant tests is deficient as it is based on the same
unjustified approximation argument mentioned before.

We next turn to Proposition 5 of \cite{Mart10}. In the last of the three
cases considered in this proposition, both test statistics are degenerate
and hence the power functions are trivially constant equal to $0$ or $1$ (a
case ruled out in \cite{Mart10}). More importantly, the proof of Proposition
5 is severely flawed for several reasons, of which we only discuss a few:
First, the proof makes use of Corollary 1 of \cite{Mart10}, the proof of
which is based on the incorrect Theorem 1 in \cite{Mart10}; it also makes
use of Lemma E.4 and Proposition 4 of \cite{Mart10} which are incorrect as
discussed before. Second, even if these results used in the proof were
correct as they stand, additional problems would arise: Lemma E.4 only
delivers $\alpha ^{\ast }=0$, and not the stronger conclusion that the
limiting power equals $1$, as would be required in the proof. Furthermore,
Proposition 4 has Gaussianity of the errors as a hypothesis, while such an
assumption is missing in Proposition 5.

We conclude by mentioning that a correct version of the part of Proposition
5 of \cite{Mart10} concerning tests of the form (8) in \cite{Mart10} can
probably be obtained by substituting our Corollary \ref{Cor1new} and
Proposition \ref{E4new} for Corollary 1 and Lemma E.4 of \cite{Mart10} in
the proof, but we have not checked the details. For the Cliff-Ord test this
does not seem to work in the same way as the corresponding case of
Proposition 4 of \cite{Mart10} is lacking a proof as discussed before.

\section{Proofs for Section \protect\ref{main}\label{app_proofs}}

\textbf{Proof of Lemma \ref{convCM}:} Let $\rho _{m}$ be a sequence in $%
[0,a) $ converging to $a$ and let 
\begin{equation*}
\sum_{j=1}^{n-1}\lambda _{j}(\Sigma (\rho _{m}))v_{j}(\rho _{m})v_{j}(\rho
_{m})^{\prime }+\lambda _{n}(\Sigma (\rho _{m}))v_{n}(\rho _{m})v_{n}(\rho
_{m})^{\prime }
\end{equation*}%
be a spectral decomposition of $\Sigma (\rho _{m})$, with $v_{j}(\rho _{m})$
($j=1,\ldots ,n$) forming an orthonormal basis of eigenvectors of $\Sigma
(\rho _{m})$ and $\lambda _{j}(\Sigma (\rho _{m}))$ for $j=1,\ldots ,n$
denoting the corresponding eigenvalues ordered from smallest to largest and
counted with their multiplicities. Because $\Sigma ^{-1}(a-)$ is
rank-deficient by assumption, we must have $\lambda _{1}(\Sigma ^{-1}(\rho
_{m}))\rightarrow 0$, or equivalently $\lambda _{n}^{-1}(\Sigma (\rho
_{m}))\rightarrow 0$. Because the kernel of $\Sigma ^{-1}(a-)$ has dimension
one and because of positive definiteness of $\Sigma (\rho _{m})$ we can
infer the existence of some $0<M<\infty $ such that $0<\lambda _{j}(\Sigma
(\rho _{m}))<M$ must hold for every $j=1,\ldots ,n-1$ and $m\in \mathbb{N}$.
As a consequence, the sum from $1$ to $n-1$ in the previous display, after
being premultiplied by $\lambda _{n}^{-1}(\Sigma (\rho _{m}))$, converges to
zero for $m\rightarrow \infty $. It remains to show that $v_{n}(\rho
_{m})v_{n}(\rho _{m})^{\prime }\rightarrow ee^{\prime }$. Let $m^{\prime }$
be an arbitrary subsequence of $m$. By norm-boundedness of the sequence $%
v_{n}(\rho _{m})$ there exists another subsequence $m^{\prime \prime }$
along which $v_{n}(\rho _{m})$ converges to some normalized vector $e^{\ast
} $, say. Clearly 
\begin{equation*}
\Sigma ^{-1}(\rho _{m^{^{\prime \prime }}})v_{n}(\rho _{m^{^{\prime \prime
}}})=\lambda _{n}^{-1}(\Sigma (\rho _{m^{^{\prime \prime }}}))v_{n}(\rho
_{m^{^{\prime \prime }}}).
\end{equation*}%
The left hand side in the previous display now converges to $\Sigma
^{-1}(a-)e^{\ast }$ while the right hand side converges to zero. Therefore $%
e^{\ast }$ is an element of the (one-dimensional) kernel of $\Sigma
^{-1}(a-) $. Since $e^{\ast }$ is normalized, we must have $e^{\ast }{%
e^{\ast }}^{\prime }=ee^{\prime }$. This proves the claim as the subsequence 
$m^{\prime }$ was arbitrary. $\blacksquare $

\begin{lemma}
\label{convPM} Let $\mathbf{v}_{m}$ be a sequence of random $n$-vectors such
that $E(\mathbf{v}_{m})=0$ and $E(\Vert \mathbf{v}{_{m}}\Vert ^{2})<\infty $
and let $\Omega _{m}=E(\mathbf{v}_{m}\mathbf{v}_{m}^{\prime })$. If $\Omega
_{m}\rightarrow ee^{\prime }$ as $m\rightarrow \infty $ for some $e\in 
\mathbb{R}^{n}$, then the sequence $\mathbf{v}_{m}$ is tight and the support
of every weak accumulation point of the sequence of distributions of $%
\mathbf{v}_{m}$ is a subset of $\limfunc{span}(e)$. If, in addition, every
weak accumulation point of the distributions of $\mathbf{v}_{m}$ has no mass
at the origin and if $e$ is normalized, then the distribution of $\mathcal{I}%
_{0,\zeta _{e}}(\mathbf{v}_{m})$ converges weakly to $\delta _{e}$.
\end{lemma}

\textbf{Proof: }Let $M$ be an arbitrary positive real number. Since the
sequence $\limfunc{trace}(\Omega _{m})$ is convergent to $\limfunc{trace}%
(ee^{\prime })$, it is bounded from above by $S$, say. Markov's inequality
gives%
\begin{equation*}
\Pr (\Vert \mathbf{v}{_{m}}\Vert \geq M)=\Pr (\Vert \mathbf{v}{_{m}}\Vert
^{2}\geq M^{2})\leq \frac{E(\mathbf{v}_{m}^{\prime }\mathbf{v}_{m})}{M^{2}}=%
\frac{\limfunc{trace}(\Omega _{m})}{M^{2}}\leq \frac{S}{M^{2}}
\end{equation*}%
for every $m\in \mathbb{N}$, which implies tightness. To prove the claim
about the support of weak accumulation points note that $\mathbf{v}_{m}=\Pi
_{\text{$\limfunc{span}$}(e)}\mathbf{v}_{m}+\Pi _{\text{$\limfunc{span}$}%
(e)^{\bot }}\mathbf{v}_{m}$ and that the support of $\Pi _{\text{$\limfunc{%
span}$}(e)}\mathbf{v}_{m}$ is certainly a subset of $\limfunc{span}(e)$,
which is a closed set. It thus suffices to show that $\Pi _{\limfunc{span}%
(e)^{\bot }}\mathbf{v}_{m}$ converges to zero in probability. But this is
again a consequence of Markov's inequality: For every $\varepsilon >0$ we
have 
\begin{equation}
\Pr (\Vert {\Pi _{\text{$\limfunc{span}$}(e)^{\bot }}\mathbf{v}_{m}}\Vert
\geq \varepsilon )\leq \frac{E(\Vert {\Pi _{\text{$\limfunc{span}$}(e)^{\bot
}}\mathbf{v}_{m}}\Vert ^{2})}{\varepsilon ^{2}}=\frac{E(\mathbf{v}%
_{m}^{\prime }\Pi _{\text{$\limfunc{span}$}(e)^{\bot }}\mathbf{v}_{m})}{%
\varepsilon ^{2}}=\frac{\limfunc{trace}(\Pi _{\text{$\limfunc{span}$}%
(e)^{\bot }}\Omega _{m})}{\varepsilon ^{2}}.  \label{markov}
\end{equation}%
Because $\Omega _{m}\rightarrow ee^{\prime }$, we obtain $\Pi _{\limfunc{span%
}(e)^{\bot }}\Omega _{m}\rightarrow 0$ and hence the upper bound in (\ref%
{markov})\ converges to zero as $m\rightarrow \infty $. To prove the final
assertion let $m^{\prime }$ be an arbitrary subsequence and $m^{\prime
\prime }$ a subsequence thereof such that $\mathbf{v}_{m^{\prime \prime }}$
converges weakly to $\mathbf{v}$, say. By what has already been established,
we may assume that $\Pi _{\limfunc{span}(e)}\mathbf{v}=\mathbf{v}$ almost
surely holds. Because $\mathcal{I}_{0,\zeta _{e}}$ is continuous at $\lambda
e$ for every $\lambda \neq 0$ and because $\Pr (\mathbf{v}=0)=0$ by the
assumptions, we can apply the continuous mapping theorem to conclude that 
\begin{equation*}
\mathcal{I}_{0,\zeta _{e}}(\mathbf{v}_{m^{\prime \prime }})\rightarrow 
\mathcal{I}_{0,\zeta _{e}}(\mathbf{v}).
\end{equation*}%
Since $\Pr (\mathbf{v}=0)=0$, we have that $\mathcal{I}_{0,\zeta _{e}}(%
\mathbf{v})$ is almost surely equal to $\zeta _{e}\left( \frac{\mathbf{v}}{%
\Vert \mathbf{v}\Vert }\right) $. But this is almost surely equal to $e$ by
definition of $\zeta _{e}$. This completes the proof because $m^{\prime }$
was an arbitrary subsequence. $\blacksquare $

\textbf{Proof of Proposition \ref{AsPexII}:} 1. Let $\rho _{m}\in \lbrack
0,a)$ be a sequence converging to $a$. Assumption \ref{ASDR} implies that $%
P_{\beta ,\sigma ,\rho _{m}}\circ M_{X\beta ,\lambda _{n}^{1/2}(\Sigma (\rho
_{m}))\sigma }$ coincides with $P_{0,\lambda _{n}^{-1/2}(\Sigma (\rho
_{m})),\rho _{m}}$, which is precisely the distribution of $\lambda
_{n}^{-1/2}(\Sigma (\rho _{m}))\Sigma ^{1/2}(\rho _{m})\mathbf{z}$. By
Assumption \ref{ASC} we have $\lambda _{n}^{-1}(\Sigma (\rho _{m}))\Sigma
(\rho _{m})\rightarrow ee^{\prime }$. By continuity of the symmetric
nonnegative definite square root we obtain%
\begin{equation*}
\lambda _{n}^{-1/2}(\Sigma (\rho _{m}))\Sigma ^{1/2}(\rho _{m})=\left(
\lambda _{n}^{-1}(\Sigma (\rho _{m}))\Sigma (\rho _{m})\right)
^{1/2}\rightarrow \left( ee^{\prime }\right) ^{1/2}=ee^{\prime }.
\end{equation*}%
Consequently, $\lambda _{n}^{-1/2}(\Sigma (\rho _{m}))\Sigma ^{1/2}(\rho
_{m})\mathbf{z}$ converges weakly to $ee^{\prime }\mathbf{z}$. Hence, the
only accumulation point $P$, say, of $P_{0,\lambda _{n}^{-1/2}(\Sigma (\rho
_{m})),\rho _{m}}$ is the distribution of $ee^{\prime }\mathbf{z}$. The
claim now follows because $P\left( \left\{ 0\right\} \right) =\Pr
(ee^{\prime }\mathbf{z}=0)=\Pr (e^{\prime }\mathbf{z}=0)=0$ by assumption

2. Let $\rho _{m}$ be as before and observe that again $P_{\beta ,\sigma
,\rho _{m}}\circ M_{X\beta ,\lambda _{n}^{1/2}(\Sigma (\rho _{m}))\sigma }$
coincides with $P_{0,\lambda _{n}^{-1/2}(\Sigma (\rho _{m})),\rho _{m}}$,
which, however, now equals the distribution of $\lambda _{n}^{-1/2}(\Sigma
(\rho _{m}))L(\rho _{m})\mathbf{z}$. Since $L(\rho _{m})$ is a square root
of $\Sigma (\rho _{m})$, there must exist an orthogonal matrix $U(\rho _{m})$
such that $L(\rho _{m})=\Sigma ^{1/2}(\rho _{m})U(\rho _{m})$. Rewrite $%
\lambda _{n}^{-1/2}(\Sigma (\rho _{m}))L(\rho _{m})$ as $\lambda
_{n}^{-1/2}(\Sigma (\rho _{m}))\Sigma ^{1/2}(\rho _{m})U(\rho _{m})$. Fix an
arbitrary subsequence $m^{\prime }$ of $m$. Along a suitable subsubsequence $%
m^{\prime \prime }$ the matrix $U(\rho _{m^{\prime \prime }})$ converges to
an orthogonal matrix $U$, say. Therefore $\lambda _{n}^{-1/2}(\Sigma (\rho
_{m^{\prime \prime }}))\Sigma ^{1/2}(\rho _{m^{\prime \prime }})U(\rho
_{m^{\prime \prime }})$ converges to $ee^{\prime }U$. Hence, the only
accumulation point $P$, say, of $P_{0,\lambda _{n}^{-1/2}(\Sigma (\rho
_{m})),\rho _{m}}$ along the subsequence $m^{\prime \prime }$ is the
distribution of $ee^{\prime }U\mathbf{z}$. But clearly $P\left( \left\{
0\right\} \right) =\Pr (ee^{\prime }U\mathbf{z}=0)=\Pr (e^{\prime }U\mathbf{z%
}=0)$. Now this is equal to $0$ in case the distribution of $\mathbf{z}$ is
dominated by $\mu _{\mathbb{R}^{n}}$ since the set $\left\{ y\in \mathbb{R}%
^{n}:e^{\prime }Uy=0\right\} $ is obviously a $\mu _{\mathbb{R}^{n}}$-null
set. Since $m^{\prime }$ was arbitrary, the proof of the first claim is
complete. To prove the second claim observe that $\Pr (e^{\prime }U\mathbf{z}%
=0)=\Pr (e^{\prime }U\left( \mathbf{z}/\left\Vert \mathbf{z}\right\Vert
\right) =0)$, which equals zero since the distribution of $\mathbf{z}%
/\left\Vert \mathbf{z}\right\Vert $ is dominated by $\upsilon _{S^{n-1}}$ by
assumption and since $\left\{ y\in S^{n-1}:e^{\prime }Uy=0\right\} $ is a $%
\upsilon _{S^{n-1}}$-null set (cf. Remark \ref{E1}(i)). $\blacksquare $

\textbf{Proof of Theorem \ref{MT}:} Let $\rho _{m}$ be a sequence in $[0,a)$
converging to $a$. Invariance of the test $\varphi $ w.r.t. $G_{X}$ implies 
\begin{align*}
E_{\beta ,\sigma ,\rho _{m}}(\varphi )& =\int_{\mathbb{R}^{n}}\varphi
(y)dP_{\beta ,\sigma ,\rho _{m}}=\int_{\mathbb{R}^{n}}\varphi \left(
M_{X\beta ,\lambda _{n}^{1/2}(\Sigma (\rho _{m}))\sigma }(y)\right)
dP_{\beta ,\sigma ,\rho _{m}} \\
& =\int_{\mathbb{R}^{n}}\varphi (y)d\left( P_{\beta ,\sigma ,\rho _{m}}\circ
M_{X\beta ,\lambda _{n}^{1/2}(\Sigma (\rho _{m}))\sigma }\right) \\
& =\int_{\mathbb{R}^{n}}\varphi \left( \mathcal{I}_{0,\zeta _{e}}(y)\right)
d\left( P_{\beta ,\sigma ,\rho _{m}}\circ M_{X\beta ,\lambda
_{n}^{1/2}(\Sigma (\rho _{m}))\sigma }\right) \\
& =\int_{\mathbb{R}^{n}}\varphi (y)d\left( \left( P_{\beta ,\sigma ,\rho
_{m}}\circ M_{X\beta ,\lambda _{n}^{1/2}(\Sigma (\rho _{m}))\sigma }\right)
\circ \mathcal{I}_{0,\zeta _{e}}\right) ,
\end{align*}%
where the last but one equality holds because of Remark \ref{MI}(ii). The
covariance matrix of $\mathbf{v}_{m}$, say, a centered random variable with
distribution $P_{\beta ,\sigma ,\rho _{m}}\circ M_{X\beta ,\lambda
_{n}^{1/2}(\Sigma (\rho _{m}))\sigma }$, is given by $\lambda
_{n}^{-1}(\Sigma (\rho _{m}))\Sigma (\rho _{m})$ which converges to $%
ee^{\prime }$ by Assumption \ref{ASC}. Note that $e$ is necessarily
normalized. By Assumption \ref{ASD} every weak accumulation point $P$ of $%
P_{\beta ,\sigma ,\rho _{m}}\circ M_{X\beta ,\lambda _{n}^{1/2}(\Sigma (\rho
_{m}))\sigma }$ satisfies $P(\left\{ 0\right\} )=0$ (note that $P_{\beta
,\sigma ,\rho _{m}}\circ M_{X\beta ,\lambda _{n}^{1/2}(\Sigma (\rho
_{m}))\sigma }$ is in fact tight by Lemma \ref{convPM}). Thus we can apply
Lemma \ref{convPM} to conclude that 
\begin{equation*}
\left( P_{\beta ,\sigma ,\rho _{m}}\circ M_{X\beta ,\lambda
_{n}^{1/2}(\Sigma (\rho _{m}))\sigma }\right) \circ \mathcal{I}_{0,\zeta
_{e}}\rightarrow \delta _{e}
\end{equation*}%
weakly as $m\rightarrow \infty $. Since $\varphi $ is bounded and is
continuous at $e$, the claim then follows from a version of the Portmanteau
theorem, cf. Theorem 30.12 in \cite{Bauer}. $\blacksquare $

\textbf{Proof of Proposition \ref{prop_bound}:} 1. Because $\emptyset \neq
\Phi \neq \mathbb{R}^{n}$ we can find $y_{0}\in \Phi $ and $y_{1}\notin \Phi 
$. By $G_{X}$-invariance we have that $\gamma y_{0}+X\theta \in \Phi $ and $%
\gamma y_{1}+X\theta \notin \Phi $ for every $\gamma \neq 0$ and for every $%
\theta \in \mathbb{R}^{k}$. Letting $\gamma $ converge to zero we see that $%
X\theta $ belongs to the closure of $\Phi $ as well as of its complement.
Thus $X\theta \in \limfunc{bd}(\Phi )$ holds for every $\theta $.

2. Suppose $y$ is an element of the boundary of the rejection region. If $%
y\in \limfunc{span}(X)$ there is nothing to prove. Hence assume $y\notin 
\limfunc{span}(X)$. If $T(y)\neq \kappa $ would hold, then by the continuity
assumption $y$ would be either in the interior or the exterior (i.e., the
complement of the closure) of the rejection region.

3. Because $T_{B}$ is continuous on $\mathbb{R}^{n}\backslash \limfunc{span}%
\left( X\right) $, Part 2 of the proposition establishes that the l.h.s. of (%
\ref{char_bd}) is contained in the r.h.s. Because of Part 1, it suffices to
show that every $y_{0}\notin \limfunc{span}\left( X\right) $ satisfying $%
T_{B}\left( y_{0}\right) =\kappa $ belongs to $\limfunc{bd}\left( \Phi
_{B,\kappa }\right) $. Obviously, $y_{0}\notin \Phi _{B,\kappa }$. It
remains to show that $y_{0}$ can be approximated by a sequence of elements
belonging to $\Phi _{B,\kappa }$: For $\lambda \in \mathbb{R}$ set $y\left(
\lambda \right) =y_{0}+\lambda y_{\ast }$ where $y_{\ast }\in \mathbb{R}^{n}$
is such that $T_{B}\left( y_{\ast }\right) >\kappa $. Such an $y_{\ast }$
exists, because $\Phi _{B,\kappa }\neq \emptyset $ by assumption.
Furthermore, $y_{\ast }\notin \limfunc{span}\left( X\right) $ must hold,
since otherwise $\lambda _{1}\left( B\right) =T_{B}\left( y_{\ast }\right)
>\kappa $ would follow, which in turn would entail $T_{B}\left( y\right)
\geq \lambda _{1}\left( B\right) >\kappa $ for all $y\in \mathbb{R}^{n}$,
i.e., $\Phi _{B,\kappa }=\mathbb{R}^{n}$, contradicting the assumptions. Set 
$A=C_{X}^{\prime }\left( B-\kappa I_{n-k}\right) C_{X}$ and note that $%
y_{0}^{\prime }Ay_{0}=0$ and $y_{\ast }^{\prime }Ay_{\ast }>0$ hold. Now%
\begin{equation*}
y\left( \lambda \right) ^{\prime }Ay\left( \lambda \right) =y_{0}^{\prime
}Ay_{0}+2\lambda y_{0}^{\prime }Ay_{\ast }+\lambda ^{2}y_{\ast }^{\prime
}Ay_{\ast }=2\lambda y_{0}^{\prime }Ay_{\ast }+\lambda ^{2}y_{\ast }^{\prime
}Ay_{\ast }\text{.}
\end{equation*}%
Choose a sequence $\lambda _{m}$ that converges to zero for $m\rightarrow
\infty $ and satisfies $\lambda _{m}>0$ for all $m$ if $y_{0}^{\prime
}Ay_{\ast }\geq 0$ and $\lambda _{m}<0$ for all $m$ if $y_{0}^{\prime
}Ay_{\ast }<0$. Then $y\left( \lambda _{m}\right) $ converges to $y_{0}$ and 
$y\left( \lambda _{m}\right) \notin \limfunc{span}\left( X\right) $ holds
for large enough $m$. Furthermore, we have $y^{\prime }\left( \lambda
_{m}\right) Ay\left( \lambda _{m}\right) >0$. But this means that $%
T_{B}\left( y\left( \lambda _{m}\right) \right) >\kappa $ holds for large $m$%
. $\blacksquare $

\textbf{Proof of Lemma \ref{AR1}: }Suppose Assumption \ref{ASCII} holds.
Then clearly 
\begin{equation}
\lim_{\rho \rightarrow a}c^{2}(\rho )\Pi _{\text{$\limfunc{span}$}(e)^{\bot
}}\Sigma (\rho )\Pi _{\text{$\limfunc{span}$}(e)^{\bot }}=\Lambda \Lambda
^{\prime }  \label{eq_79}
\end{equation}%
holds. Set $V=\Lambda \Lambda ^{\prime }$. Furthermore, the above relation
clearly implies $Ve=0$ and hence $\limfunc{span}(e)\subseteq \ker (V)$.
Because $Vy=0$ if and only if $\Lambda ^{\prime }y=0$, and because $\func{%
rank}(\Lambda ^{\prime })=\func{rank}(\Lambda )=$ $n-1$, it follows that $%
\ker (V)$ must be one-dimensional. Hence $\ker (V)=\limfunc{span}(e)$ must
hold. Since $V$ maps $\limfunc{span}(e)^{\bot }$ into $\limfunc{span}%
(e)^{\bot }$ in view of (\ref{eq_79}), it follows that $V$ is injective on $%
\limfunc{span}(e)^{\bot }$. To prove the converse, note that $V$ given by (%
\ref{scaled_limit_2})\ is by construction a bijection from $\limfunc{span}%
(e)^{\bot }$ to itself and is symmetric and nonnegative definite. Thus its
symmetric nonnegative definite square root $V^{1/2}$ exists and is a
bijective map from $\limfunc{span}(e)^{\bot }$ to itself. Furthermore, the
symmetric nonnegative square root of $c^{2}(\rho )\Pi _{\limfunc{span}%
(e)^{\bot }}\Sigma (\rho )\Pi _{\limfunc{span}(e)^{\bot }}$ can be written
in the form $c(\rho )\Pi _{\text{$\limfunc{span}$}(e)^{\bot }}\Sigma
^{1/2}(\rho )U(\rho )$ for a suitable choice of an orthogonal matrix $U(\rho
)$. By continuity of the symmetric nonnegative square root we obtain 
\begin{equation*}
c(\rho )\Pi _{\text{$\limfunc{span}$}(e)^{\bot }}\Sigma ^{1/2}(\rho )U(\rho
)\rightarrow V^{1/2}.
\end{equation*}%
It remains to set $\Lambda =V^{1/2}$ and $L_{\ast }(\rho )=\Sigma
^{1/2}(\rho )U(\rho )$. $\blacksquare $

\textbf{Proof of Theorem \ref{BT}:} A.1. By $G_{X}$-invariance of $\varphi $
and Assumption \ref{ASDR} the power function does neither depend on $\beta $
nor $\sigma $ (cf. Remark \ref{invariance}), and thus it suffices to
consider the case $\beta =0$ and $\sigma =1$. By Assumption \ref{ASDR} we
furthermore have%
\begin{equation}
E_{0,1,\rho }(\varphi )=\int_{\mathbb{R}^{n}}\varphi dP_{0,1,\rho }=\int
\varphi (L(\rho )\mathbf{z})d\Pr =\int \varphi (L_{\ast }(\rho )U\left( \rho
\right) \mathbf{z})d\Pr  \label{expect_1}
\end{equation}%
where $U\left( \rho \right) =L_{\ast }^{-1}(\rho )L(\rho )$ is an orthogonal
matrix. Observe that $\varphi (y+\gamma e)=\varphi (y)$ holds for every $y$
and for every $\gamma \in \mathbb{R}$: This is trivial for $\gamma =0$ and
follows for $\gamma \neq 0$ from%
\begin{equation}
\varphi (y+\gamma e)=\varphi (\gamma ^{-1}y+e)=\varphi (\gamma
^{-1}y)=\varphi (y),  \label{identity}
\end{equation}%
where we have made use of $G_{X}$-invariance of $\varphi $ as well as of (%
\ref{INV}). Observing that $\Pi _{\text{$\limfunc{span}$}(e)}L_{\ast }(\rho
)U\left( \rho \right) \mathbf{z}$ as well as $\Pi _{\text{$\limfunc{span}$}%
(e)}U\left( \rho \right) \mathbf{z}$ belong to $\limfunc{span}$$(e)$, using
relation (\ref{identity}) as well as $G_{X}$-invariance of $\varphi $ leads
to%
\begin{eqnarray}
\int \varphi (L_{\ast }(\rho )U\left( \rho \right) \mathbf{z})d\Pr &=&\int
\varphi \left( \Pi _{\text{$\limfunc{span}$}(e)^{\bot }}L_{\ast }(\rho
)U\left( \rho \right) \mathbf{z}\right) d\Pr  \notag \\
&=&\int \varphi \left( c(\rho )\Pi _{\text{$\limfunc{span}$}(e)^{\bot
}}L_{\ast }(\rho )U\left( \rho \right) \mathbf{z}\right) d\Pr  \notag \\
&=&\int \varphi \left( A\left( \rho \right) U\left( \rho \right) \mathbf{z}%
\right) d\Pr ,  \label{expect_2}
\end{eqnarray}%
where $A\left( \rho \right) $ is shorthand for $\Pi _{\text{$\limfunc{span}$}%
(e)}+c(\rho )\Pi _{\text{$\limfunc{span}$}(e)^{\bot }}L_{\ast }(\rho )$.
Since the image of $\Lambda $ is $\limfunc{span}(e)^{\bot }$ and $\Lambda $
is injective when restricted to $\limfunc{span}(e)^{\bot }$ it follows that $%
A:=\Pi _{\limfunc{span}(e)}+\Lambda $ is bijective as a map from $\mathbb{R}%
^{n}$ to $\mathbb{R}^{n}$. [To see this suppose that $Ay=0$. Because $%
\Lambda y\in \limfunc{span}(e)^{\bot }$ this implies $\Pi _{\limfunc{span}%
(e)}y=0$ as well as $\Lambda y=0$. The first equality now implies $y\in 
\limfunc{span}(e)^{\bot }$. Bijectivity of $\Lambda $ on $\limfunc{span}%
(e)^{\bot }$ then implies $y=0$.] By Assumption \ref{ASCII} the matrix $%
A\left( \rho \right) $ converges to $A$ for $\rho \rightarrow a$ and thus $%
A\left( \rho \right) $ is bijective as a map from $\mathbb{R}^{n}$ to $%
\mathbb{R}^{n}$ whenever $\rho $ is sufficiently close to $a$, say $\rho
\geq \rho _{0}$. If now $\omega $ is an accumulation point of $E_{0,1,\rho
}(\varphi )$, we can find a sequence $\rho _{m}$ that converges to $a$ such
that $E_{0,1,\rho _{m}}(\varphi )$ converges to $\omega $. By passing to a
suitable subsequence, we may also assume that $U\left( \rho _{m}\right) $
converges to an orthogonal matrix $U$, say. W.l.o.g. we may furthermore
assume that $\rho _{m}\geq \rho _{0}$ holds and thus $A\left( \rho
_{m}\right) $ is nonsingular. By the transformation formula for densities
the $\mu _{\mathbb{R}^{n}}$-density of the random vector $A\left( \rho
_{m}\right) U(\rho _{m})\mathbf{z}$ is given by%
\begin{equation*}
\left\vert \det \left( A^{-1}\left( \rho _{m}\right) \right) \right\vert
p\left( U^{\prime }\left( \rho _{m}\right) A^{-1}\left( \rho _{m}\right)
y\right) .
\end{equation*}%
Because of $A\left( \rho _{m}\right) \rightarrow A$, $U\left( \rho
_{m}\right) \rightarrow U$, and because $p$ is continuous $\mu _{\mathbb{R}%
^{n}}$-almost everywhere, this expression converges for $\mu _{\mathbb{R}%
^{n}}$-almost every $y\in \mathbb{R}^{n}$ to 
\begin{equation}
\left\vert \det \left( A^{-1}\right) \right\vert p\left( U^{\prime
}A^{-1}y\right)  \label{dens}
\end{equation}%
as $\rho \rightarrow a$, which is the density of the random vector $AU%
\mathbf{z}$. Scheff{\'{e}}'s lemma thus implies that the distribution of $%
A\left( \rho _{m}\right) U(\rho _{m})\mathbf{z}$ converges in total
variation norm to $Q_{A,U}$, the distribution of $AU\mathbf{z}$. It now
follows in view of (\ref{expect_1}) and (\ref{expect_2})\ that 
\begin{equation*}
E_{0,1,\rho _{m}}(\varphi )\rightarrow E_{Q_{A,U}}\left( \varphi \right)
=\int \varphi dQ_{A,U}=\int \varphi \left( AU\mathbf{z}\right) d\Pr .
\end{equation*}%
Now%
\begin{equation*}
\varphi \left( AU\mathbf{z}\right) =\varphi \left( \Pi _{\limfunc{span}(e)}U%
\mathbf{z}+\Lambda U\mathbf{z}\right) =\varphi \left( \Lambda U\mathbf{z}%
\right)
\end{equation*}%
holds because of (\ref{identity}), implying $E_{Q_{A,U}}\left( \varphi
\right) =E_{Q_{\Lambda ,U}}\left( \varphi \right) $. This shows that $\omega
=E_{Q_{\Lambda ,U}}\left( \varphi \right) $ must hold. Conversely, given $%
U\in \mathcal{U}\left( L_{\ast }^{-1}L\right) $ we can find a sequence $\rho
_{m}\rightarrow a$ such that $U\left( \rho _{m}\right) =L_{\ast }^{-1}(\rho
_{m})L(\rho _{m})$ converges to the given $U$. Repeating the argument given
above then shows that $E_{Q_{\Lambda ,U}}\left( \varphi \right) $, for the
given $U$, arises as an accumulation point of $E_{0,1,\rho }(\varphi )$ for $%
\rho \rightarrow a$.

A.2. The claim follows immediately from the already established Part 1.

A.3. Recall that $E_{Q_{\Lambda ,U}}\left( \varphi \right)
=E_{Q_{A,U}}\left( \varphi \right) $. Hence,%
\begin{equation*}
E_{Q_{\Lambda ,U}}\left( \varphi \right) =\int_{\mathbb{R}^{n}}\varphi
\left( AUy\right) p\left( y\right) dy=\int_{\mathbb{R}^{n}\backslash \left\{
0\right\} }\varphi \left( AUy\right) p\left( y\right) dy=\int_{\left(
0,\infty \right) \times S^{n-1}}\varphi \left( rAUs\right) p\left( rs\right)
dH\left( r,s\right)
\end{equation*}%
where $H$ is the pushforward measure of $\mu _{\mathbb{R}^{n}}$ (restricted
to $\mathbb{R}^{n}\backslash \left\{ 0\right\} $) under the map $y\mapsto
\left( \left\Vert y\right\Vert ,y/\left\Vert y\right\Vert \right) $. Now $H$
is nothing else than the product of the measure on $\left( 0,\infty \right) $
with density $r^{n-1}$ and the surface measure $c\upsilon _{S^{n-1}}$ on $%
S^{n-1}$ with the constant $c$ given by $2\pi ^{n/2}/\Gamma \left(
n/2\right) $ (cf. \cite{Stroock}). In view of Fubini's theorem (observe all
functions involved are nonnegative) and invariance of $\varphi $ we then
obtain%
\begin{equation*}
E_{Q_{\Lambda ,U}}\left( \varphi \right) =c\int_{S^{n-1}}\varphi \left(
AUs\right) \left( \int_{\left( 0,\infty \right) }p_{s}\left( r\right)
r^{n-1}dr\right) d\upsilon _{S^{n-1}}.
\end{equation*}%
If $\varphi \left( \cdot \right) $ is not equal to zero $\mu _{\mathbb{R}%
^{n}}$-almost everywhere, then so is $\varphi \left( AU\cdot \right) $
because $AU$ is nonsingular. Now scale invariance of $\varphi $ translates
into scale invariance of $\varphi \left( AU\cdot \right) $, and hence $%
\varphi \left( AU\cdot \right) $ restricted to $S^{n-1}$ is not equal to
zero $\upsilon _{S^{n-1}}$-almost everywhere, cf. Remark \ref{E1}(i) in
Appendix \ref{auxil}. Since the inner integral in the preceding display is
positive $\upsilon _{S^{n-1}}$-almost everywhere by the assumption on $p$,
we conclude that $E_{Q_{\Lambda ,U}}\left( \varphi \right) $ must be
positive. The claim that $E_{Q_{\Lambda ,U}}\left( \varphi \right) <1$ is
proved by applying the above to $1-\varphi $. Hence, if $\varphi $ is
neither $\mu _{\mathbb{R}^{n}}$-almost everywhere equal to zero nor $\mu _{%
\mathbb{R}^{n}}$-almost everywhere equal to one, we have established that $%
E_{Q_{\Lambda ,U}}\left( \varphi \right) $ is strictly between $0$ and $1$.
Next observe that $\mathcal{U}\left( L_{\ast }^{-1}L\right) $ is a compact
set. It thus suffices to establish that the map $U\rightarrow
E_{Q_{A,U}}\left( \varphi \right) $ is continuous on $\mathcal{U}\left(
L_{\ast }^{-1}L\right) $. But this follows from (\ref{dens}), $\mu _{\mathbb{%
R}^{n}}$-almost sure continuity of $p$, and Scheff\'{e}'s Lemma.

B. By the assumptions on $\mathfrak{P}$ the random vector $\mathbf{z}$ is
spherically symmetric with $\Pr \left( \mathbf{z}=0\right) =0$, and hence is
almost surely equal to $\mathbf{r}\mathbf{E}$ where $\mathbf{r}=\left\Vert 
\mathbf{z}\right\Vert $ is a random variable satisfying $\Pr (\mathbf{r}%
>0)=1 $ and where $\mathbf{E}=\mathbf{z}/\left\Vert \mathbf{z}\right\Vert $
is independent of $\mathbf{r}$ and is uniformly distributed on the unit
sphere $S^{n-1}$ (cf. Lemma 1 in \cite{Cambanis}). Possibly after enlarging
the underlying probability space we can find a random variable $\mathbf{r}%
_{0}$ which is independent of $\mathbf{E}$ and which is distributed as the
square root of a chi-square with $n$ degrees of freedom. By $G_{X}$%
-invariance of $\varphi $ we have 
\begin{eqnarray}
E_{0,1,\rho }(\varphi ) &=&\int \varphi (L(\rho )\mathbf{z})d\Pr =\int
\varphi (L(\rho )\mathbf{r}\mathbf{E})d\Pr =\int \varphi (L(\rho )\mathbf{E}%
)d\Pr  \notag \\
&=&\int \varphi (L(\rho )\mathbf{r}_{0}\mathbf{E})d\Pr =\int \varphi (L(\rho
)\mathbf{G})d\Pr  \label{reduction}
\end{eqnarray}%
where $\mathbf{G}=\mathbf{r}_{0}\mathbf{E}$ has a standard multivariate
Gaussian distribution. Again using $G_{X}$-invariance of $\varphi $ we
similarly obtain%
\begin{equation*}
E_{Q_{\Lambda ,U}}\left( \varphi \right) =E\varphi \left( \Lambda U\mathbf{z}%
\right) =E\varphi \left( \Lambda U\mathbf{r}\mathbf{E}\right) =E\varphi
\left( \Lambda U\mathbf{E}\right) =E\varphi \left( \Lambda U\mathbf{r}_{0}%
\mathbf{E}\right) =E\varphi \left( \Lambda U\mathbf{G}\right) =E_{Q_{\Lambda
,U}^{0}}\left( \varphi \right)
\end{equation*}%
where $Q_{\Lambda ,U}^{0}$ denotes the distribution of $\Lambda U\mathbf{G}$%
. This shows that we may act as if $\mathbf{z}$ were Gaussian. Consequently,
the results in A.1-A.3 apply. Furthermore, under elliptical symmetry $%
Q_{\Lambda ,U}=Q_{\Lambda ,I_{n}}$ holds for every orthogonal matrix $U$.
Hence, there exists only one accumulation point which is given by $%
E_{Q_{\Lambda ,I_{n}}}\left( \varphi \right) $. [Alternatively, under
elliptical symmetry we may choose w.l.o.g. $L(\cdot )$ to be any square root
of $\Sigma \left( \cdot \right) $, and thus equal to $L_{\ast }(\cdot )$,
and then apply Part A.2.] That $0<E_{Q,I_{n}}\left( \varphi \right) <1$
holds under the additional assumption on $\varphi $ follows from Part A.3. $%
\blacksquare $

\begin{lemma}
\label{off-diag} Suppose Assumptions \ref{ASC} and \ref{ASCII} hold with the
same vector $e$. Then 
\begin{equation*}
\Xi \left( \rho \right) :=\lambda _{n}^{-1/2}\left( \Sigma \left( \rho
\right) \right) c\left( \rho \right) \Pi _{\limfunc{span}\left( e\right)
^{\bot }}\Sigma \left( \rho \right) \Pi _{\limfunc{span}\left( e\right) }
\end{equation*}%
is bounded for $\rho \rightarrow a$ and the set of all accumulation points
of $\Xi \left( \rho \right) $ for $\rho \rightarrow a$ is given by 
\begin{equation*}
\left\{ \Lambda U_{0}^{\prime }ee^{\prime }:U_{0}\in \mathcal{U}\left(
\Sigma ^{-1/2}L_{\ast }\right) \right\} .
\end{equation*}%
The same statements hold if $\Xi \left( \rho \right) $ is replaced by $\Xi
_{1}\left( \rho \right) :=\lambda _{n}^{-1/2}\left( \Sigma \left( \rho
\right) \right) c\left( \rho \right) \Pi _{\limfunc{span}\left( e\right)
^{\bot }}\Sigma \left( \rho \right) $ or $\Xi _{2}\left( \rho \right)
:=c\left( \rho \right) \Pi _{\limfunc{span}\left( e\right) ^{\bot }}\Sigma
^{1/2}\left( \rho \right) \Pi _{\limfunc{span}\left( e\right) }$.
\end{lemma}

\textbf{Proof: }Rewrite $\Xi \left( \rho \right) $ as $A_{1}\left( \rho
\right) U^{\prime }\left( \rho \right) A_{2}\left( \rho \right) \Pi _{%
\limfunc{span}\left( e\right) }$ where $A_{1}\left( \rho \right) =c\left(
\rho \right) \Pi _{\limfunc{span}\left( e\right) ^{\bot }}L_{\ast }\left(
\rho \right) $, $U\left( \rho \right) =\Sigma ^{-1/2}\left( \rho \right)
L_{\ast }\left( \rho \right) $ is orthogonal, and $A_{2}\left( \rho \right)
=\lambda _{n}^{-1/2}\left( \Sigma \left( \rho \right) \right) \Sigma
^{1/2}\left( \rho \right) $. Now $A_{1}\left( \rho \right) $ and $%
A_{2}\left( \rho \right) $ converge to $\Lambda $ and $ee^{\prime }$,
respectively, by Assumptions \ref{ASC} and \ref{ASCII}. Since $U\left( \rho
\right) $ is clearly bounded, boundedness of $\Xi \left( \rho \right) $
follows. The claim concerning the set of accumulation points also now
follows immediately. The proofs for $\Xi _{1}$ and $\Xi _{2}$ are completely
analogous. $\blacksquare $

\textbf{Proof of Theorem \ref{BT2}:} 1. Using invariance w.r.t. $G_{X}$,
Equation (\ref{ATaylor}), and homogeneity of $D$ we obtain for every $\gamma
\neq 0$%
\begin{equation}
T(\gamma e+h)=T(e+\gamma ^{-1}h)=T(e)+\gamma ^{-q}D(h)+R(\gamma ^{-1}h)
\label{Taylor}
\end{equation}%
for every $h\in \mathbb{R}^{n}$. Let $\omega $ be an accumulation point of $%
P_{\beta ,\sigma ,\rho }\left( \left\{ y\in \mathbb{R}^{n}:T(y)>T(e)\right\}
\right) $ for $\rho \rightarrow a$. Then we can find a sequence $\rho
_{m}\in \lbrack 0,a)$ with $\rho _{m}\rightarrow a$ along which the
rejection probability converges to $\omega $. W.l.o.g. (possibly after
passing to a suitable subsequence) we may also assume that along this
sequence the orthogonal matrices $U\left( \rho _{m}\right) =L_{\ast
}^{-1}(\rho _{m})L(\rho _{m})$ and $U_{0}\left( \rho _{m}\right) =\Sigma
^{-1/2}(\rho _{m})L(\rho _{m})$ converge to orthogonal matrices $U$ and $%
U_{0}$, respectively. Using $\Pi _{\limfunc{span}(e)}=ee^{\prime }$ and
invariance w.r.t. $G_{X}$ we obtain 
\begin{equation*}
T\left( X\beta +\sigma L(\rho _{m})\mathbf{z}\right) =T\left( ee^{\prime
}L(\rho _{m})\mathbf{z}+\Pi _{\text{$\limfunc{span}$}(e)^{\bot }}L(\rho _{m})%
\mathbf{z}\right) .
\end{equation*}%
Observe that $e^{\prime }L(\rho _{m})\mathbf{z}$ is nonzero with probability 
$1$ because $e\neq 0$, $L(\rho _{m})$ is nonsingular, and $\mathbf{z}$
possesses a density. Hence, combining the previous display and equation (\ref%
{Taylor}) with $\gamma =e^{\prime }L(\rho _{m})\mathbf{z}$ and $h=\Pi _{%
\limfunc{span}(e)^{\bot }}L(\rho _{m})\mathbf{z}$ and then multiplying by $%
c^{q}(\rho _{m})\lambda ^{q/2}(m)$, where $\lambda (m)$ is shorthand for $%
\lambda _{n}(\Sigma (\rho _{m}))$, we obtain that%
\begin{eqnarray}
&&c^{q}(\rho _{m})\lambda ^{q/2}(m)\left( T\left( X\beta +\sigma L(\rho _{m})%
\mathbf{z}\right) -T(e)\right)  \notag \\
&=&\left( \lambda ^{-1/2}(m)e^{\prime }L(\rho _{m})\mathbf{z}\right)
^{-q}D\left( c(\rho _{m})\Pi _{\text{$\limfunc{span}$}(e)^{\bot }}L_{\ast
}(\rho _{m})U\left( \rho _{m}\right) \mathbf{z}\right)  \notag \\
&&+c^{q}(\rho _{m})\lambda ^{q/2}(m)R\left( \left( \lambda
^{-1/2}(m)e^{\prime }L(\rho _{m})\mathbf{z}\right) ^{-1}\lambda
^{-1/2}(m)\Pi _{\text{$\limfunc{span}$}(e)^{\bot }}L(\rho _{m})\mathbf{z}%
\right)  \label{rep}
\end{eqnarray}%
holds almost surely. Next observe that by Assumption \ref{ASC}, continuity
of the symmetric nonnegative definite square root, and $\left( ee^{\prime
}\right) ^{1/2}=ee^{\prime }$ we have%
\begin{equation}
\lambda ^{-1/2}(m)e^{\prime }L(\rho _{m})\mathbf{z}=\lambda
^{-1/2}(m)e^{\prime }\Sigma ^{1/2}(\rho _{m})U_{0}(\rho _{m})\mathbf{z}%
\rightarrow e^{\prime }U_{0}\mathbf{z}  \label{conv_1}
\end{equation}%
and 
\begin{equation}
\lambda ^{-1/2}(m)\Pi _{\text{$\limfunc{span}$}(e)^{\bot }}L(\rho _{m})%
\mathbf{z}=\lambda ^{-1/2}(m)\Pi _{\text{$\limfunc{span}$}(e)^{\bot }}\Sigma
^{1/2}(\rho _{m})U_{0}(\rho _{m})\mathbf{z}\rightarrow \Pi _{\text{$\limfunc{%
span}$}(e)^{\bot }}ee^{\prime }U_{0}\mathbf{z}=0,  \label{conv_2}
\end{equation}%
where the convergence holds for every realization of $\mathbf{z}$. Note that 
$e^{\prime }U_{0}\mathbf{z}\neq 0$ holds almost surely. Relation (\ref%
{conv_1}) together with Assumption \ref{ASCII} then implies that the first
term on the r.h.s. of (\ref{rep}) converges almost surely to $\left(
e^{\prime }U_{0}\mathbf{z}\right) ^{-q}D\left( \Lambda U\mathbf{z}\right) $
since $D$ is clearly continuous. We next show that the second term on the
r.h.s. of (\ref{rep}) converges to zero almost surely: Let $h_{m}$ denote
the argument of $R$ in (\ref{rep}). Fix a realization of $\mathbf{z}$ such
that $e^{\prime }U_{0}\mathbf{z}\neq 0$. Then $h_{m}$ is well-defined for
large enough $m$, and it converges to zero because of (\ref{conv_1}) and (%
\ref{conv_2}). Since $R(0)=0$ holds as a consequence of (\ref{ATaylor}), we
only need to consider subsequences along which $h_{m}\neq 0$. For notational
convenience we denote such subsequences again by $h_{m}$. Because of the
assumptions on $R$ it suffices to show that $c^{q}(\rho _{m})\lambda
^{q/2}(m)\left\Vert h_{m}\right\Vert ^{q}$ is bounded. Now%
\begin{eqnarray*}
&&c^{q}(\rho _{m})\lambda ^{q/2}(m)\left\Vert h_{m}\right\Vert ^{q} \\
&=&\left\Vert \left( \lambda ^{-1/2}(m)e^{\prime }L(\rho _{m})\mathbf{z}%
\right) ^{-1}c(\rho _{m})\Pi _{\text{$\limfunc{span}$}(e)^{\bot }}L(\rho
_{m})\mathbf{z}\right\Vert ^{q} \\
&=&\left\Vert \left( \lambda ^{-1/2}(m)e^{\prime }L(\rho _{m})\mathbf{z}%
\right) ^{-1}c(\rho _{m})\Pi _{\text{$\limfunc{span}$}(e)^{\bot }}L_{\ast
}(\rho _{m})U\left( \rho _{m}\right) \mathbf{z}\right\Vert ^{q} \\
&\rightarrow &\left\Vert \left( e^{\prime }U_{0}\mathbf{z}\right)
^{-1}\Lambda U\mathbf{z}\right\Vert ^{q}<\infty ,
\end{eqnarray*}%
where we have made use of (\ref{conv_1}) and Assumption \ref{ASCII}. We have
thus established that 
\begin{equation}
c^{q}(\rho _{m})\lambda ^{q/2}(m)\left( T\left( X\beta +\sigma L(\rho _{m})%
\mathbf{z}\right) -T(e)\right) \rightarrow \left( e^{\prime }U_{0}\mathbf{z}%
\right) ^{-q}D\left( \Lambda U\mathbf{z}\right)  \label{conv_3}
\end{equation}%
almost surely. Note that the range of $\Lambda $ is $\limfunc{span}(e)^{\bot
}$, and that $\Lambda $ is bijective as a map from $\limfunc{span}(e)^{\bot
} $ to itself. Hence, the random variable $\Lambda U\mathbf{z}$ takes its
values in $\limfunc{span}(e)^{\bot }$ and possesses a density on this
subspace (w.r.t. $n-1$ dimensional Lebesgue measure on this subspace). Since 
$D$ restricted to $\limfunc{span}(e)^{\bot }$ can be expressed as a
multivariate polynomial (in $n-1$ variables) and does not vanish identically
on $\limfunc{span}(e)^{\bot }$, it vanishes at most on a subset of $\limfunc{%
span}(e)^{\bot }$ that has $\left( n-1\right) $-dimensional Lebesgue measure
zero. It follows that $D\left( \Lambda U\mathbf{z}\right) $, and hence the
limit in (\ref{conv_3}), is nonzero almost surely. Observe that 
\begin{eqnarray*}
P_{\beta ,\sigma ,\rho _{m}}\left( \left\{ y\in \mathbb{R}%
^{n}:T(y)>T(e)\right\} \right) &=&\Pr \left( T\left( X\beta +\sigma L\left(
\rho _{m}\right) \mathbf{z}\right) >T\left( e\right) \right) \\
&=&\Pr \left( c^{q}(\rho _{m})\lambda ^{q/2}(m)\left( T\left( X\beta +\sigma
L(\rho _{m})\mathbf{z}\right) -T(e)\right) >0\right)
\end{eqnarray*}%
since $c(\rho _{m})$ and $\lambda (m)$ are positive. By an application of
the Portmanteau theorem we can thus conclude from (\ref{conv_3}) that for $%
m\rightarrow \infty $%
\begin{equation}
P_{\beta ,\sigma ,\rho _{m}}\left( \left\{ y\in \mathbb{R}%
^{n}:T(y)>T(e)\right\} \right) \rightarrow \Pr (\left( e^{\prime }U_{0}%
\mathbf{z}\right) ^{-q}D(\Lambda U\mathbf{z})>0).  \label{accu}
\end{equation}%
The limit in the preceding display obviously reduces to (\ref{even}) and (%
\ref{odd}), respectively, and clearly $\left( U,U_{0}\right) \in \mathcal{U}%
\left( L_{\ast }^{-1}L,\Sigma ^{-1/2}L\right) $ implies $U\in \mathcal{U}%
\left( L_{\ast }^{-1}L\right) $. This together then proves that every
accumulation point $\omega $ has the claimed form. To prove the converse,
observe first that for every $U\in \mathcal{U}\left( L_{\ast }^{-1}L\right) $
we can find an $U_{0}$ such that $\left( U,U_{0}\right) \in \mathcal{U}%
\left( L_{\ast }^{-1}L,\Sigma ^{-1/2}L\right) $ holds (exploiting
compactness of the set of orthogonal matrices). Now, let $\left(
U,U_{0}\right) \in \mathcal{U}\left( L_{\ast }^{-1}L,\Sigma ^{-1/2}L\right) $
be given. Then we can find a sequence $\rho _{m}\in \lbrack 0,a)$ with $\rho
_{m}\rightarrow a$ such that $U\left( \rho _{m}\right) =L_{\ast }^{-1}(\rho
_{m})L(\rho _{m})$ and $U_{0}\left( \rho _{m}\right) =\Sigma ^{-1/2}(\rho
_{m})L(\rho _{m})$ converge to $U$ and $U_{0}$, respectively. Repeating the
preceding arguments, then shows that $\Pr (\left( e^{\prime }U_{0}\mathbf{z}%
\right) ^{-q}D(\Lambda U\mathbf{z})>0)$ is the limit of $P_{\beta ,\sigma
,\rho _{m}}\left( \left\{ y\in \mathbb{R}^{n}:T(y)>T(e)\right\} \right) $.
The final claim is now obvious.

2. If $\mathfrak{P}$ is an elliptically symmetric family we can w.l.o.g. set 
$L(\cdot )=L_{\ast }(\cdot )$, implying that $\mathcal{U}\left( L_{\ast
}^{-1}L,\Sigma ^{-1/2}L\right) $ reduces to $\left\{ I_{n}\right\} \times 
\mathcal{U}\left( \Sigma ^{-1/2}L_{\ast }\right) $. Furthermore, as $\mathbf{%
z}$ is then spherically symmetric and satisfies $\Pr (\mathbf{z}=0)=0$, it
is almost surely equal to $\mathbf{r}\mathbf{E}$ where $\mathbf{r}$ must
satisfy $\Pr (\mathbf{r}>0)=1$ and where $\mathbf{E}$ is independent of $%
\mathbf{r}$ and is uniformly distributed on the unit sphere in $\mathbb{R}%
^{n}$. Let $\mathbf{r}_{0}$ be a random variable which is independent of $%
\mathbf{E}$ and which is distributed as the square root of a chi-square with 
$n$ degrees of freedom (this may require enlarging the underlying
probability space) and define $\mathbf{G=r}_{0}\mathbf{E}$ which clearly is
a multivariate Gaussian random vector with mean zero and covariance matrix $%
I_{n}$. Define $\mathfrak{P}_{0}$ in the same way as $\mathfrak{P}$, but
with $\mathbf{G}$ replacing $\mathbf{z}$ in Assumption \ref{ASDR}. Observe
that the rejection probabilities of the test considered are the same whether
they are calculated under the experiment $\mathfrak{P}$ or $\mathfrak{P}_{0}$
because of $G_{X}$-invariance of the test statistic. Applying the already
established Part 1 in the context of the experiment $\mathfrak{P}_{0}$ thus
shows that the accumulation points of the rejection probabilities calculated
under $\mathfrak{P}_{0}$ as well as under $\mathfrak{P}$ equal $\Pr
(D(\Lambda \mathbf{G})>0)$ for even $q$ and equal $\Pr (D(\Lambda \mathbf{G}%
)>0,e^{\prime }U_{0}\mathbf{G}>0)+\Pr (D(\Lambda \mathbf{G})<0,e^{\prime
}U_{0}\mathbf{G}<0)$ for odd $q$. In view of homogeneity of $D$ and the fact
that $\mathbf{r}$ as well as $\mathbf{r}_{0}$ are almost surely positive,
these probabilities do not change their value if we replace $\mathbf{G}$ by $%
\mathbf{z}$. This proves (\ref{even2}) and (\ref{odd2}). To prove the last
but one claim observe that $E\left( \left( e^{\prime }U_{0}\mathbf{G}\right)
\Lambda \mathbf{G}\right) =\Lambda U_{0}^{\prime }e=0$. Consequently, $%
e^{\prime }U_{0}\mathbf{G}$ and $\Lambda \mathbf{G}$ are independent. Hence
the accumulation point can be written as 
\begin{equation*}
\Pr (D(\Lambda \mathbf{G})>0)\Pr (e^{\prime }U_{0}\mathbf{G}>0)+\Pr
(D(\Lambda \mathbf{G})<0)\Pr (e^{\prime }U_{0}\mathbf{G}<0).
\end{equation*}%
This reduces to $1/2$, because then obviously $\Pr (e^{\prime }U_{0}\mathbf{G%
}>0)=\Pr (e^{\prime }U_{0}\mathbf{G}<0)=1/2$ (note that $\Pr (e^{\prime
}U_{0}\mathbf{G}=0)=0$) and because $\Pr (D(\Lambda \mathbf{G})=0)=0$ (which
is proved by arguments similar to the ones given below (\ref{conv_3})). The
final claim follows because by the assumed symmetry $\Lambda U_{0}^{\prime
}e=U_{0}\Lambda ^{\prime }e=0$, the last equality following from the
definition of $\Lambda $.

3. Lemma \ref{off-diag} shows that under the additional assumption we have $%
\Lambda U_{0}^{\prime }ee^{\prime }=0$ for every $U_{0}\in \mathcal{U}\left(
\Sigma ^{-1/2}L_{\ast }\right) $, and hence $\Lambda U_{0}^{\prime }e=0$.
The claim then follows from Part 2. $\blacksquare $

\begin{lemma}
\label{T}Suppose $T$ is a test statistic that satisfies the conditions
imposed on $T$ in Theorem \ref{BT2} for some normalized vector $e$. Then:

\begin{enumerate}
\item $D\left( h\right) =D\left( \Pi _{\limfunc{span}(e)^{\bot }}h\right) $
holds for every $h\in \mathbb{R}^{n}$. In particular, $D$ vanishes on all of 
$\limfunc{span}(e)$.

\item If $D\left( h\right) <0$ holds for every $h\in \limfunc{span}(e)^{\bot
}$ with $h\neq 0$, then there exists a neighborhood of $e$ in $\mathbb{R}%
^{n} $ such that $T\left( y\right) \leq T\left( e\right) $ holds for every $%
y $ in that neighborhood.
\end{enumerate}
\end{lemma}

\textbf{Proof: }1. Write $h$ as $\gamma e+h_{2}$ with $h_{2}=\Pi _{\limfunc{%
span}(e)^{\bot }}h$. Then for every sufficiently small real $c>0$ we have $%
1+c\gamma \neq 0$, and hence exploiting $G_{X}$-invariance of $T$ we obtain%
\begin{equation*}
T\left( e+ch\right) =T\left( \left( 1+c\gamma \right) e+ch_{2}\right)
=T\left( e+\left( 1+c\gamma \right) ^{-1}ch_{2}\right) .
\end{equation*}%
Applying (\ref{ATaylor}) to both sides of the above equation, using
homogeneity of $D$, and dividing by $c^{-q}$ we arrive at%
\begin{equation*}
D\left( h\right) +c^{-q}R\left( ch\right) =\left( 1+c\gamma \right)
^{-q}D\left( h_{2}\right) +c^{-q}R\left( \left( 1+c\gamma \right)
^{-1}ch_{2}\right) .
\end{equation*}%
Now observe that $c^{-q}R\left( ch\right) $ is zero for $h=0$, and converges
to zero for $c\rightarrow 0$ for $h\neq 0$. A similar statement holds for $%
c^{-q}R\left( \left( 1+c\gamma \right) ^{-1}ch_{2}\right) $ as well. Since $%
1+c\gamma \rightarrow 1$, we obtain $D\left( h\right) =D\left( h_{2}\right) $
which proves the first claim. The second claim is then an immediate
consequence since $D\left( 0\right) =0$ by homogeneity.

2. Suppose the claim were false. We could then find a sequence $%
h_{m}\rightarrow 0$ with $T\left( e+h_{m}\right) >T\left( e\right) $.
Rewrite $h_{m}$ as $\gamma _{m}e+h_{m2}$ with $h_{m2}=\Pi _{\limfunc{span}%
(e)^{\bot }}h_{m}$. Clearly, $\gamma _{m}\rightarrow 0$ would have to hold,
implying $1+\gamma _{m}>0$ for all sufficiently large $m$. Using $G_{X}$%
-invariance we obtain $T\left( e+h_{m}\right) =T\left( e+\left( 1+\gamma
_{m}\right) ^{-1}h_{m2}\right) $ for all large $m$. In particular, we
conclude that $h_{m2}\neq 0$ would have to hold for all large $m$. Applying (%
\ref{ATaylor}) to the r.h.s. of the preceding equation we thus obtain for
all large $m$%
\begin{equation*}
0<T\left( e+h_{m}\right) -T\left( e\right) =D\left( \left( 1+\gamma
_{m}\right) ^{-1}h_{m2}\right) +R\left( \left( 1+\gamma _{m}\right)
^{-1}h_{m2}\right) .
\end{equation*}%
Using homogeneity of $D$ we then have for all large $m$%
\begin{equation*}
0<D\left( h_{m2}/\left\Vert h_{m2}\right\Vert \right) +R\left( \left(
1+\gamma _{m}\right) ^{-1}h_{m2}\right) /\left\Vert \left( 1+\gamma
_{m}\right) ^{-1}h_{m2}\right\Vert ^{q}=D\left( h_{m2}/\left\Vert
h_{m2}\right\Vert \right) +o\left( 1\right) .
\end{equation*}%
Note that $h_{m2}/\left\Vert h_{m2}\right\Vert $ is an element of the
compact set $S^{n-1}\cap \limfunc{span}(e)^{\bot }$ on which $D$ is
continuous and negative. Hence, the r.h.s. of the preceding display is
eventually bounded from above by zero, a contradiction. $\blacksquare $

Inspection of the proof of Part 1 of the preceding lemma shows that this
proof in fact does not make use of the property that $D$ does not vanish on
all of $\limfunc{span}(e)^{\bot }$.

\textbf{Proof of Corollary \ref{Lem_Illust_1}:} 1. Clearly $T_{B}(e)>\kappa
\geq \lambda _{1}\left( B\right) $ implies $e\notin \limfunc{span}(X)$ in
view of the definition of $T_{B}$. In view of the assumption on $\kappa $,
the rejection region satisfies $\emptyset \neq \Phi _{B,\kappa }\neq \mathbb{%
R}^{n}$. Consequently $T_{B}(e)>\kappa $ implies $e\notin \limfunc{bd}\left(
\Phi _{B,\kappa }\right) $, cf. Proposition \ref{prop_bound}. But $e\in \Phi
_{B,\kappa }$ clearly holds, implying that $e\in \limfunc{int}\left( \Phi
_{B,\kappa }\right) $. The result then follows immediately from Theorem \ref%
{MT} combined with the observation that $\mathbf{1}_{\Phi _{B,\kappa }}$ is
continuous at $e$ if and only if $e\notin \limfunc{bd}\left( \Phi _{B,\kappa
}\right) $.

2. Since $e\notin \limfunc{span}(X)$ by assumption, we conclude similarly as
above that $T_{B}(e)<\kappa $ implies $e\notin \limfunc{bd}\left( \Phi
_{B,\kappa }\right) $. But $e\notin \Phi _{B,\kappa }$ clearly holds,
implying that $e\notin \limfunc{cl}\left( \Phi _{B,\kappa }\right) $. As
before, the result then follows from Theorem \ref{MT}. $\blacksquare $

\textbf{Proof of Corollary \ref{Lem_Illust_2}:} Observe that (\ref{INV}) is
satisfied for $\mathbf{1}_{\Phi _{B,\kappa }}$ since $T_{B}$ is $G_{X}$%
-invariant and $e\in \limfunc{span}(X)$ by assumption. Hence, all
assumptions of Part B of Theorem \ref{BT} are satisfied and thus the
existence and the form of the limit follows. If $\kappa >\lambda _{1}\left(
B\right) $ the test $\mathbf{1}_{\Phi _{B,\kappa }}$ is neither $\mu _{%
\mathbb{R}^{n}}$-almost everywhere equal to zero nor $\mu _{\mathbb{R}^{n}}$%
-almost everywhere equal to one, whereas $\mathbf{1}_{\Phi _{B,\kappa }}$ is 
$\mu _{\mathbb{R}^{n}}$-almost everywhere equal to one if $\kappa =\lambda
_{1}\left( B\right) $ as discussed in Remark \ref{range_for_kappa}. Part B
of Theorem \ref{BT} and Remark \ref{RBT}(iv) then deliver the remaining
claims. $\blacksquare $

\textbf{Proof of Corollary \ref{Lem_Illust_3}:} All assumptions for Part 2
of Theorem \ref{BT2} (including the elliptic symmetry assumption) except for
(\ref{ATaylor}) are obviously satisfied. We first consider the situation of
Part 1 of the corollary: That $\lambda =T_{B}\left( e\right) $ follows
immediately from $e\notin \limfunc{span}(X)$ and the definition of $T_{B}$.
Furthermore, it was shown in Example \ref{EX_quad} that (\ref{ATaylor})
holds with $q=2$ and $D$ given by (\ref{D_2}), and that $D$ satisfies all
conditions required in Theorem \ref{BT2}. Applying the second part of
Theorem \ref{BT2} with $q=2$ then immediately gives (\ref{even_special}).
Furthermore, observe that 
\begin{equation}
\Lambda ^{\prime }\left( C_{X}^{\prime }BC_{X}-\lambda C_{X}^{\prime
}C_{X}\right) \Lambda =A^{\prime }\left( C_{X}^{\prime }BC_{X}-\lambda
C_{X}^{\prime }C_{X}\right) A  \label{matrix}
\end{equation}%
where $A=\Lambda +ee^{\prime }$ is nonsingular (cf. the proof of Theorem \ref%
{BT}). By the general assumptions we have $\lambda <\lambda _{n-k}\left(
B\right) $. If now $\lambda >\lambda _{1}\left( B\right) $ holds, we see
that the matrix in (\ref{matrix}) is not equal to the zero matrix and is
indefinite. Consequently, the r.h.s. of (\ref{even_special}) is strictly
between zero and one. In case $\lambda =\lambda _{1}\left( B\right) $ the
matrix in (\ref{matrix}) is again not equal to the zero matrix, but is now
nonnegative definite, which shows that the r.h.s. of (\ref{even_special})
equals $1$.

Next consider the situation of Part 2 of the corollary: As shown in Example %
\ref{EX_quad}, now condition (\ref{ATaylor}) holds with $q=1$ and $D$ given
by (\ref{D_1}), and $D$ satisfies all conditions required in Theorem \ref%
{BT2}. Applying the second part of Theorem \ref{BT2} now with $q=1$ then
immediately gives (\ref{odd_special}). The claim regarding (\ref{odd_special}%
) falling into $(0,1)$ then follows immediately from Remark \ref{rem_accu}%
(iii), while the final claim follows from this in conjunction with Remark %
\ref{rem_accu}(ii). The claim in parenthesis follows from the second part of
Theorem \ref{BT2} and the following observation: Note that $\Lambda
U_{0}^{\prime }e=0$ implies that%
\begin{equation*}
a_{1}:=\left( e^{\prime }C_{X}^{\prime }BC_{X}-\left\Vert C_{X}e\right\Vert
^{-2}\left( e^{\prime }C_{X}^{\prime }BC_{X}e\right) e^{\prime
}C_{X}^{\prime }C_{X}\right) \Lambda
\end{equation*}%
and $e^{\prime }U_{0}$ are orthogonal. Furthermore, $a_{1}\neq 0$ since the
matrix in parentheses in the definition of $a_{1}$ does not vanish on all of 
$\limfunc{span}(e)^{\bot }$ (see Example \ref{EX_quad}). Since also $%
e^{\prime }U_{0}\neq 0$, we conclude that $a_{1}$ and $e^{\prime }U_{0}$ are
not collinear.

Finally, Part 3 of the corollary follows immediately from Part 3 of Theorem %
\ref{BT2} observing that $q=1$ as shown by Example \ref{EX_quad}. $%
\blacksquare $

\textbf{Proof of Lemma \ref{D2new}:} Let $\kappa $ be a real number such
that $\kappa <T(e)$ and $0<|\kappa -T(e)|<\delta $. Then $e\in \Phi _{\kappa
}$ and $e\notin \limfunc{bd}(\Phi _{\kappa })$ hold, implying that $e\in 
\limfunc{int}(\Phi _{\kappa })$. Theorem \ref{MT} and Remark \ref{RMT} then
entail $\lim_{\rho \rightarrow a}P_{0,1,\rho }(\Phi _{\kappa })=1$. If $%
\kappa <T(e)$ but $|\kappa -T(e)|\geq \delta $ the same conclusion can be
drawn since $\Phi _{\kappa _{1}}\supseteq \Phi _{\kappa _{2}}$ for $\kappa
_{1}\leq \kappa _{2}$. Therefore, we have $\lim_{\rho \rightarrow
a}P_{0,1,\rho }(\Phi _{\kappa })=1$ for every $\kappa <T(e)$. Next, let $%
\kappa $ be a real number such that $\kappa >T(e)$ and $|\kappa
-T(e)|<\delta $ hold. This implies $e\notin \Phi _{\kappa }$ and $e\notin 
\limfunc{bd}(\Phi _{\kappa })$, and hence $e\notin \limfunc{cl}(\Phi
_{\kappa })$. Theorem \ref{MT} and Remark \ref{RMT} now give $\lim_{\rho
\rightarrow a}P_{0,1,\rho }(\Phi _{\kappa })=0$ for those values of $\kappa $%
. Monotonicity of $\Phi _{\kappa }$ w.r.t. $\kappa $ shows that this
relation must hold for all $\kappa >T(e)$. From (\ref{astar}) and the just
established results we obtain%
\begin{equation*}
\alpha ^{\ast }(T)=%
\begin{cases}
\inf_{\kappa <T(e)}P_{0,1,0}(\Phi _{\kappa }) & \text{ if }\liminf_{\rho
\rightarrow a}P_{0,1,\rho }(\Phi _{T(e)})=0, \\ 
\inf_{\kappa \leq T(e)}P_{0,1,0}(\Phi _{\kappa }) & \text{ if }\liminf_{\rho
\rightarrow a}P_{0,1,\rho }(\Phi _{T(e)})>0.%
\end{cases}%
\end{equation*}%
The function $\kappa \mapsto P_{0,1,0}(\Phi _{\kappa })$ is precisely one
minus the cumulative distribution function of $P_{0,1,0}\circ T$, and hence
is continuous at $T\left( e\right) $ by assumption. Since it is clearly also
decreasing in $\kappa $, we may conclude that 
\begin{equation*}
\alpha ^{\ast }(T)=\inf_{\kappa <T(e)}P_{0,1,0}(\Phi _{\kappa
})=\inf_{\kappa \leq T(e)}P_{0,1,0}(\Phi _{\kappa })=P_{0,1,0}(\Phi _{T(e)}).
\end{equation*}%
Finally note that the claim in parenthesis is an immediate consequence of
the second part of Proposition \ref{prop_bound}. $\blacksquare $

\begin{lemma}
\label{DT} Suppose that $Q$ is a probability measure on $\mathbb{R}^{n}$
which is absolutely continuous w.r.t. $\mu _{\mathbb{R}^{n}}$. Let $T_{B}$
be given by (\ref{T_quadratic}).

\begin{enumerate}
\item Then the support of $Q\circ T_{B}$ is contained in $[\lambda
_{1}(B),\lambda _{n-k}(B)]$. Furthermore, if $\lambda _{1}(B)<\lambda
_{n-k}(B)$, the cumulative distribution function of $Q\circ T_{B}$ is
continuous on the real line.

\item If the density of $Q$ is positive on an open neighborhood of the
origin except possibly for a $\mu _{\mathbb{R}^{n}}$-null set, then the
support of $Q\circ T_{B}$ is $[\lambda _{1}(B),\lambda _{n-k}(B)]$.
\end{enumerate}
\end{lemma}

\textbf{Proof:} 1. Observe that the image of $\mathbb{R}^{n-k}\backslash
\left\{ 0\right\} $ under the map $v\mapsto v^{\prime }Bv/v^{\prime }v$ is $%
[\lambda _{1}(B),\lambda _{n-k}(B)]$. Because $T_{B}$ is defined to be $%
\lambda _{1}(B)$ on $\limfunc{span}(X)$, it follows that the range of $T_{B}$
is contained in $[\lambda _{1}(B),\lambda _{n-k}(B)]$, implying that the
support of $Q\circ T_{B}$ is contained in the same interval. [We note for
later use that the range of $T_{B}$ actually coincides with all of $[\lambda
_{1}(B),\lambda _{n-k}(B)]$, because $C_{X}:\mathbb{R}^{n}\rightarrow 
\mathbb{R}^{n-k}$ is surjective.] Next assume that $\lambda _{1}(B)<\lambda
_{n-k}(B)$. To prove the continuity of the cumulative distribution function $%
c\mapsto \left( Q\circ T_{B}\right) \left( (-\infty ,c]\right) $ it suffices
to show that $\left( Q\circ T_{B}\right) \left( \left\{ c\right\} \right) $
is equal to zero for every $c\in \mathbb{R}$. Note that $Q\left( \limfunc{%
span}(X)\right) =0$ since $Q$ is absolutely continuous w.r.t. $\mu _{\mathbb{%
R}^{n}}$ and $k<n$ holds. Consequently, we have for every $c\in \mathbb{R}$%
\begin{equation*}
\left( Q\circ T_{B}\right) \left( \left\{ c\right\} \right) =Q\left( \left\{
y\in \mathbb{R}^{n}:y^{\prime }C_{X}^{\prime }(B-cI_{n-k})C_{X}y=0\right\}
\right) .
\end{equation*}%
To show that $\left( Q\circ T_{B}\right) \left( \left\{ c\right\} \right) =0$
it suffices to show that $\mu _{\mathbb{R}^{n}}\left( \left\{ y\in \mathbb{R}%
^{n}:y^{\prime }C_{X}^{\prime }(B-cI_{n-k})C_{X}y=0\right\} \right) =0$. The
set under consideration is obviously an algebraic set. Hence, it is a $\mu _{%
\mathbb{R}^{n}}$-null set if we can show that the quadratic form in the
definition of this set does not vanish everywhere. Suppose the contrary,
i.e., $y^{\prime }C_{X}^{\prime }(B-cI_{n-k})C_{X}y=0$ for every $y\in 
\mathbb{R}^{n}$ would hold. Because $C_{X}:\mathbb{R}^{n}\rightarrow \mathbb{%
R}^{n-k}$ is surjective, $v^{\prime }(B-cI_{n-k})v=0$ for every $v\in 
\mathbb{R}^{n-k}$ would have to hold. Since $B-cI_{n-k}$ is symmetric, this
would imply $B-cI_{n-k}=0$, contradicting $\lambda _{1}(B)<\lambda _{n-k}(B)$%
. This establishes $\left( Q\circ T_{B}\right) \left( \left\{ c\right\}
\right) =0$.

2. If $\lambda _{1}(B)=\lambda _{n-k}(B)$ this is trivial. Hence assume $%
\lambda _{1}(B)<\lambda _{n-k}(B)$. Let $\lambda $ be an element in the
interior of $[\lambda _{1}(B),\lambda _{n-k}(B)]$ and let $\varepsilon >0$
arbitrary. Without loss of generality assume that $\varepsilon $ is
sufficiently small such that $\left( \lambda -\varepsilon ,\lambda
+\varepsilon \right) \subseteq \lbrack \lambda _{1}(B),\lambda _{n-k}(B)]$.
Let $y\in \mathbb{R}^{n}$ be such that $T_{B}\left( y\right) =\lambda $.
Such an $y$ exists, because the range of $T_{B}$ is all of $[\lambda
_{1}(B),\lambda _{n-k}(B)]$ as noted in the proof of Part 1. But then $%
y\notin \limfunc{span}(X)$ must hold (since $\lambda >\lambda _{1}(B)$), and
hence $T_{B}$ is continuous at $y$. Consequently, there is an open ball that
is mapped into $\left( \lambda -\varepsilon ,\lambda +\varepsilon \right) $
by $T_{B}$. By $G_{X}$-invariance of $T_{B}$ any open neighborhood of the
origin contains such a ball. Because $Q$ has a density that is almost
everywhere positive on a sufficiently small open neighborhood of the origin,
we see that $Q\circ T_{B}$ puts positive mass on $\left( \lambda
-\varepsilon ,\lambda +\varepsilon \right) $. $\blacksquare $

\textbf{Proof of Proposition \ref{D3new}:} 1. Noting that $e\notin \limfunc{%
span}(X)$ and that $T_{B}$ is continuous on $\mathbb{R}^{n}\backslash 
\limfunc{span}$$(X)$, we may use Lemma \ref{D2new} in conjunction with the
preceding Lemma \ref{DT} with $Q=P_{0,1,0}$ to conclude that $\alpha ^{\ast
}(T_{B})=P_{0,1,0}(\Phi _{B,T_{B}(e)})$. Note that this quantity can also be
written as $1-\left( P_{0,1,0}\circ T_{B}\right) \left( (-\infty
,T_{B}(e)]\right) $. Thus $\alpha ^{\ast }\left( T_{B}\right) =0$ is
equivalent to the cumulative distribution function of $T_{B}$ under $%
P_{0,1,0}$ being equal to one when evaluated at $T_{B}(e)$. Lemma \ref{DT}
implies that this is in turn equivalent to $T_{B}(e)=\lambda _{n-k}(B)$
(since $T_{B}\left( e\right) >\lambda _{n-k}(B)$ is clearly impossible). But 
$T_{B}(e)=\lambda _{n-k}(B)$ is clearly equivalent to $C_{X}e\in \limfunc{Eig%
}\left( B,\lambda _{n-k}(B)\right) $. This proves the first claim of Part 1.
Next observe that for every $\kappa \in (-\infty ,\lambda _{n-k}(B))$ the
assumptions on $P_{0,1,0}$ together with Part 2 of Lemma \ref{DT} imply $%
P_{0,1,0}(\Phi _{B,\kappa })>0=\alpha ^{\ast }(T_{B})$. The second claim
then follows from Lemma \ref{D2new}. For the claim in parenthesis see Remark %
\ref{range_for_kappa}.

2. By the same reasoning as in the proof of Part 1 we see that $\alpha
^{\ast }\left( T_{B}\right) =1$ is then equivalent to $T_{B}(e)=\lambda
_{1}(B)$. Since $e\notin \limfunc{span}(X)$ by assumption, this is in turn
equivalent to $C_{X}e\in \limfunc{Eig}\left( B,\lambda _{1}(B)\right) $.
This proves the first claim of Part 2. The second claim follows directly
from Lemma \ref{D2new} because $P_{0,1,0}(\Phi _{B,\kappa })<1=\alpha ^{\ast
}(T_{B})$ holds for $\kappa $ in the specified range in view of Lemma \ref%
{DT} and the assumptions on $P_{0,1,0}$. The remaining claims follow from
Remark \ref{range_for_kappa}.

3. The first claim is obvious in light of Parts 1 and 2, and the remaining
claims follows from Lemma \ref{DT} and Lemma \ref{D2new}. $\blacksquare $

\textbf{Proof of Proposition \ref{D4new}: }The test is obviously invariant
w.r.t. $G_{X}$, and the additional invariance condition (\ref{INV}) in
Theorem \ref{BT} is satisfied because of $e\in \limfunc{span}(X)$. If $%
\kappa \in (\lambda _{1}(B),\lambda _{n-k}(B))$, the rejection region as
well as its complement have positive $\mu _{\mathbb{R}^{n}}$-measure,
whereas $\Phi _{B,\kappa }$ is $\mathbb{R}^{n}$ or the complement of a $\mu
_{\mathbb{R}^{n}}$-null set in case $\kappa \leq \lambda _{1}(B)$, see
Remark \ref{range_for_kappa}. The second claim then follows from Theorem \ref%
{BT}, Part A.3, and Remark \ref{RBT}(i) in case $\kappa \in (\lambda
_{1}(B),\lambda _{n-k}(B))$, and is obvious otherwise. Now, the just
established claim implies that $\alpha ^{\ast }\left( T_{B}\right) $ is not
larger than $\inf \left\{ P_{0,1,0}\left( \Phi _{B,\kappa }\right) :\kappa
<\lambda _{n-k}(B)\right\} =P_{0,1,0}\left( \left\{ T_{B}\geq \lambda
_{n-k}(B)\right\} \right) $. Since the set $\left\{ T_{B}\geq \lambda
_{n-k}(B)\right\} $ is a $\mu _{\mathbb{R}^{n}}$-null set (as $\lambda
_{1}(B)<\lambda _{n-k}(B)$ is assumed) and since $P_{0,1,0}$ is absolutely
continuous by the assumptions of the lemma, $\alpha ^{\ast }\left(
T_{B}\right) =0$ then follows. The claim in parentheses is trivial. $%
\blacksquare $

\begin{lemma}
\label{auxid} Let $A\in \mathbb{R}^{n\times n}$ be a symmetric positive
definite matrix and let $\delta \geq 0$. Then, the following statements are
equivalent:

(i) $C_{X}AC_{X}^{\prime }=\delta I_{n-k}$ for some matrix $C_{X}$
satisfying $C_{X}C_{X}^{\prime }=I_{n-k}$ and $C_{X}^{\prime }C_{X}=\Pi _{%
\limfunc{span}(X)^{\bot }}$,

(ii) $C_{X}AC_{X}^{\prime }=\delta I_{n-k}$ for any matrix $C_{X}$
satisfying $C_{X}C_{X}^{\prime }=I_{n-k}$ and $C_{X}^{\prime }C_{X}=\Pi _{%
\limfunc{span}(X)^{\bot }}$,

(iii) $\Pi _{\limfunc{span}(X)^{\bot }}A\Pi _{\limfunc{span}(X)^{\bot
}}=\delta \Pi _{\limfunc{span}(X)^{\bot }}$,

(iv) there exists a matrix $D$ such that $DD^{\prime }=A$ and $\Pi _{%
\limfunc{span}(X)^{\bot }}D=\delta ^{1/2}\Pi _{\limfunc{span}(X)^{\bot }}$
holds.
\end{lemma}

\textbf{Proof:} That (i), (ii) and (iii) are equivalent is obvious from the
relations $C_{X}C_{X}^{\prime }=I_{n-k}$ and $C_{X}^{\prime }C_{X}=\Pi _{%
\limfunc{span}(X)^{\bot }}$. That (iv) implies (iii) is obvious. To see that
(iii) implies (iv), note that $\Pi _{\limfunc{span}(X)^{\bot }}$ is
symmetric and idempotent, thus%
\begin{equation*}
\Pi _{\text{$\limfunc{span}$}(X)^{\bot }}A^{1/2}A^{1/2}\Pi _{\text{$\limfunc{%
span}$}(X)^{\bot }}=(\delta ^{1/2}\Pi _{\text{$\limfunc{span}$}(X)^{\bot
}})(\delta ^{1/2}\Pi _{\text{$\limfunc{span}$}(X)^{\bot }})^{\prime }.
\end{equation*}%
In other words, $\Pi _{\limfunc{span}(X)^{\bot }}A^{1/2}$ and $\delta
^{1/2}\Pi _{\limfunc{span}(X)^{\bot }}$ are both square roots of the same
matrix, which implies existence of an orthogonal matrix, $U$ say, such that 
\begin{equation*}
\Pi _{\text{$\limfunc{span}$}(X)^{\bot }}A^{1/2}U=\delta ^{1/2}\Pi _{\text{$%
\limfunc{span}$}(X)^{\bot }}.
\end{equation*}%
Setting $D=A^{1/2}U$ then completes the proof. $\blacksquare $

\textbf{Proof of Theorem \ref{ID}:} 1. Clearly, as $\Sigma (\rho ^{\ast })$
is positive definite, we must have $C_{X}\Sigma (\rho ^{\ast })C_{X}^{\prime
}=\delta I_{n-k}$ with $\delta =\delta \left( \rho ^{\ast }\right) >0$. By
Lemma \ref{auxid}, there exists an $n\times n$ matrix $D=D\left( \rho ^{\ast
}\right) $ such that $DD^{\prime }=\Sigma (\rho ^{\ast })$ and $\Pi _{%
\limfunc{span}(X)^{\bot }}D=\delta ^{1/2}\Pi _{\limfunc{span}(X)^{\bot }}$.
Since $D$ is a square root of $\Sigma (\rho ^{\ast })$ there exists an
orthogonal matrix $U(\rho ^{\ast })$ such that $K(\rho ^{\ast })=DU(\rho
^{\ast })$. Now observe that 
\begin{equation}
\Pi _{\text{$\limfunc{span}$}(X)^{\bot }}(X\beta +\sigma K(\rho ^{\ast
})z)=\sigma \Pi _{\text{$\limfunc{span}$}(X)^{\bot }}K(\rho ^{\ast
})z=\sigma \delta ^{1/2}\Pi _{\text{$\limfunc{span}$}(X)^{\bot }}U(\rho
^{\ast })z.  \label{eq1ip}
\end{equation}%
This immediately gives the last equality in (\ref{eq0ip}). Now, if $\Pi _{%
\limfunc{span}(X)^{\bot }}(X\beta +\sigma K(\rho ^{\ast })z)\neq 0$, then we
can use the equation in the previous display to obtain 
\begin{equation*}
\mathcal{I}_{X}(X\beta +\sigma K(\rho ^{\ast })z)=\left\langle \Pi _{\text{$%
\limfunc{span}$}(X)^{\bot }}U(\rho ^{\ast })z/\Vert {\Pi _{\text{$\limfunc{%
span}$}(X)^{\bot }}U(\rho ^{\ast })z}\Vert \right\rangle =\mathcal{I}%
_{X}(U(\rho ^{\ast })z)
\end{equation*}%
and%
\begin{equation*}
\mathcal{I}_{X}^{+}(X\beta +\sigma K(\rho ^{\ast })z)=\Pi _{\text{$\limfunc{%
span}$}(X)^{\bot }}U(\rho ^{\ast })z/\Vert {\Pi _{\text{$\limfunc{span}$}%
(X)^{\bot }}U(\rho ^{\ast })z}\Vert =\mathcal{I}_{X}^{+}(U(\rho ^{\ast })z).
\end{equation*}%
If $\Pi _{\limfunc{span}(X)^{\bot }}(X\beta +\sigma K(\rho ^{\ast })z)=0$,
then also $\Pi _{\limfunc{span}(X)^{\bot }}U(\rho ^{\ast })z=0$ in view of (%
\ref{eq1ip}). Hence, also in this case we obtain $\mathcal{I}_{X}(X\beta
+\sigma K(\rho ^{\ast })z)=\mathcal{I}_{X}(U(\rho ^{\ast })z)$ and $\mathcal{%
I}_{X}^{+}(X\beta +\sigma K(\rho ^{\ast })z)=\mathcal{I}_{X}^{+}(U(\rho
^{\ast })z)$. This proves Part 1.

2. Observe that under the assumption on $\mathfrak{P}$ the distribution of $%
\mathcal{I}_{X}^{1}$ under $P_{\beta ,\sigma ,\rho ^{\ast }}$ is the
distribution of $\mathcal{I}_{X}^{1}(X\beta +\sigma L(\rho ^{\ast })\mathbf{z%
})=\mathcal{I}_{X}^{1}(\sigma \delta ^{1/2}U(\rho ^{\ast })\mathbf{z})$
(upon choosing $K(\rho ^{\ast })=L(\rho ^{\ast })$) which coincides with the
distribution of $\mathcal{I}_{X}^{1}(\sigma \delta ^{1/2}\mathbf{z})$ by the
implied spherical symmetry of the distribution of $\mathbf{z}$. But clearly,
the distribution of $\mathcal{I}_{X}^{1}(\sigma \delta ^{1/2}\mathbf{z})$
coincides with the distribution of $\mathcal{I}_{X}^{1}(\sigma \delta
^{1/2}L\left( 0\right) \mathbf{z})$ by spherical symmetry and since $\Sigma
(0)=I_{n}$ implies that $L\left( 0\right) $ is an orthogonal matrix. In
turn, the distribution of $\mathcal{I}_{X}^{1}(\sigma \delta ^{1/2}L\left(
0\right) \mathbf{z})$ coincides with the distribution of $\mathcal{I}%
_{X}^{1} $ under $P_{0,\sigma \delta ^{1/2},0}$ since $\mathfrak{P}$, in
particular, satisfies Assumption \ref{ASDR}. This proves that $P_{\beta
,\sigma ,\rho ^{\ast }}\circ \mathcal{I}_{X}^{1}=P_{0,\sigma \delta
^{1/2}\left( \rho ^{\ast }\right) ,0}\circ \mathcal{I}_{X}^{1}$. That $%
P_{\beta ,\sigma ,0}\circ \mathcal{I}_{X}^{1}=P_{0,\sigma ,0}\circ \mathcal{I%
}_{X}^{1}$ can be proved in the same way observing that $C_{X}\Sigma
(0)C_{X}^{\prime }=C_{X}C_{X}^{\prime }=I_{n-k}$. [Alternatively, it follows
immediately from $G_{X}^{1}$-invariance and the fact that the distribution
of $\mathbf{z}$ does not depend on $\beta $.] The proofs for the
corresponding statements regarding the distributions of $\mathcal{I}_{X}$
and $\mathcal{I}_{X}^{+}$ are analogous. Since every other invariant
statistic can be represented as a function of $\mathcal{I}_{X}$, $\mathcal{I}%
_{X}^{+}$, and $\mathcal{I}_{X}^{1}$, respectively, the second claim of Part
2 follows. The third claim is now obvious. $\blacksquare $

\section{Proofs for Sections \protect\ref{SEM} and \protect\ref{SPLM}\label%
{app_proofs2}}

\textbf{Proof of Lemma \ref{ISAR1}:} Suppose $\sigma _{1}^{2}\Sigma
_{SEM}(\rho _{1})=\sigma _{2}^{2}\Sigma _{SEM}(\rho _{2})$ and set $\tau
=\sigma _{2}^{2}/\sigma _{1}^{2}$. This implies 
\begin{equation}
(\tau -1)I_{n}=(\tau \rho _{1}-\rho _{2})(W^{\prime }+W)+(\rho _{2}^{2}-\tau
\rho _{1}^{2})W^{\prime }W.  \label{repre}
\end{equation}

If $\tau =1$ inspection of the diagonal elements in (\ref{repre}) shows that
all diagonal elements of $(\rho _{2}^{2}-\rho _{1}^{2})W^{\prime }W$ must be
zero, which is only possible if $\rho _{2}^{2}-\rho _{1}^{2}=0$ since $W$
can not be the zero matrix. But then we arrive at $\rho _{1}=\rho _{2}$ and $%
\sigma _{1}=\sigma _{2}$. Now suppose $\tau \neq 1$ would hold. Then
inspection of the diagonal elements in (\ref{repre}) shows that the diagonal
elements of $W^{\prime }W$ are all identical equal to $b>0$, say, and must
satisfy $\tau -1=(\rho _{2}^{2}-\tau \rho _{1}^{2})b$, which can
equivalently be written as 
\begin{equation*}
\tau \left( 1+\rho _{1}^{2}b\right) =1+\rho _{2}^{2}b
\end{equation*}%
Furthermore, multiplying (\ref{repre}) by $f_{\max }^{\prime }$ from the
left and by $f_{\max }$ from the right and noting that $f_{\max }^{\prime
}f_{\max }=1$ holds, gives after a rearrangement%
\begin{equation*}
\tau \left( 1-\rho _{1}\lambda _{\max }\right) ^{2}=\left( 1-\rho
_{2}\lambda _{\max }\right) ^{2}.
\end{equation*}%
Expressing $\tau $ from the last equation (note that $1-\rho _{1}\lambda
_{\max }>0$), and substituting into the last but one equation gives%
\begin{equation*}
\left( 1+\rho _{1}^{2}b\right) /\left( 1-\rho _{1}\lambda _{\max }\right)
^{2}=\left( 1+\rho _{2}^{2}b\right) /\left( 1-\rho _{2}\lambda _{\max
}\right) ^{2}.
\end{equation*}%
But the function $\rho \mapsto \left( 1+\rho ^{2}b\right) /\left( 1-\rho
\lambda _{\max }\right) ^{2}$ is obviously strictly increasing on $%
[0,\lambda _{\max }^{-1})$ since $b>0$ holds. This gives $\rho _{1}=\rho
_{2} $ and consequently also $\tau =1$ would hold, a contradiction. $%
\blacksquare $

\textbf{Proof of Lemma \ref{SARL}:} Clearly, $\Sigma _{SEM}^{-1}\left(
\left( \lambda _{\max }^{-1}\right) -\right) =(I_{n}-\lambda _{\max
}^{-1}W^{\prime })(I_{n}-\lambda _{\max }^{-1}W)$ and its kernel equals the
kernel of $I_{n}-\lambda _{\max }^{-1}W$ which obviously contains $f_{\max }$
and which is one-dimensional by the assumptions on $W$. Therefore the kernel
equals $\limfunc{span}(f_{\max })$, which together with Lemma \ref{convCM}
proves the first claim. To prove the second claim we need to show that $%
\Lambda $ in the formulation of the lemma is well-defined, is injective when
restricted to $\limfunc{span}(f_{\max })^{\bot }$, and satisfies%
\begin{equation}
\Pi _{\text{$\limfunc{span}$}(f_{\max })^{\bot }}(I_{n}-\rho
W)^{-1}\rightarrow \Lambda  \label{claim_1}
\end{equation}%
for $\rho \rightarrow \lambda _{\max }^{-1}$ with $\rho \in \lbrack
0,\lambda _{\max }^{-1})$. Observe that for every $0\leq \rho <\lambda
_{\max }^{-1}$ we can find a $\delta \left( \rho \right) <1$ such that $\rho
\lambda _{\max }<\delta \left( \rho \right) $ holds. Noting that $\rho
\lambda _{\max }$ is the spectral radius of $\rho W$ by our assumptions on $%
W $, we can conclude that $\left\Vert \left( \rho W\right) ^{j}\right\Vert
^{1/j}\rightarrow \rho \lambda _{\max }<\delta \left( \rho \right) $ for $%
j\rightarrow \infty $ (where $\left\Vert \cdot \right\Vert $ denotes an
arbitrary matrix norm), cf. \cite{HJ1985}, Corollary 5.6.14. But then it
follows that $\left( I_{n}-\rho W\right) ^{-1}$ can be written as the
norm-convergent series $\sum_{j=0}^{\infty }\rho ^{j}W^{j}$ for every $0\leq
\rho <\lambda _{\max }^{-1}$. Thus we obtain 
\begin{equation}
\Pi _{\text{$\limfunc{span}$}(f_{\max })^{\bot }}\left( I_{n}-\rho W\right)
^{-1}=\Pi _{\text{$\limfunc{span}$}(f_{\max })^{\bot }}\sum_{j=0}^{\infty
}\rho ^{j}W^{j}=\Pi _{\text{$\limfunc{span}$}(f_{\max })^{\bot
}}+\sum_{j=1}^{\infty }\rho ^{j}\Pi _{\text{$\limfunc{span}$}(f_{\max
})^{\bot }}W^{j}.  \label{Neumann}
\end{equation}%
Let $g_{2},\ldots ,g_{n}$ be an orthonormal basis of $\limfunc{span}$$%
(f_{\max })^{\bot }$ and define the $n\times (n-1)$ matrix $U_{2}=\left(
g_{2},\ldots ,g_{n}\right) $. Then $U=\left( f_{\max }:U_{2}\right) $ is an
orthogonal matrix. Set $D=U^{\prime }WU$ and observe that $D$ takes the form%
\begin{equation*}
D=%
\begin{pmatrix}
\lambda _{\max } & b^{\prime } \\ 
0 & A%
\end{pmatrix}%
.
\end{equation*}%
For later use we note that $\lambda _{\max }$ is not an eigenvalue of $A$
since the eigenvalues of $D$ and $W$ coincide, since the eigenvalues of $D$
are made up of $\lambda _{\max }$ and the eigenvalues of $A$, and because $%
\lambda _{\max }$ has algebraic multiplicity $1$ by assumption. Now clearly%
\begin{equation}
\Pi _{\limfunc{span}(f_{\max })^{\bot }}W^{j}=\Pi _{\limfunc{span}(f_{\max
})^{\bot }}UD^{j}U^{\prime }=U_{2}A^{j}U_{2}^{\prime }  \label{proj_W_j}
\end{equation}%
holds for $j\geq 1$, which implies%
\begin{equation*}
\Pi _{\text{$\limfunc{span}$}(f_{\max })^{\bot
}}W^{j}=U_{2}A^{j}U_{2}^{\prime }=(U_{2}AU_{2}^{\prime })^{j}=(\Pi _{\text{$%
\limfunc{span}$}(f_{\max })^{\bot }}W)^{j}
\end{equation*}%
for $j\geq 1$. Consequently,%
\begin{eqnarray}
\Pi _{\text{$\limfunc{span}$}(f_{\max })^{\bot }}\left( I_{n}-\rho W\right)
^{-1} &=&\Pi _{\text{$\limfunc{span}$}(f_{\max })^{\bot
}}+\sum_{j=1}^{\infty }\rho ^{j}(\Pi _{\text{$\limfunc{span}$}(f_{\max
})^{\bot }}W)^{j}  \notag \\
&=&-\Pi _{\text{$\limfunc{span}$}(f_{\max })}+\sum_{j=0}^{\infty }\rho
^{j}(\Pi _{\text{$\limfunc{span}$}(f_{\max })^{\bot }}W)^{j}  \notag \\
&=&\left( I_{n}-\rho \Pi _{\text{$\limfunc{span}$}(f_{\max })^{\bot
}}W\right) ^{-1}-\Pi _{\text{$\limfunc{span}$}(f_{\max })},  \label{proj_W_2}
\end{eqnarray}%
observing that the infinite sum in the second line of (\ref{proj_W_2}) is
norm-convergent because of (\ref{Neumann}), and thus necessarily equals the
inverse matrix in the last line of (\ref{proj_W_2}). Because $\lambda _{\max
}$ is not an eigenvalue of $\Pi _{\limfunc{span}(f_{\max })^{\bot }}W$ in
view of (\ref{proj_W_j}) with $j=1$, the matrix $I_{n}-\lambda _{\max
}^{-1}\Pi _{\text{$\limfunc{span}$}(f_{\max })^{\bot }}W$ is invertible,
showing that $\Lambda $ is well-defined. Furthermore, from (\ref{proj_W_2})
we see that (\ref{claim_1}) indeed holds. Finally, $\Lambda $ is injective
on $\limfunc{span}(f_{\max })^{\bot }$ since $\Lambda $ coincides with $%
\left( I_{n}-\lambda _{\max }^{-1}\Pi _{\limfunc{span}(f_{\max })^{\bot
}}W\right) ^{-1}$ on this subspace. $\blacksquare $

\textbf{Proof of Lemma \ref{SARL_2}:} The first claim is an obvious
consequence of the maintained assumptions for the SEM. The second claim
follows from Proposition \ref{AsPexII} together with the already established
first claim, since Assumption \ref{ASC} holds for the SEM as shown in Lemma %
\ref{SARL}. $\blacksquare $

\textbf{Proof of Corollary \ref{Cor1new}:} Parts 1-3 follow from combining
Lemmata \ref{SARL}, \ref{SARL_2}, Theorem \ref{MT}, Remark \ref{RMT},
Theorem \ref{BT}, and Remark \ref{RBT}(i), noting that here $L=L_{\ast }$.
Part 4 is then a simple consequence of Part 3 in view of Proposition \ref%
{prop_bound}, Remark \ref{range_for_kappa}, and Remark \ref{RBT}(iv); cf.
also the proof of Corollary \ref{Lem_Illust_1}. $\blacksquare $

\textbf{Proof of Corollary \ref{Cor2new}:} The first part follows
immediately from Part 1 of Corollary \ref{Lem_Illust_3}. The second part
follows from Part 3 of the same corollary if we can verify that the
additional condition assumed there is satisfied. First observe that $f_{\max
}$ is an eigenvector of $I_{n}-\rho W$ to the eigenvalue $1-\rho \lambda
_{\max }$ and that $I_{n}-\rho W$ is nonsingular for $0\leq \rho <\lambda
_{\max }^{-1}$. Because $I_{n}-\rho W$ is symmetric, $f_{\max }$ is then
also an eigenvector of $\Sigma _{SEM}\left( \rho \right) =\left( I_{n}-\rho
W\right) ^{-2}$ with eigenvalue $\left( 1-\rho \lambda _{\max }\right) ^{-2}$%
. Next observe that $\Pi _{\text{$\limfunc{span}$}(f_{\max })}=f_{\max
}f_{\max }^{\prime }$. But then we have for $0\leq \rho <\lambda _{\max
}^{-1}$%
\begin{eqnarray*}
\Pi _{\text{$\limfunc{span}$}(f_{\max })^{\bot }}\Sigma _{SEM}\left( \rho
\right) \Pi _{\text{$\limfunc{span}$}(f_{\max })} &=&\Pi _{\text{$\limfunc{%
span}$}(f_{\max })^{\bot }}(I_{n}-\rho W)^{-2}f_{\max }f_{\max }^{\prime } \\
&=&\left( 1-\rho \lambda _{\max }\right) ^{-2}\Pi _{\text{$\limfunc{span}$}%
(f_{\max })^{\bot }}f_{\max }f_{\max }^{\prime }=0.
\end{eqnarray*}%
$\blacksquare $

\textbf{Proof of Proposition \ref{prop1new}:} Lemmata \ref{SARL} and \ref%
{SARL_2} show that Assumptions \ref{ASC} and \ref{ASD} are satisfied. In
view of the assumptions of the proposition, $P_{0,1,0}$ is clearly
absolutely continuous w.r.t. $\mu _{\mathbb{R}^{n}}$ with a density that is
positive on an open neighborhood of the origin except possibly for a $\mu _{%
\mathbb{R}^{n}}$-null set, and $e=f_{\max }\notin \limfunc{span}\left(
X\right) $ is trivially satisfied since $k=0$. Obviously, $W+W^{\prime }$ is
not a multiple of the identity matrix. [If it were, inspection of the
diagonal elements shows that $W+W^{\prime }$ would have to be the zero
matrix. However, this is also impossible since $f_{\max }^{\prime }\left(
W+W^{\prime }\right) f_{\max }=2f_{\max }^{\prime }Wf_{\max }=2\lambda
_{\max }>0$.] Also, $-\Sigma _{SEM}^{-1}(\bar{\rho})$ cannot be a multiple
of the identity matrix in view of Lemma \ref{ISAR1}. Hence, in both cases we
have $\lambda _{1}\left( B\right) <\lambda _{n}\left( B\right) $.
Proposition \ref{D3new} and the observation that the rejection probabilities
are monotonically decreasing in $\kappa $ now establishes the first claim of
the proposition. It remains to show that $f_{\max }\in \limfunc{Eig}\left(
B,\lambda _{n}(B)\right) $ is equivalent to $f_{\max }$ being an eigenvector
of $W^{\prime }$ under the additional assumptions on $W$. We may assume that 
$f_{\max }$ is entrywise positive. We argue here similarly as in the proof
of Proposition 1 in \cite{Mart10}. Consider first the case where $B=-\Sigma
_{SEM}^{-1}(\bar{\rho})$. If $f_{\max }\in \limfunc{Eig}\left( B,\lambda
_{n}(B)\right) $ then%
\begin{equation*}
\lambda _{n}(B)f_{\max }=Bf_{\max }=-\left( 1-\bar{\rho}\lambda _{\max
}\right) \left( I-\bar{\rho}W^{\prime }\right) f_{\max }
\end{equation*}%
from which it follows that $f_{\max }$ is an eigenvector of $W^{\prime }$.
Conversely, if $f_{\max }$ is an eigenvector of $W^{\prime }$ then $f_{\max
} $ is easily seen to be also an eigenvector of $B=-\Sigma _{SEM}^{-1}(\bar{%
\rho})$ and hence also of $-B^{-1}$. Now, $-B^{-1}=\Sigma _{SEM}(\bar{\rho})$
is an entrywise positive matrix by a result in \cite{Gant}, p. 69.
Consequently, the eigenspace corresponding to its largest eigenvalue is
one-dimensional and is spanned by a unique normalized and entrywise positive
eigenvector $g$, say. Since $-B^{-1}$ is symmetric and $f_{\max }$ is an
entrywise positive eigenvector of $-B^{-1}$, it must correspond to the
largest eigenvalue of $-B^{-1}$ (because otherwise it would have to be
orthogonal to $g$, which is impossible as $f_{\max }$ and $g$ are both
entrywise nonnegative). Hence, $f_{\max }\in \limfunc{Eig}\left(
-B^{-1},\lambda _{n}\left( -B^{-1}\right) \right) =\limfunc{Eig}\left(
B,\lambda _{n}\left( B\right) \right) $. Next consider the case $%
B=W+W^{\prime }$. As before, $f_{\max }\in \limfunc{Eig}\left( B,\lambda
_{n}(B)\right) $ implies hat $f_{\max }$ is an eigenvector of $W^{\prime }$.
Conversely, $f_{\max }$ being an eigenvector of $W^{\prime }$ implies that $%
f_{\max }$ is an eigenvector of $B$. Since $W+W^{\prime }$ is symmetric,
entrywise nonnegative, and irreducible (since $W$ is so) the same argument
as in the first case can be applied. $\blacksquare $

\textbf{Proof of Proposition \ref{E4new}: }As in the proof of Proposition %
\ref{prop1new} it follows that Assumptions \ref{ASC} and \ref{ASD} are
satisfied and that $P_{0,1,0}$ is absolutely continuous w.r.t. $\mu _{%
\mathbb{R}^{n}}$ with a density that is positive on an open neighborhood of
the origin except possibly for a $\mu _{\mathbb{R}^{n}}$-null set. By
assumption $f_{\max }\notin \limfunc{span}(X)$ holds. Consider first case
(ii): Observe that the eigenspaces of $C_{X}\Sigma _{SEM}(\rho
)C_{X}^{\prime }$ and $C_{X}\lambda _{n}^{-1}\left( \Sigma _{SEM}(\rho
)\right) \Sigma _{SEM}(\rho )C_{X}^{\prime }$ are identical. By Assumption %
\ref{ASC} we have for $\rho \rightarrow \lambda _{\max }^{-1}$, $\rho \in
\lbrack 0,\lambda _{\max }^{-1})$,%
\begin{equation*}
C_{X}\lambda _{n}^{-1}\left( \Sigma _{SEM}(\rho )\right) \Sigma _{SEM}(\rho
)C_{X}^{\prime }\rightarrow C_{X}f_{\max }f_{\max }^{\prime }C_{X}^{\prime }
\end{equation*}%
the limiting matrix being a matrix of rank exactly equal to $1$ since $%
C_{X}f_{\max }\neq 0$ by the assumption $f_{\max }\notin \limfunc{span}(X)$.
Hence, its largest eigenvalue is positive and has algebraic multiplicity $1$%
, while all other eigenvalues are zero. It follows from \cite{Tyler81}, p.
726, Lemma 2.1, that then the eigenspace corresponding to the largest
eigenvalue of $C_{X}\lambda _{n}^{-1}\left( \Sigma _{SEM}(\rho )\right)
\Sigma _{SEM}(\rho )C_{X}^{\prime }$ (and thus the eigenspace corresponding
to the largest eigenvalue of $C_{X}\Sigma _{SEM}(\rho )C_{X}^{\prime }$)
converges to the eigenspace of the limiting matrix corresponding to its
largest eigenvalue (in the sense that the corresponding projection matrices
onto these spaces converge). The latter space is obviously given by $%
\limfunc{span}(C_{X}f_{\max })$. Because the eigenspaces of $C_{X}\Sigma
_{SEM}(\rho )C_{X}^{\prime }$ corresponding to the largest eigenvalue are
independent of $\rho $ by assumption, it follows that these eigenspaces all
coincide with $\limfunc{span}(C_{X}f_{\max })$. Consequently, also $\limfunc{%
Eig}\left( B,\lambda _{n-k}(B)\right) =\limfunc{span}(C_{X}f_{\max })$ holds
for $B=-\left( C_{X}\Sigma _{SEM}(\bar{\rho})C_{X}^{\prime }\right) ^{-1}$.
In particular, $\lambda _{1}\left( B\right) <\lambda _{n-k}\left( B\right) $
follows, as $n-k>1$ has been assumed. The result now follows from the first
part of Proposition \ref{D3new}.

Next consider case (i): By assumption $\lambda _{1}\left( B\right) <\lambda
_{n-k}\left( B\right) $ holds. Hence we may apply the first part of
Proposition \ref{D3new} and it remains to show that $C_{X}f_{\max }$ belongs
to $\limfunc{Eig}\left( B,\lambda _{n-k}(B)\right) $. Now, observe that $%
D\left( \rho \right) :=C_{X}\left( \Sigma _{SEM}(\rho )-I_{n}\right)
C_{X}^{\prime }/\rho \rightarrow B$ for $\rho \rightarrow 0$, $\rho >0$.
Because $C_{X}f_{\max }$ is an eigenvector of $C_{X}\Sigma _{SEM}(\rho
)C_{X}^{\prime }$ corresponding to its largest eigenvalue, $\nu \left( \rho
\right) $ say, as was shown above, it is also an eigenvector of $D\left(
\rho \right) $ corresponding to its largest eigenvalue, namely $\left( \nu
\left( \rho \right) -1\right) /\rho $. Because $D\left( \rho \right)
\rightarrow B$ for $\rho \rightarrow 0$, it follows that $C_{X}f_{\max }$ is
an eigenvector of $B$ corresponding to the limit of $\left( \nu \left( \rho
\right) -1\right) /\rho $, which necessarily then needs to coincide with the
largest eigenvalue of $B$. $\blacksquare $

\textbf{Proof of Theorem \ref{theorem_SPLM}:} Observe that the covariance
matrix of $\mathbf{y}$ under $P_{\beta ,\sigma ,\rho }$ is given by $\sigma
^{2}\Sigma _{SEM}(\rho )$. Now, for $\rho \in \lbrack 0,\lambda _{\max
}^{-1})$ we have 
\begin{equation*}
\lambda _{\max }^{-1/2}\left( \Sigma _{SEM}(\rho )\right) \left( I_{n}-\rho
W\right) ^{-1}=\left[ \lambda _{\max }^{-1}\left( \Sigma _{SEM}(\rho
)\right) \Sigma _{SEM}(\rho )\right] ^{1/2}U(\rho ),
\end{equation*}%
for a suitable orthogonal $n\times n$ matrix $U(\rho )$. From Lemma \ref%
{SARL} we know that $\lambda _{\max }^{-1}\left( \Sigma _{SEM}(\rho )\right)
\Sigma _{SEM}(\rho )$ converges to $f_{\max }f_{\max }^{\prime }$ as $\rho
\rightarrow \lambda _{\max }^{-1}$, $\rho \in \lbrack 0,\lambda _{\max
}^{-1})$. Continuity and uniqueness of the symmetric square root hence gives 
\begin{equation*}
\left[ \lambda _{\max }^{-1}\left( \Sigma _{SEM}(\rho )\right) \Sigma
_{SEM}(\rho )\right] ^{1/2}\rightarrow \left( f_{\max }f_{\max }^{\prime
}\right) ^{1/2}=f_{\max }f_{\max }^{\prime }.
\end{equation*}%
Now, let $\rho _{m}\rightarrow \lambda _{\max }^{-1}$, $\rho _{m}\in \lbrack
0,\lambda _{\max }^{-1})$ be an arbitrary sequence. Then we can always find
a subsequence $m^{\prime }$ such that along this subsequence $U(\rho _{m})$
converges to an orthogonal matrix $U$. Consequently, $\lambda _{\max
}^{-1/2}\left( \Sigma _{SEM}(\rho _{m^{\prime }})\right) \left( I_{n}-\rho
_{m^{\prime }}W\right) ^{-1}$ converges to $f_{\max }f_{\max }^{\prime }U$.
Under $P_{\beta ,\sigma ,\rho _{m^{\prime }}}$ the random vector $\mathbf{y}%
/\left\Vert \mathbf{y}\right\Vert $ clearly has the same distribution as%
\begin{equation*}
\lambda _{\max }^{-1/2}\left( \Sigma _{SEM}(\rho _{m^{\prime }})\right)
\left( I_{n}-\rho _{m^{\prime }}W\right) ^{-1}\left( X\beta +\sigma \mathbf{z%
}\right) /\left\Vert \lambda _{\max }^{-1/2}\left( \Sigma _{SEM}(\rho
_{m^{\prime }})\right) \left( I_{n}-\rho _{m^{\prime }}W\right) ^{-1}\left(
X\beta +\sigma \mathbf{z}\right) \right\Vert
\end{equation*}%
where $\mathbf{z}$ is a fixed random vector distributed according to the
distribution of $\mathbf{\varepsilon }$, which is independent of the
parameters by assumption. Observing that the random variable $f_{\max
}^{\prime }U\left( X\beta +\sigma \mathbf{z}\right) $ is almost surely
nonzero by the assumption on the distribution of $\mathbf{\varepsilon }$,
the expression in the preceding display is now seen to converge in
distribution as $m^{\prime }\rightarrow \infty $ to 
\begin{equation*}
f_{\max }f_{\max }^{\prime }U\left( X\beta +\sigma \mathbf{z}\right)
/\left\Vert f_{\max }f_{\max }^{\prime }U\left( X\beta +\sigma \mathbf{z}%
\right) \right\Vert =\mathbf{c}f_{\max }
\end{equation*}%
where $\mathbf{c}$ is a random variable with values in $\left\{ -1,1\right\} 
$. It then follows from the continuous mapping theorem that $\mathcal{I}%
_{0,\zeta _{f_{\max }}}(\mathbf{y})$ converges in distribution under $%
P_{\beta ,\sigma ,\rho _{m^{\prime }}}$ to $\zeta _{f_{\max }}\left( \mathbf{%
c}f_{\max }\right) =f_{\max }$. In other words, $P_{\beta ,\sigma ,\rho
_{m^{\prime }}}\circ \mathcal{I}_{0,\zeta _{f_{\max }}}$ converges weakly to
pointmass $\delta _{f_{\max }}$. Now observe that%
\begin{eqnarray*}
E_{\beta ,\sigma ,\rho _{m^{\prime }}}(\varphi ) &=&\int_{\mathbb{R}%
^{n}}\varphi (y)dP_{\beta ,\sigma ,\rho _{m^{\prime }}}\left( y\right)
=\int_{\mathbb{R}^{n}}\varphi (\mathcal{I}_{0,\zeta _{f_{\max
}}}(y))dP_{\beta ,\sigma ,\rho _{m^{\prime }}}\left( y\right) \\
&=&\int_{\mathbb{R}^{n}}\varphi (y)d(P_{\beta ,\sigma ,\rho _{m^{\prime
}}}\circ \mathcal{I}_{0,\zeta _{f_{\max }}})\left( y\right) .
\end{eqnarray*}%
But the r.h.s. of the preceding display converges to $\varphi \left( f_{\max
}\right) $ because $P_{\beta ,\sigma ,\rho _{m^{\prime }}}\circ \mathcal{I}%
_{0,\zeta _{f_{\max }}}$ converges weakly to pointmass $\delta _{f_{\max }}$
and because $\varphi $ is bounded and is continuous at $f_{\max }$, cf.
Theorem 30.12 in \cite{Bauer}. A standard subsequence argument then shows
that the limit of $E_{\beta ,\sigma ,\rho }(\varphi )$ for $\rho \rightarrow
\lambda _{\max }^{-1}$, $\rho \in \lbrack 0,\lambda _{\max }^{-1})$ is as
claimed. The second claim is an immediate consequence of the first one. $%
\blacksquare $

\textbf{Proof of Lemma \ref{SEMID}:} From (\ref{eig}) we obtain $W^{\prime
}\Pi _{\text{$\limfunc{span}$}(X)^{\bot }}=\lambda \Pi _{\text{$\limfunc{span%
}$}(X)^{\bot }}$. For $0\leq \rho <\lambda _{\max }^{-1}$ we thus obtain 
\begin{equation*}
\left( I_{n}-\rho W^{\prime }\right) ^{-1}\Pi _{\text{$\limfunc{span}$}%
(X)^{\bot }}=\left( 1-\rho \lambda \right) ^{-1}\Pi _{\text{$\limfunc{span}$}%
(X)^{\bot }},
\end{equation*}%
which after transposition establishes the first claim. An immediate
consequence of the first claim is 
\begin{equation*}
\Pi _{\text{$\limfunc{span}$}(X)^{\bot }}\Sigma _{SEM}(\rho )\Pi _{\text{$%
\limfunc{span}$}(X)^{\bot }}=\left( 1-\rho \lambda \right) ^{-2}\Pi _{\text{$%
\limfunc{span}$}(X)^{\bot }}
\end{equation*}%
which establishes the second claim in view of Lemma \ref{auxid}. $%
\blacksquare $

\section{Auxiliary Results\label{auxil}}

\begin{lemma}
\label{Proj_1}Let $\mathbf{z}$ be a random $n$-vector with a density, $p$
say, w.r.t. $\mu _{\mathbb{R}^{n}}$. Then $\mathbf{s}=\mathbf{z}/\left\Vert 
\mathbf{z}\right\Vert $ is well-defined with probability $1$ and has a
density, $\bar{p}$ say, w.r.t. the uniform probability measure $\upsilon
_{S^{n-1}}$ on $S^{n-1}$. The density $\bar{p}$ satisfies 
\begin{equation*}
\bar{p}\left( s\right) =c\int_{\left( 0,\infty \right) }p\left( rs\right)
r^{n-1}d\mu _{(0,\infty )}(r)
\end{equation*}%
$\upsilon _{S^{n-1}}$-almost everywhere, where $c=2\pi ^{n/2}/\Gamma \left(
n/2\right) $. Furthermore, if $p$ is positive on an open neighborhood of the
origin except possibly for a $\mu _{\mathbb{R}^{n}}$-null set (which is, in
particular, the case if $p$ is positive $\mu _{\mathbb{R}^{n}}$-almost
everywhere), then $\bar{p}$ is positive $\upsilon _{S^{n-1}}$-almost
everywhere.
\end{lemma}

\textbf{Proof: }Let $B$ be a Borel set in $S^{n-1}$ and let $\chi :\mathbb{R}%
^{n}\backslash \left\{ 0\right\} \rightarrow S^{n-1}$ be given by $\chi
\left( z\right) =z/\left\Vert z\right\Vert $. Then 
\begin{equation*}
\Pr (\mathbf{s}\in B)=\Pr (\mathbf{z}\in \chi ^{-1}\left( B\right) )=\int_{%
\mathbb{R}^{n}\backslash \left\{ 0\right\} }\boldsymbol{1}_{\chi ^{-1}\left(
B\right) }\left( z\right) p\left( z\right) dz=\int_{\left( 0,\infty \right)
\times S^{n-1}}\boldsymbol{1}_{\chi ^{-1}\left( B\right) }\left( rs\right)
p\left( rs\right) dH(r,s)
\end{equation*}%
where $H$ is the pushforward measure of $\mu _{\mathbb{R}^{n}}$ (restricted
to $\mathbb{R}^{n}\backslash \left\{ 0\right\} $) under the map $z\mapsto
\left( \left\Vert z\right\Vert ,z/\left\Vert z\right\Vert \right) $. But $H$
is nothing else than the product of the measure on $\left( 0,\infty \right) $
with Lebesgue density $r^{n-1}$ and the surface measure $c\upsilon
_{S^{n-1}} $ on $S^{n-1}$ where $c$ is given in the lemma (cf. \cite{Stroock}%
). In view of Tonelli's theorem (observe all functions involved are
nonnegative) and since $\boldsymbol{1}_{\chi ^{-1}\left( B\right) }\left(
rs\right) =\boldsymbol{1}_{\chi ^{-1}\left( B\right) }\left( s\right) =%
\boldsymbol{1}_{B}\left( s\right) $ clearly holds for $s\in S^{n-1}$, we
obtain%
\begin{equation*}
\Pr (\mathbf{s}\in B)=\int_{S^{n-1}}\boldsymbol{1}_{B}\left( s\right) \left(
c\int_{\left( 0,\infty \right) }p\left( rs\right) r^{n-1}d\mu _{(0,\infty
)}(r)\right) d\upsilon _{S^{n-1}}\left( s\right) ,
\end{equation*}%
which establishes the claims except for the last one. We next prove the
final claim. First, observe that for every Borel set $B$ in $S^{n-1}$ we
have $\upsilon _{S^{n-1}}\left( B\right) >0$ if and only if $\mu _{\mathbb{R}%
^{n}}\left( \chi ^{-1}\left( B\right) \right) >0$. [This is seen as follows:
Specializing what has been proved so far to the case where $\mathbf{z}$
follows a standard Gaussian distribution, shows that in this case $\mathbf{s}
$ is uniformly distributed on $S^{n-1}$. Hence, $\upsilon _{S^{n-1}}\left(
B\right) =\Pr (\mathbf{s}\in B)=\Pr (\mathbf{z}\in \chi ^{-1}\left( B\right)
)$. But then the equivalence of the Gaussian measure with $\mu _{\mathbb{R}%
^{n}}$ establishes that $\upsilon _{S^{n-1}}\left( B\right) >0$ if and only
if $\mu _{\mathbb{R}^{n}}\left( \chi ^{-1}\left( B\right) \right) >0$.] Let
now $B$ satisfy $\upsilon _{S^{n-1}}\left( B\right) >0$. Clearly, $\Pr (%
\mathbf{s}\in B)=\Pr (\mathbf{z}\in \chi ^{-1}\left( B\right) )\geq \Pr (%
\mathbf{z}\in \chi ^{-1}\left( B\right) \cap V)$ where $V$ is an open
neighborhood of the origin on which $p$ is positive $\mu _{\mathbb{R}^{n}}$%
-almost everywhere. But then we must have $\mu _{\mathbb{R}^{n}}\left( \chi
^{-1}\left( B\right) \cap V\right) >0$, because $\mu _{\mathbb{R}^{n}}\left(
\chi ^{-1}\left( B\right) \right) >0$ follows as a consequence of $\upsilon
_{S^{n-1}}\left( B\right) >0$ as just shown above and because $\chi
^{-1}\left( B\right) $ can be written as a countable union of the sets $%
j\left( \chi ^{-1}\left( B\right) \cap V\right) $ with $j\in \mathbb{N}$. By
the assumption on $p$ we can now conclude that $\Pr (\mathbf{z}\in \chi
^{-1}\left( B\right) \cap V)>0$ holds. Hence, we have established that $\Pr (%
\mathbf{s}\in B)>0$ holds whenever $\upsilon _{S^{n-1}}\left( B\right) >0$
is satisfied. $\blacksquare $

\begin{remark}
\label{E1}(i) In the proof we have shown that for $B$, a Borel subset of the
unit sphere, we have $\upsilon _{S^{n-1}}\left( B\right) >0$ if and only if $%
\mu _{\mathbb{R}^{n}}\left( \chi ^{-1}\left( B\right) \right) >0$, a fact
that we shall freely use in various places.

(ii) Let $\mathbf{z}$ be a random $n$-vector such that $\Pr (\mathbf{z}=0)=0$%
. Assume that $\mathbf{z}/\left\Vert \mathbf{z}\right\Vert $ has a density
w.r.t. $\upsilon _{S^{n-1}}$ (which is, in particular, the case if $\mathbf{z%
}$ is spherically symmetric). Let $A$ be a $G_{0}^{+}$-invariant Borel set
in $\mathbb{R}^{n}$ with $\mu _{\mathbb{R}^{n}}\left( A\right) =0$. Then $%
\Pr (\mathbf{z}\in A)=0$ holds. To see this use $G_{0}^{+}$-invariance and
the fact that $\mathbf{z}$ has no atom at the origin to obtain $\Pr (\mathbf{%
z}\in A)=\Pr (\mathbf{z}\in A\backslash \left\{ 0\right\} )=\Pr (\mathbf{z}%
/\left\Vert \mathbf{z}\right\Vert \in A\backslash \left\{ 0\right\} )=\Pr (%
\mathbf{z}/\left\Vert \mathbf{z}\right\Vert \in B)$, where $B=\chi \left(
A\backslash \left\{ 0\right\} \right) $. Note that $B$ is a Borel subset of $%
S^{n-1}$ satisfying $\chi ^{-1}\left( B\right) =A\backslash \left\{
0\right\} $. Hence $\mu _{\mathbb{R}^{n}}\left( \chi ^{-1}\left( B\right)
\right) =0$ holds. But then $\upsilon _{S^{n-1}}\left( B\right) =0$ by what
was shown in (i). Since $\mathbf{s}=\mathbf{z}/\left\Vert \mathbf{z}%
\right\Vert $ possesses a density w.r.t. $\upsilon _{S^{n-1}}$ by
assumption, we conclude that $\Pr (\mathbf{z}/\left\Vert \mathbf{z}%
\right\Vert \in B)=0$, and thus also $\Pr (\mathbf{z}\in A)=0$ must hold.

(iii) Let $\mathbf{z}$ be as in (ii) and let $A$ be a $G_{X}^{+}$-invariant
Borel set in $\mathbb{R}^{n}$ with $\mu _{\mathbb{R}^{n}}\left( A\right) =0$%
. Then for every $\beta \in \mathbb{R}^{n}$, $0<\sigma <\infty $, and every
nonsingular $n\times n$ matrix $L$ we have $\Pr (X\beta +\sigma L\mathbf{z}%
\in A)=\Pr (L\mathbf{z}\in A)=\Pr (\mathbf{z}\in L^{-1}\left( A\right) )=0$
in view of (ii) since $L^{-1}\left( A\right) $ is a $G_{0}^{+}$-invariant $%
\mu _{\mathbb{R}^{n}}$-null set.
\end{remark}

\begin{lemma}
\label{Proj_2}Let $\mathbf{z}$ be a random $n$-vector satisfying $\Pr (%
\mathbf{z}=0)=0$. Then $\mathbf{s}=\mathbf{z}/\left\Vert \mathbf{z}%
\right\Vert $ is well-defined with probability $1$. Assume further that the
distribution of $\mathbf{s}$ has a density, $g$ say, w.r.t. $\upsilon
_{S^{n-1}}$. Suppose $\mathbf{r}$ is a random variable taking values in $%
(0,\infty )$ that is independent of $\mathbf{z}/\left\Vert \mathbf{z}%
\right\Vert $ and that has a density, $h$ say, w.r.t. $\mu _{(0,\infty )}$.
Define $\mathbf{z}^{\dag }=\mathbf{rz}/\left\Vert \mathbf{z}\right\Vert $ on
the event $\mathbf{z\neq 0}$ and assign arbitrary values to $\mathbf{z}%
^{\dag }$ on the event $\mathbf{z=0}$ in a measurable way. Then, the
following holds:

\begin{enumerate}
\item $\Pr (\mathbf{z}^{\dag }=0)=0$ and $\mathbf{z}^{\dag }/\Vert {\mathbf{z%
}}^{\dag }\Vert =\mathbf{z}/\Vert {\mathbf{z}}\Vert $ for $\mathbf{z}\neq 0$%
, $\mathbf{z}^{\dag }\neq 0$.

\item $\mathbf{z}^{\dag }$ possesses a density $g^{\dag }$ w.r.t. Lebesgue
measure $\mu _{\mathbb{R}^{n}}$ which is given by%
\begin{equation*}
g^{\dag }(z)=%
\begin{cases}
c^{-1}g\left( z/\Vert {z}\Vert \right) \frac{h\left( \Vert {z}\Vert \right) 
}{\Vert {z}\Vert ^{n-1}} & \text{ if }z\neq 0 \\ 
0 & \text{ if }z=0,%
\end{cases}%
\end{equation*}%
where $c$ has been given in Lemma \ref{Proj_1}.

\item If $g$ is $\upsilon _{S^{n-1}}$-almost everywhere continuous and $h$
is $\mu _{(0,\infty )}$-almost everywhere continuous, then $g^{\dag }$ is $%
\mu _{\mathbb{R}^{n}}$-almost everywhere continuous.

\item If $g$ is $\upsilon _{S^{n-1}}$-almost everywhere positive and $h$ is $%
\mu _{(0,\infty )}$-almost everywhere positive, then $g^{\dag }$ is $\mu _{%
\mathbb{R}^{n}}$-almost everywhere positive.

\item If $g$ is constant $\upsilon _{S^{n-1}}$-almost everywhere [which is,
in particular, the case if $\mathbf{z}$ is spherically symmetric] and if $%
\mathbf{r}$ is distributed as the square root of a $\chi ^{2}$-distributed
random variable with $n$ degrees of freedom, then $\mathbf{z}^{\dag }$ is
Gaussian with mean zero and covariance matrix $I_{n}$.
\end{enumerate}
\end{lemma}

\textbf{Proof:} Part 1 is obvious. To prove Part 2 we denote the
distribution of $\mathbf{z}/\Vert {\mathbf{z}}\Vert $ by $G$ and the
distribution of $\mathbf{r}$ by $H$. Because $\mathbf{z}/\left\Vert \mathbf{z%
}\right\Vert $ and $\mathbf{r}$ are independent, the joint distribution of $%
\mathbf{z}/\Vert {\mathbf{z}}\Vert $ and $\mathbf{r}$ on $S^{n-1}\times
(0,\infty )$, equipped with the product $\sigma $-field, is given by the
product measure $G\otimes H$. Therefore, the distribution of $\mathbf{z}%
^{\dag }$ is the push-forward measure of $G\otimes H$ under the mapping $%
m(s,r)=rs$. Hence for every $A\in \mathcal{B}(\mathbb{R}^{n})$ we have,
using Tonelli's theorem and the fact that $G$ and $H$ have densities $g$ and 
$h$, respectively, that%
\begin{align*}
\Pr (\mathbf{z}^{\dag }& \in A)=\int_{S^{n-1}\times (0,\infty )}\mathbf{1}%
_{A}(rs)d(G\otimes H)(s,r)=\int_{(0,\infty )}\int_{S^{n-1}}\mathbf{1}%
_{A}(rs)dG(s)dH(r) \\
& =\int_{(0,\infty )}\int_{S^{n-1}}\mathbf{1}_{A}(rs)g(s)d\upsilon
_{S^{n-1}}(s)h(r)d\mu _{(0,\infty )}(r) \\
& =\int_{(0,\infty )}r^{n-1}\int_{S^{n-1}}\mathbf{1}_{A}(rs)g(s)r^{1-n}h%
\left( r\right) d\upsilon _{S^{n-1}}(s)d\mu _{(0,\infty )}(r) \\
& =\int_{(0,\infty )}r^{n-1}\int_{S^{n-1}}f\left( rs\right) d\upsilon
_{S^{n-1}}(s)d\mu _{(0,\infty )}(r),
\end{align*}%
where for $x\in \mathbb{R}^{n}$ the function $f$ is given by%
\begin{equation*}
f(x)=%
\begin{cases}
\mathbf{1}_{A}(x)g(x/\Vert {x}\Vert )\Vert {x}\Vert ^{1-n}h\left( \Vert {x}%
\Vert \right) & \text{ if }x\neq 0 \\ 
0 & \text{ if }x=0.%
\end{cases}%
\end{equation*}%
Since $f$ is clearly a non-negative and Borel-measurable function, we can
apply Theorem 5.2.2 in Stroock (1999) to see that 
\begin{align*}
\Pr (\mathbf{z}^{\dag }& \in A)=\int_{(0,\infty
)}r^{n-1}\int_{S^{n-1}}f(rs)d\upsilon _{S^{n-1}}(s)d\mu _{(0,\infty )}(r) \\
& =\int_{\mathbb{R}^{n}}c^{-1}f(x)d\mu _{\mathbb{R}^{n}}(x)=\int_{\mathbb{R}%
^{n}}\mathbf{1}_{A}(x)g^{\dag }(x)d\mu _{\mathbb{R}^{n}}(x).
\end{align*}%
This establishes the second part of the lemma. To prove the third part
denote by $D_{g^{\dag }}\subseteq \mathbb{R}^{n}$, $D_{g}\subseteq S^{n-1}$
and $D_{h}\subseteq (0,\infty )$ the discontinuity points of $g^{\dag }$, $g$%
, and $h$, respectively, which are measurable. Using Part 2 of the lemma we
see that $x\neq 0$, $x/\Vert {x}\Vert \in \mathbb{R}^{n}\backslash D_{g}$,
and $\Vert {x}\Vert \in \mathbb{R}^{n}\backslash D_{h}$ imply $x\in \mathbb{R%
}^{n}\backslash D_{g^{\dag }}$. Therefore, negating the statement, we see
that $\mathbf{1}_{D_{g^{\dag }}}(x)\leq \mathbf{1}_{\left\{ 0\right\} }(x)+%
\mathbf{1}_{D_{g}}(x/\Vert {x}\Vert )+\mathbf{1}_{D_{h}}(\Vert {x}\Vert )$
must hold which implies 
\begin{equation}
\mu _{\mathbb{R}^{n}}(D_{g^{\dag }})=\int_{\mathbb{R}^{n}}\mathbf{1}%
_{D_{g^{\dag }}}(x)d\mu _{\mathbb{R}^{n}}(x)\leq \int_{\mathbb{R}^{n}}%
\mathbf{1}_{D_{g}}(x/\Vert {x}\Vert )d\mu _{\mathbb{R}^{n}}(x)+\int_{\mathbb{%
R}^{n}}\mathbf{1}_{D_{h}}(\Vert {x}\Vert )d\mu _{\mathbb{R}^{n}}(x).
\label{DG*}
\end{equation}%
Using again Theorem 5.2.2 in Stroock (1999) we see that 
\begin{align*}
\int_{\mathbb{R}^{n}}\mathbf{1}_{D_{g}}(x/\Vert {x}\Vert )d\mu _{\mathbb{R}%
^{n}}(x)& =\int_{(0,\infty )}r^{n-1}\int_{S^{n-1}}\mathbf{1}%
_{D_{g}}(s)cd\upsilon _{S^{n-1}}(s)d\mu _{(0,\infty )}(r) \\
& =\int_{(0,\infty )}c\upsilon _{S^{n-1}}(D_{g})r^{n-1}d\mu _{(0,\infty
)}(r)=0,
\end{align*}%
because $\upsilon _{S^{n-1}}(D_{g})=0$ holds by assumption. Similarly, we
obtain 
\begin{equation*}
\int_{\mathbb{R}^{n}}\mathbf{1}_{D_{h}}(\Vert {x}\Vert )d\mu _{\mathbb{R}%
^{n}}(x)=\int_{S^{n-1}}\int_{(0,\infty )}r^{n-1}\mathbf{1}_{D_{h}}(r)d\mu
_{(0,\infty )}(r)cd\upsilon _{S^{n-1}}(s)=0,
\end{equation*}%
because the inner integral is zero as a consequence of the assumption that $%
\mu _{(0,\infty )}(D_{h})=0$. Together with Equation (\ref{DG*}) the last
two displays establish $\mu _{\mathbb{R}^{n}}(D_{g^{\dag }})=0$. To prove
Part 4 denote by $Z_{g^{\dag }}\subseteq \mathbb{R}^{n}$, $Z_{g}\subseteq
S^{n-1}$, and $Z_{h}\subseteq (0,\infty )$ the zero sets of $g^{\dag }$, $g$%
, and $h$, respectively, which are obviously measurable. Replacing $%
D_{g^{\dag }}$, $D_{g}$, and $D_{h}$ with $Z_{g^{\dag }}$, $Z_{g}$, and $%
Z_{h}$, respectively, in the argument used above then establishes Part 4. To
prove the last part, we observe that $g$ being constant $\upsilon _{S^{n-1}}$%
- almost everywhere implies that $\mathbf{z}/\Vert {\mathbf{z}}\Vert $ is
uniformly distributed on $S^{n-1}$. Since $\mathbf{z}/\Vert {\mathbf{z}}%
\Vert $ is independent of $\mathbf{r}$, which is distributed as the square
root of a $\chi ^{2}$ with $n$ degrees of freedom, it is now obvious that $%
\mathbf{z}^{\dag }$ is Gaussian with mean zero and covariance matrix $I_{n}$%
. $\blacksquare $

\begin{remark}
As long as we are only concerned with distributional properties of $\mathbf{z%
}$ we can assume w.l.o.g. that the probability space supporting $\mathbf{z}$
is rich enough to allow independent random variables $\mathbf{r}$ that have
the required properties. In particular, we can then always choose $\mathbf{r}
$ such that the density is simultaneously $\mu _{(0,\infty )}$-almost
everywhere continuous and $\mu _{(0,\infty )}$-almost everywhere positive
(e.g., by choosing $\mathbf{r}$ to follow a $\chi ^{2}$-distribution).
\end{remark}

{\small {
\bibliographystyle{ims}
\bibliography{refs}

\begin{thebibliography}{30}
\expandafter\ifx\csname natexlab\endcsname\relax\def\natexlab#1{#1}\fi
\expandafter\ifx\csname url\endcsname\relax
  \def\url#1{\texttt{#1}}\fi
\expandafter\ifx\csname urlprefix\endcsname\relax\def\urlprefix{URL }\fi
\providecommand{\eprint}[2][]{\url{#2}}

\bibitem[{Anselin(2001)}]{Anselin}
\textsc{Anselin, L.} (2001).
\newblock Spatial econometrics.
\newblock In \textit{A companion to theoretical econometrics}. Blackwell
  Companions Contemp. Econ., Blackwell, Malden, MA, 310--330.

\bibitem[{Arnold(1979)}]{Arnold79}
\textsc{Arnold, S.~F.} (1979).
\newblock Linear models with exchangeably distributed errors.
\newblock \textit{Journal of the American Statistical Association}, \textbf{74}
  194--199.

\bibitem[{Bartels(1992)}]{bartels1992}
\textsc{Bartels, R.} (1992).
\newblock On the power function of the {D}urbin-{W}atson test.
\newblock \textit{Journal of Econometrics}, \textbf{51} 101--112.

\bibitem[{Bauer(2001)}]{Bauer}
\textsc{Bauer, H.} (2001).
\newblock \textit{Measure and integration theory}, vol.~26 of \textit{de
  Gruyter Studies in Mathematics}.
\newblock Walter de Gruyter \& Co., Berlin.

\bibitem[{Cambanis et~al.(1981)Cambanis, Huang and Simons}]{Cambanis}
\textsc{Cambanis, S.}, \textsc{Huang, S.} and \textsc{Simons, G.} (1981).
\newblock On the theory of elliptically contoured distributions.
\newblock \textit{Journal of Multivariate Analysis}, \textbf{11} 368--385.

\bibitem[{Gantmacher(1959)}]{Gant}
\textsc{Gantmacher, F.~R.} (1959).
\newblock \textit{The theory of matrices, vol 2}.
\newblock AMS Chelsea Publishing, Providence, RI.

\bibitem[{Horn and Johnson(1985)}]{HJ1985}
\textsc{Horn, R.~A.} and \textsc{Johnson, C.~R.} (1985).
\newblock \textit{Matrix analysis}.
\newblock Cambridge University Press, Cambridge.

\bibitem[{Kadiyala(1970)}]{Kadiyala}
\textsc{Kadiyala, K.~R.} (1970).
\newblock Testing for the independence of regression disturbances.
\newblock \textit{Econometrica}, \textbf{38} 97--117.

\bibitem[{Kariya(1980)}]{Kariya80JASA}
\textsc{Kariya, T.} (1980).
\newblock Note on a condition for equality of sample variances in a linear
  model.
\newblock \textit{Journal of the American Statistical Association}, \textbf{75}
  701--703.

\bibitem[{King(1985)}]{King1985}
\textsc{King, M.~L.} (1985).
\newblock A point optimal test for autoregressive disturbances.
\newblock \textit{Journal of Econometrics}, \textbf{27} 21--37.

\bibitem[{King(1987)}]{King1987a}
\textsc{King, M.~L.} (1987).
\newblock Testing for autocorrelation in linear regression models: a survey.
\newblock In \textit{Specification analysis in the linear model}. Internat.
  Lib. Econom., Routledge \& Kegan Paul, London, 19--73.

\bibitem[{King and Hillier(1985)}]{KingHillier1985}
\textsc{King, M.~L.} and \textsc{Hillier, G.~H.} (1985).
\newblock Locally best invariant tests of the error covariance matrix of the
  linear regression model.
\newblock \textit{Journal of the Royal Statistical Society. Series B
  (Methodological)}, \textbf{47} 98--102.

\bibitem[{Kleiber and Kr{\"a}mer(2005)}]{kleiber2005}
\textsc{Kleiber, C.} and \textsc{Kr{\"a}mer, W.} (2005).
\newblock Finite-sample power of the {D}urbin-{W}atson test against
  fractionally integrated disturbances.
\newblock \textit{Econometrics Journal}, \textbf{8} 406--417.

\bibitem[{Kr{\"a}mer(1985)}]{kramer1985}
\textsc{Kr{\"a}mer, W.} (1985).
\newblock The power of the {D}urbin-{W}atson test for regressions without an
  intercept.
\newblock \textit{Journal of Econometrics}, \textbf{28} 363--370.

\bibitem[{Kr{\"a}mer(2005)}]{kramer2005}
\textsc{Kr{\"a}mer, W.} (2005).
\newblock Finite sample power of {C}liff-{O}rd-type tests for spatial
  disturbance correlation in linear regression.
\newblock \textit{Journal of Statistical Planning and Inference}, \textbf{128}
  489--496.

\bibitem[{Kr{\"a}mer and Zeisel(1990)}]{KramerZeisel1990}
\textsc{Kr{\"a}mer, W.} and \textsc{Zeisel, H.} (1990).
\newblock Finite sample power of linear regression autocorrelation tests.
\newblock \textit{Journal of Econometrics}, \textbf{43} 363--372.

\bibitem[{Lehmann and Romano(2005)}]{LR05}
\textsc{Lehmann, E.~L.} and \textsc{Romano, J.~P.} (2005).
\newblock \textit{Testing statistical hypotheses}.
\newblock 3rd ed. Springer Texts in Statistics, Springer, New York.

\bibitem[{L{\"o}bus and Ritter(2000)}]{lobus2000}
\textsc{L{\"o}bus, J.-U.} and \textsc{Ritter, L.} (2000).
\newblock The limiting power of the {D}urbin-{W}atson test.
\newblock \textit{Communications in Statistics-Theory and Methods}, \textbf{29}
  2665--2676.

\bibitem[{Martellosio(2010)}]{Mart10}
\textsc{Martellosio, F.} (2010).
\newblock Power properties of invariant tests for spatial autocorrelation in
  linear regression.
\newblock \textit{Econometric Theory}, \textbf{26} 152--186.

\bibitem[{Martellosio(2011{\natexlab{a}})}]{Mart11SPL}
\textsc{Martellosio, F.} (2011{\natexlab{a}}).
\newblock Efficiency of the {OLS} estimator in the vicinity of a spatial unit
  root.
\newblock \textit{Statistics \& Probability Letters}, \textbf{81} 1285 -- 1291.

\bibitem[{Martellosio(2011{\natexlab{b}})}]{Mart11}
\textsc{Martellosio, F.} (2011{\natexlab{b}}).
\newblock Nontestability of equal weights spatial dependence.
\newblock \textit{Econometric Theory}, \textbf{27} 1369--1375.

\bibitem[{Martellosio(2012)}]{Mart12}
\textsc{Martellosio, F.} (2012).
\newblock Testing for spatial autocorrelation: The regressors that make the
  power disappear.
\newblock \textit{Econometric Reviews}, \textbf{31} 215--240.

\bibitem[{Mynbaev(2012)}]{mynbaev2012}
\textsc{Mynbaev, K.} (2012).
\newblock Distributions escaping to infinity and the limiting power of the
  {C}liff-{O}rd test for autocorrelation.
\newblock \textit{ISRN Probability and Statistics}, \textbf{2012}.

\bibitem[{Preinerstorfer(2014)}]{Prein2014}
\textsc{Preinerstorfer, D.} (2014).
\newblock How to avoid the zero-power trap in testing for correlation.
\newblock \textit{Working paper, in preparation.}

\bibitem[{Preinerstorfer and P{\"o}tscher(2013)}]{PP13}
\textsc{Preinerstorfer, D.} and \textsc{P{\"o}tscher, B.~M.} (2013).
\newblock On size and power of heteroskedasticity and autocorrelation robust
  tests.
\newblock \textit{Econometric Theory, forthcoming.}

\bibitem[{Small(1993)}]{small1993}
\textsc{Small, J.~P.} (1993).
\newblock The limiting power of point optimal autocorrelation tests.
\newblock \textit{Communications in Statistics--Theory and Methods},
  \textbf{22} 3907--3916.

\bibitem[{Stroock(1999)}]{Stroock}
\textsc{Stroock, D.~W.} (1999).
\newblock \textit{A concise introduction to the theory of integration}.
\newblock 3rd ed. Birkh\"auser Boston Inc., Boston, MA.

\bibitem[{Tillman(1975)}]{Tillman1975}
\textsc{Tillman, J.} (1975).
\newblock The power of the {D}urbin-{W}atson test.
\newblock \textit{Econometrica}, \textbf{43} 959--974.

\bibitem[{Tyler(1981)}]{Tyler81}
\textsc{Tyler, D.~E.} (1981).
\newblock Asymptotic inference for eigenvectors.
\newblock \textit{Annals of Statistics}, \textbf{9} pp. 725--736.

\bibitem[{Zeisel(1989)}]{zeisel1989}
\textsc{Zeisel, H.} (1989).
\newblock On the power of the {D}urbin-{W}atson test under high
  autocorrelation.
\newblock \textit{Communications in Statistics-Theory and Methods}, \textbf{18}
  3907--3916.

\end{thebibliography}
}}

\end{document}